\documentclass[10pt]{article}
\usepackage{mathrsfs}
\usepackage{bbm}
\usepackage[leqno]{amsmath}
\usepackage{amsfonts}
\usepackage{graphicx}
\usepackage{amsmath}
\usepackage{amssymb}
\usepackage{latexsym}
\usepackage{amsmath, amsfonts,amssymb, amsthm, euscript,makeidx,color,mathrsfs}
\usepackage{indentfirst}

\usepackage{appendix}

\usepackage{stmaryrd}

\def\sqr#1#2{{\vcenter{\vbox{\hrule height.#2pt
				\hbox{\vrule width.#2pt height#1pt \kern#1pt \vrule width.#2pt}
				\hrule height.#2pt}}}}

\def\3n{\negthinspace \negthinspace \negthinspace }
\def\2n{\negthinspace \negthinspace }
\def\1n{\negthinspace }

\def\dbE{\mathbb{E}}
\def\dbF{\mathbb{F}}

\def\dbP{\mathbb{P}}

\def\dbR{\mathbb{R}}

\def\={\buildrel \triangle \over =}

\def\ds{\displaystyle}

\def\e{\varepsilon}

\def\l{\lambda}

\def\i{\infty}
\def\cD{{\cal D}}

\def\cS{{\cal S}}

\def\no{\noindent}

\def\ss{\smallskip}
\def\ms{\medskip}

\def\q{\quad}
\def\qq{\qquad}

\def\limsup{\mathop{\overline{\rm lim}}}
\def\liminf{\mathop{\underline{\rm lim}}}

\def\esssup{\mathop{\rm esssup}}
\def\essinf{\mathop{\rm essinf}}

\def\cd{\cdot}

\def\les{\leqslant}
\def\ges{\geqslant}

\def\({\Big (}
\def\){\Big )}
\def\[{\Big[}
\def\]{\Big]}
\def\bde{\begin{definition}}
	\def\ede{\end{definition}}
\def\be{\begin{equation}}
\def\bel{\begin{equation}\label}
\def\ee{\end{equation}}
\def\bt{\begin{theorem}}
	\def\et{\end{theorem}}
\def\bc{\begin{corollary}}
	\def\ec{\end{corollary}}
\def\bl{\begin{lemma}}
	\def\el{\end{lemma}}
\def\bp{\begin{proposition}}
	\def\ep{\end{proposition}}
\def\bas{\begin{assumption}}
	\def\eas{\end{assumption}}
\def\br{\begin{remark}}
	\def\er{\end{remark}}
\def\ba{\begin{array}}
	\def\ea{\end{array}}
\def\ed{\end{document}}

\def\square#1{\vbox{\hrule\hbox{\vrule height#1%
		\kern#1\vrule}\hrule}}
\def\rectangle#1#2{\vbox{\hrule\hbox{\vrule height#1%
		\kern#2\vrule}\hrule}}

\font\tenbb=msbm10 \font\sevenbb=msbm7 \font\fivebb=msbm5

\newfam\bbfam
\scriptscriptfont\bbfam=\fivebb \textfont\bbfam=\tenbb
\scriptfont\bbfam=\sevenbb

\newtheorem{lemma}{Lemma}[section]
\newtheorem{remark}{Remark}[section]

\newtheorem{theorem}{Theorem}[section]
\newtheorem{corollary}{Corollary}[section]

\newtheorem{definition}{Definition}[section]
\newtheorem{proposition}{Proposition}[section]
\newtheorem{assumption}{Assumption}[section]

\makeatletter

\@addtoreset{equation}{section}
\makeatother

\numberwithin{equation}{section}
\title{Infinite-Horizon Stochastic Recursive Differential Games
with Jumps: Nonautonomous Coefficients and Nonzero Terminal Values\footnote{This work is supported by  the National Key R\&D Program of China (No. 2023YFA1009002),  the
National Natural Science Foundation of China (No. 12371443) and the Changbai
Talent Program of Jilin Province.}}

\author{Sheng Luo\thanks{School of Computer Science and Technology (School of Software),
Harbin Institute of Technology, Weihai, China
(Email: {\tt 15263105996@163.com}).} \and
  Qingmeng Wei\thanks{School of Mathematics and Statistics, Northeast Normal University, Changchun 130024, China
          (Email: {\tt weiqm100@nenu.edu.cn}).}
            \and
           Yaxin Wen%
   \thanks{{School of Mathematics and Statistics, Northeast Normal University, Changchun 130024, China  (Email: {\tt wenyx110@163.com}).}}
}

\date{}
\begin{document}
\maketitle
\vspace{-2em}
\begin{abstract}
In this paper, we investigate a class of decoupled infinite-horizon forward--backward stochastic differential equations (FBSDEs) with time-dependent coefficients and jumps. Under suitable time-integrability conditions on the coefficient growth and stochastic Lipschitz moduli, we establish an $L^p$ theory ($p\ges 2$) for these systems, including well-posedness, moment estimates, conditional stability, and tail estimates, without imposing dissipativity of the forward dynamics or monotonicity of the backward generator in its value variable. In particular, the forward state converges in $L^p$ and almost surely to a limiting random variable $X_\infty$, while the backward solution converges to the generally nonzero terminal value $\Phi(X_\infty)$.
Within this framework, we study infinite-horizon two-player zero-sum stochastic recursive differential games. We prove that the lower and upper value functions are deterministic and establish the dynamic programming principle (DPP) using backward stochastic semigroups and nonanticipative strategies. Both value functions converge locally uniformly to $\Phi$ as the initial time tends to infinity and are characterized as viscosity solutions of the corresponding non-autonomous, nonlocal Hamilton--Jacobi--Bellman--Isaacs (HJBI) equations. Under an additional integrable Lyapunov-barrier condition, we establish a comparison principle for continuous viscosity subsolutions and supersolutions of uniform linear growth, provided that at least one is uniformly Lipschitz in space. This yields uniqueness for both HJBI equations in the linear-growth class. Finally, we establish a weaker, contact-wise Isaacs criterion ensuring that the lower and upper value functions coincide and hence that the game has a value.
\end{abstract}

\noindent\textbf{Keywords:}
infinite-horizon decoupled FBSDEs with jumps;
stochastic recursive differential games;
time-integrable Lipschitz {conditions};
nonzero terminal values;
dynamic programming principle;
nonlocal HJBI equations.

\noindent\textbf{MSC 2020:}  91A15, 49N70, 60H10, 49L25.

 \section{Introduction}

Since the pioneering work of Pardoux and Peng \cite{PardouxPeng1990}, backward stochastic differential equations (BSDEs) have become an important tool in stochastic control and differential games, as well as in the probabilistic study of nonlinear partial differential equations.
Poisson random measures provide a natural framework for describing discontinuities arising from random jumps. In this setting, Tang and Li \cite{TangLi1994} introduced BSDEs driven by a Brownian motion and an independent compensated Poisson random measure. These equations, commonly referred to as BSDEs with jumps, were subsequently used by Barles, Buckdahn and Pardoux \cite{BarlesBuckdahnPardoux1997} to establish probabilistic representations for nonlinear partial integro-differential equations. Further developments of BSDE-based dynamic programming methods for stochastic control and differential games with jumps can be found in \cite{LiPeng2009,BuckdahnHuLi2011,LiWei2015,LuoLiLiWei2025,LiuWei2026} and the references therein.

For stochastic differential equations (SDEs) and BSDEs on an infinite horizon, the analysis requires establishing global well-posedness and stability, understanding the long-time behavior of the forward state, and specifying an appropriate terminal condition for the backward equation at infinity. These requirements make the analysis more involved than in the finite-horizon setting.
  Monotonicity, dissipativity, discounting, and related structural conditions have been widely used to establish well-posedness and stability for BSDEs and forward--backward stochastic systems over unbounded horizons. Early results on BSDEs with random terminal times, possibly taking infinite values, can be traced back to Peng \cite{Peng1991}.  {Subsequent studies of BSDEs and forward--backward systems with random terminal times or infinite horizons} can be found in \cite{DarlingPardoux1997,BriandHu1998,Pardoux1998,HuTessitore2007,PengShi2000,WeiYu2021,HW2026} and the references therein. In the presence of jumps, related results under monotonicity or dissipativity conditions can be found in \cite{Royer2006, Yu2017,LuoLiWei2026}.

Another line of research concerns infinite-horizon BSDEs under time-integrable Lipschitz conditions, without requiring monotonicity of the generator in its backward variable. Chen and Wang \cite{ChenWang2000} established well-posedness for infinite-horizon BSDEs with general square-integrable terminal values under deterministic time-dependent Lipschitz conditions. Liu, Li and Fan \cite{LiuLiFan2020} further established $L^p$ well-posedness under stochastic Lipschitz conditions, which was subsequently used by Li and Yu \cite{LiYu2025} in their study of infinite-horizon FBSDEs. In the presence of jumps, Wu \cite{Wu2003} investigated BSDEs with possibly infinite random terminal times under time-integrable Lipschitz conditions. Both   \cite{LiYu2025} and \cite{Wu2003} also studied fully coupled FBSDEs with nonzero terminal values under additional domination--monotonicity and joint monotonicity conditions, respectively.

 In the finite-horizon setting, the dynamic programming principle (DPP) has been extensively studied for stochastic control and differential games with jumps, including those with BSDE-defined recursive payoffs; see \cite{LiPeng2009,BuckdahnHuLi2011,LiWei2015,Biswas2012,LuoLiLiWei2025,LiuWei2026} and the references therein.
 In comparison, corresponding studies on infinite horizons remain relatively limited, particularly for recursive problems involving jumps. Li and Zhao \cite{LiZhao2019} and Buckdahn, Li and Zhao \cite{BuckdahnLiZhao2021} investigated infinite-horizon recursive control and differential games under nonexpansivity conditions. Luo, Li and Wei \cite{LLW2025} studied infinite-horizon stochastic recursive control problems with jumps under dissipativity of the forward dynamics and strict monotonicity of the backward generator in its value variable. More recently, Huang and Wei \cite{HW2026} investigated non-autonomous infinite-horizon stochastic recursive differential games under monotonicity of the backward generator and time-dependent discounting conditions.

Motivated by the non-autonomous framework of Huang and Wei \cite{HW2026}, a natural question is whether infinite-horizon stochastic recursive differential games can be studied without imposing dissipativity of the forward dynamics or monotonicity of the backward generator in its value variable. We address this question in a jump-diffusion setting under time-integrable growth and stochastic Lipschitz conditions, allowing generally nonzero terminal payoffs determined by the limiting forward state.

We begin by developing an $L^p$ theory ($p\ges 2$) for decoupled infinite-horizon FBSDEs driven by a Brownian motion and an independent compensated Poisson random measure. The results include well-posedness, moment estimates, conditional stability, and tail estimates. In particular, the forward state converges in $L^p$ and almost surely to a limiting random variable $X_\infty$, while the backward solution converges to the generally nonzero terminal value $\Phi(X_\infty)$.
Based on these results, we establish the DPP for two-player zero-sum stochastic recursive differential games using Elliott--Kalton nonanticipative strategies and prove that the lower and upper value functions are deterministic. Both value functions converge locally uniformly to $\Phi$ as the initial time tends to infinity and are characterized as viscosity solutions of the corresponding non-autonomous, nonlocal HJBI equations. Under an additional integrable Lyapunov-barrier condition, we establish a comparison principle requiring uniform spatial Lipschitz continuity of only one competing function, yielding uniqueness in the class of continuous functions with uniform linear growth. Finally, we obtain a weaker, contact-wise Isaacs criterion ensuring that the game has a value.

Two technical difficulties deserve particular attention in our analysis.
First, the viscosity characterization requires precise short-time freezing estimates for the underlying jump-diffusion system. Although only short-time estimates are needed, unlike the framework with uniform Lipschitz bounds in \cite{LLW2025}, the time-dependent Lipschitz moduli in our setting are merely integrable and need not be locally bounded. Consequently, standard freezing estimates relying on uniform Lipschitz constants cannot be applied directly. To overcome this difficulty, we combine local eighth-moment estimates for the controlled forward state process, derived from the jump-size condition, with nested exit-time arguments, continuity on compact sets, and the absolute continuity of time integrals. These ingredients yield the required control-uniform $o(h)$ freezing error (see Lemma \ref{lem:jump-short-freeze}).

Second, establishing comparison and uniqueness for the infinite-horizon, non-autonomous, nonlocal HJBI equations requires additional care. Since the original generator need not be monotone in its backward value variable, we introduce a bounded exponential time transformation to obtain a strictly decreasing transformed generator, together with a compatible Lyapunov barrier to control the spatial penalization. Moreover, the time-integrable Lipschitz moduli do not guarantee that the coefficient and generator discrepancies vanish at the maximizing time in the doubling-of-variables argument, even when the corresponding errors can be controlled in an integrated sense. To overcome this difficulty, we construct a nonnegative temporal error envelope whose $L^1$ norm tends to zero and introduce an auxiliary time penalization whose derivative cancels the corresponding error. Combined with the nonlocal Jensen--Ishii method \cite{BarlesImbert2008} and suitable estimates for small and large jumps, this approach establishes comparison under the assumption that only one competing function is uniformly Lipschitz in space.

The remainder of the paper is organized as follows. Section 2 introduces the probabilistic framework and preliminary results. Section 3 establishes well-posedness, moment estimates, stability, and convergence results for infinite-horizon SDEs with jumps. Section 4 develops the corresponding $L^p$ theory for infinite-horizon BSDEs with generally nonzero terminal values. Section 5 formulates the stochastic recursive differential game and establishes the determinacy of the lower and upper value functions, the DPP, and their viscosity characterizations. It further proves comparison and uniqueness for the associated nonlocal HJBI equations and provides a contact-wise Isaacs criterion for the existence of a game value. Appendix A contains the proof of the short-time freezing estimate stated in Lemma \ref{lem:jump-short-freeze}.

\section{Preliminaries}

In this study, we work on a canonical Wiener--Poisson space
$(\Omega,\mathcal F,\dbP)$, constructed as follows.
\begin{itemize}
\item
Let
$\Omega_1:=C_0(\dbR;\dbR^d)
=\{\omega_1\in C(\dbR;\dbR^d):\omega_1(0)=0\}$,
endowed with the topology of uniform convergence on compact subsets of
$\dbR$. On $\Omega_1$, define the canonical process
$B_s(\omega_1):=\omega_1(s)$, $s\in\dbR$.
Let $\dbP_1$ be the Wiener measure under which
$(B_s)_{s\ges0}$ and $(B_{-s})_{s\ges0}$ are two independent
$d$-dimensional standard Brownian motions, and let $\mathcal F_1$ be
the $\dbP_1$-completion of the Borel $\sigma$-field on $\Omega_1$.
For $t\ges0$, set
$$
\mathcal F_t^B
:=
\sigma\{B_r:-\infty<r\les t\}\vee\mathcal N_{\dbP_1},
$$
where $\mathcal N_{\dbP_1}$ denotes the collection of all
$\dbP_1$-null sets.

\item
Let $E:=\dbR^\ell\setminus\{0\}$, equipped with its Borel
$\sigma$-field $\mathcal B(E)$, and let $\lambda$ be a
$\sigma$-finite L\'evy measure on $E$ satisfying
$\ds\int_E(1\wedge|e|^2)\lambda(\mathrm d e)<\infty$.
Let $\Omega_2$ be the collection of all point functions
$p:D_p\subset\dbR\to E$ with countable domain $D_p$.
Each $p\in\Omega_2$ is identified with the counting measure
$\ds
\mu(p,\mathrm d s,\mathrm d e)
:=
\sum_{r\in D_p}\delta_{(r,p(r))}(\mathrm d s,\mathrm d e)
$
on $\dbR\times E$.
Let $\mathcal F_2^0$ be the smallest $\sigma$-field on $\Omega_2$
such that the mapping
$p\mapsto\mu(p,(s,t]\times A)$ is measurable for every
$-\infty<s<t<\infty$ and $A\in\mathcal B(E)$.
Let $\dbP_2$ be the probability measure on
$(\Omega_2,\mathcal F_2^0)$ under which the canonical random measure
$\mu(\mathrm d s,\mathrm d e)$ is a Poisson random measure on $\dbR\times E$ with
intensity measure $\mathrm d s\,\lambda(\mathrm d e)$, and let $\mathcal F_2$ be the
$\dbP_2$-completion of $\mathcal F_2^0$.
The associated compensated Poisson random measure is denoted by
$\widetilde\mu(\mathrm d s,\mathrm d e):=\mu(\mathrm d s,\mathrm d e)-\mathrm d s\,\lambda(\mathrm d e)$.
For $t\ges0$, define
$$
\dot{\mathcal F}_t^\mu
:=
\sigma\big\{
\mu((s,r]\times A):
-\infty<s\les r\les t,\ A\in\mathcal B(E)
\big\},
$$
and set
$\ds
\mathcal F_t^\mu
:=
\bigcap_{r>t}\dot{\mathcal F}_r^\mu
\vee\mathcal N_{\dbP_2}.
$

\end{itemize}
Finally, we set
$\Omega:=\Omega_1\times\Omega_2$,
$\dbP:=\dbP_1\otimes\dbP_2$,
and  denote by $\mathcal F$ the $\dbP$-completion of
$\mathcal F_1\otimes\mathcal F_2$, and by
$\mathbb F=\{\mathcal F_t\}_{t\ges0}$ the $\dbP$-completed product
filtration
$\{\mathcal F_t^B\otimes\mathcal F_t^\mu\}_{t\ges0}$.
Under $\dbP$, the canonical Brownian motion $B$ and the canonical
Poisson random measure $\mu$ are independent.
 In particular, the initial $\sigma$-field
$\mathcal F_0$ is in general nontrivial.

Throughout this paper, $\langle\cdot,\cdot\rangle$ and
$|\cdot|$ denote the standard inner product and the
corresponding norm in Euclidean spaces, respectively.
For $z_1,z_2\in\dbR^{m\times d}$, we use the
Hilbert-Schmidt   inner product and norm, defined by
$$
\langle z_1,z_2\rangle
:=
\operatorname{tr}(z_1^{\mathsf T}z_2),
\qquad
|z_1|
:=
\bigl(
\operatorname{tr}(z_1^{\mathsf T}z_1)
\bigr)^{{\frac{1}{2}}},
$$
where $z_1^{\mathsf T}$ denotes the transpose of $z_1$.
For $0\les t<T\les\infty$,  $\mathcal P_{t,T}$ {denotes} the $\sigma$-field of  $\dbF$-predictable
subsets of $\Omega\times[t,T]$.
For $1\les p,q<\infty$  and
an integer $m\ges1$,
we introduce the following spaces.
\begin{itemize}

\item
$L_{\mathcal F_t}^p(\Omega;\dbR^m)$
is the space of all $\dbR^m$-valued,
$\mathcal F_t$-measurable random variables $\xi$ such that
$\ds
\mathbb E[|\xi|^p]<\infty.
$

 \item
$L_\lambda^q(E;\mathbb R^m)$ consists of
$\mathcal B(E)$-measurable functions
$\varphi:E\to\mathbb R^m$ such that
$ \ds
|\varphi|_{\lambda,q}^q
:=
\int_E|\varphi(e)|^q\,\lambda(\mathrm de)
<\infty.
$

\item
${\mathcal S^p}(t,T;\dbR^m)$
is the space of all $\dbR^m$-valued,
$\mathbb F$-adapted c\`adl\`ag processes
$\phi=\{\phi_s\}_{s\in[t,T]}$ such that
$\ds
\|\phi\|_{{\mathcal S^p}}^p
:=
\mathbb E\Big[
\sup_{s\in[t,T]}|\phi_s|^p
\Big]
<\infty.
$

\item
$\mathcal H^{p,q}(t,T;\dbR^m)$
is the space of all $\dbR^m$-valued,
$\mathbb F$-progressively measurable processes
$\phi=\{\phi_s\}_{s\in[t,T]}$ such that
$\ds
\|\phi\|_{\mathcal H^{p,q}}^p
:=
\mathbb E\Big[
\Big(
\int_t^T|\phi_s|^q\,\mathrm ds
\Big)^{\frac pq}
\Big]
<\infty.
$

\item
$\mathcal K_{\lambda,q}^p(t,T;\dbR^m)$
is the space of all
${\mathcal P_{t,T}}\otimes\mathcal B(E)$-measurable processes
$
\phi:\Omega\times[t,T]\times E\longrightarrow\dbR^m
$
such that
$\ds
\|\phi\|_{\mathcal K_{\lambda,q}^p}^p
:=
\mathbb E\Big[
\Big(
\int_t^T\int_E
|\phi_s(e)|^q\lambda(\mathrm d e)\,\mathrm ds
\Big)^{\frac pq}
\Big]
<\infty.
$

\item
$L_{\mathbb F}^{p,\infty}(t,T;\dbR^m)$
is the space of all $\dbR^m$-valued,
$\mathbb F$-progressively measurable processes
$\phi=\{\phi_s\}_{s\in[t,T]}$ such that
$\ds
\|\phi\|_{L_{\mathbb F}^{p,\infty}}^p
:=
\int_t^T
\|\phi_s\|_\infty^p\,\mathrm ds
<\infty.
$

\item
$L_{\mathbb F,\lambda}^{p,\infty}(t,T;\dbR^m)$
is the space of all
$\mathcal P_{t,T}\otimes\mathcal B(E)$-measurable processes
$
\phi:\Omega\times[t,T]\times E\longrightarrow\dbR^m
$
such that
$\ds
\|\phi\|_{L_{\mathbb F,\lambda}^{p,\infty}}^p
:=
\int_t^T\int_E
\|\phi_s(e)\|_\infty^p
\lambda(\mathrm d e)\,\mathrm ds
<\infty.
$

 When $m=1$, we abbreviate
$L_\lambda^q(E;\mathbb R)$,
$L_{\mathbb F}^{p,\infty}(t,T;\mathbb R)$,
and
$L_{\mathbb F,\lambda}^{p,\infty}(t,T;\mathbb R)$
as
$L_\lambda^q(E)$,
$L_{\mathbb F}^{p,\infty}(t,T)$,
and
$L_{\mathbb F,\lambda}^{p,\infty}(t,T)$,
respectively.

\item
Set
$
\mathbb B_\lambda^p(t,T)
:=
\mathcal S^p(t,T;\dbR)
\times\mathcal H^{p,2}(t,T;\dbR^d) \times\bigl(
\mathcal K_{\lambda,2}^p(t,T;\dbR)
\cap\mathcal K_{\lambda,p}^p(t,T;\dbR)
\bigr).
$

\end{itemize}
In the above, for any $\dbR^m$-valued random variable $\xi$,
we  use
$\ds
\|\xi\|_\infty
:=
\operatorname*{ess\,sup}_{\omega\in\Omega}
|\xi(\omega)|.
$
For $T=\infty$, all time integrals and suprema
are taken over $[t,\infty)$.
We emphasize that, for $p>2$, the spaces
$\mathcal K_{\lambda,2}^p(t,T;\dbR^m)$
and $\mathcal K_{\lambda,p}^p(t,T;\dbR^m)$
are generally incomparable, particularly since
$\lambda(E)$ is not assumed to be finite.
Throughout this paper, $C>0$ denotes a generic
constant that may vary from line to line but is
independent of $T$. Dependence on $T$ is indicated
by $C_T$, and other dependencies are specified
when necessary.

 Before concluding this section, we recall two moment
estimates for Poisson random measures that will be
used repeatedly in the sequel.

\begin{lemma}
\label{lem:poisson-moment-estimates}\sl
Let $p\ges2$ and $t<T\les\infty$. For any $\dbF$-predictable
process $K$ such that
$$
{\mathbb E\(
\int_t^T\int_E |K_r(e)|^2\,\mu(\mathrm dr,\mathrm de)
\)^{{\frac{p}{2}}}<\infty,}
$$
the following estimates hold, $\dbP$-a.s.,
\begin{equation}\label{eq:poisson-p-moment}
\mathbb E_t\[
\int_t^T\int_E
|K_r(e)|^p\,\lambda(\mathrm de)\,\mathrm dr
\]
\les
\mathbb E_t\[
\(
\int_t^T\int_E
|K_r(e)|^2\,\mu(\mathrm dr,\mathrm de)
\)^{{\frac{p}{2}}}
\],
\end{equation}
and
\begin{equation}\label{eq:poisson-compensator-moment}
 \mathbb E_t\[
\(
\int_t^T\int_E
|K_r(e)|^2\,\lambda(\mathrm de)\,\mathrm dr
\)^{{\frac{p}{2}}}
\] \les
\Big(\frac p2\Big)^{{\frac{p}{2}}}
\mathbb E_t\[
\(
\int_t^T\int_E
|K_r(e)|^2\,\mu(\mathrm dr,\mathrm de)
\)^{{\frac{p}{2}}}
\].
\end{equation}

\end{lemma}

The first inequality follows from the conditional
compensation formula and the elementary inequality
$$
\int_t^T\int_E |K_r(e)|^p\,
\mu(\mathrm dr,\mathrm de)
\les
\(
\int_t^T\int_E |K_r(e)|^2\,
\mu(\mathrm dr,\mathrm de)
\)^{\frac p2}.
$$
The second inequality follows directly from
\cite[Lemma~3.1]{LiWei2014}.
The case $T=\infty$ follows by monotone convergence,
and the corresponding unconditional estimates
are obtained by taking expectations.

 \section{Infinite-horizon SDEs with jumps}

Let $p\ges 2$ be fixed.
For any $t\ges 0$ and
 $\xi\in
L_{\mathcal F_t}^{p}(\Omega;\mathbb R^n),
$
we consider the following infinite-horizon SDE with jumps:
\begin{equation}\label{eq:infinite-horizon-sde}
\begin{cases}
\displaystyle
\mathrm d X_s
=
b(s,X_s) \mathrm d s
+
\sigma(s,X_s)\mathrm d B_s
+
\int_E\gamma(s,e,X_{s-})\,
\widetilde{\mu}(\mathrm d s,\mathrm d e),
\qquad s\ges t,
\\[1.2ex]
X_t=\xi.
\end{cases}
\end{equation}
Here, we assume
$
  b :\Omega\times[0,\infty)\times\mathbb R^n
   \longrightarrow\mathbb R^n,$ $
\sigma :\Omega\times[0,\infty)\times\mathbb R^n
   \longrightarrow\mathbb R^{n\times d},$ $
\gamma :\Omega\times[0,\infty)\times E\times\mathbb R^n
   \longrightarrow\mathbb R^n
$ to satisfy the following conditions.

\begin{description}
  \item[(A1)] For every $x\in\mathbb R^n$, the processes
$b(\cdot,x)$ and $\sigma(\cdot,x)$ are $\dbF$-progressively
measurable, while $\gamma(\cdot,\cdot,x)$ is $\dbF$-predictable in
$(\omega,s)$ and measurable in $e$.

\item[(A2)] There exist nonnegative processes
$
l_b\in L_{\mathbb F}^{1,\infty}(0,\infty),
$
$
l_\sigma\in L_{\mathbb F}^{2,\infty}(0,\infty),
$
and
$
l_\gamma
\in
L_{\mathbb F,\lambda}^{2,\infty}(0,\infty)
\cap
L_{\mathbb F,\lambda}^{p,\infty}(0,\infty),
$
such that, for all $x,x'\in\mathbb R^n$,
$$\begin{aligned}&
|b(s,x)-b(s,x')|
\les l_b(s)|x-x'|,\q
|\sigma(s,x)-\sigma(s,x')|
\les l_\sigma(s)|x-x'|,\q \mathrm d\mathbb P\otimes\mathrm d s\mbox{-a.e.}\\
&
|\gamma(s,e,x)-\gamma(s,e,x')|
\les l_\gamma(s,e)|x-x'|,\q \mathrm d\mathbb P\otimes\mathrm d s\otimes\lambda(\mathrm d e)\mbox{-a.e.}
\end{aligned}$$

\item[(A3)] The coefficients at the origin satisfy
$
b(\cdot,0)
 \in \mathcal H^{p,1}(0,\infty;\dbR^n),$ $
\sigma(\cdot,0)
 \in \mathcal H^{p,2}(0,\infty;\dbR^{n\times d}),$ $
\gamma(\cdot,\cdot,0)
 \in \mathcal K_{\lambda,2}^p(0,\infty;\dbR^n)
\cap \mathcal K_{\lambda,p}^p(0,\infty;\dbR^n).$

\end{description}

For notational convenience and later use, we introduce the following
quantities. For $r\ges0$ and $0\les t<s\les\infty$, set
\be\label{nota: 2.1}
\begin{aligned}
&q_p(r)
:={}
\|l_b(r)\|_\infty
+\|l_\sigma(r)\|_\infty^2 +\int_E
\big(
\|l_\gamma(r,e)\|_\infty^2
+\|l_\gamma(r,e)\|_\infty^p
\big)\lambda(\mathrm de),\\
&{q_p^0(s)}:=
\|b(s,0)\|_\infty+\|\sigma(s,0)\|_\infty^2
+\int_E\big(
\|\gamma(s,e,0)\|_\infty^2
+\|\gamma(s,e,0)\|_\infty^p\big)\lambda(\mathrm d e),\\
&\Theta_p(t,s):=
\Big(\int_t^s\|l_b(r)\|_\infty\,\mathrm d r\Big)^p
+\Big(\int_t^s\|l_\sigma(r)\|_\infty^2\,\mathrm d r\Big)^{\frac p2}\\
&\qq\qq\q+\Big(\int_t^s\int_E
\|l_\gamma(r,e)\|_\infty^2\lambda(\mathrm d e)\,\mathrm d r\Big)^{\frac p2}
+\int_t^s\int_E
\|l_\gamma(r,e)\|_\infty^p\lambda(\mathrm d e)\,\mathrm d r,\\
&\eta_p(t,s):=
\Big(\int_t^s|b(r,0)|\,\mathrm d r\Big)^p
+\Big(\int_t^s|\sigma(r,0)|^2\,\mathrm d r\Big)^{\frac p2}+\Big(\int_t^s
|\gamma(r,\cdot,0)|_{\lambda,2}^2\,\mathrm d r\Big)^{\frac p2}
+\int_t^s
|\gamma(r,\cdot,0)|_{\lambda,p}^p\,\mathrm d r.
\end{aligned}
\ee
By {\bf (A2)} and {\bf (A3)}, we have
$q_p\in L^1(0,\infty)$ and
$\mathbb E\big[\eta_p(0,\infty)\big]<\infty$.
Consequently,
\bel{}
\lim_{t\to\infty}\Theta_p(t,\infty)=0,
\qquad
\lim_{t\to\infty}
\mathbb E\big[\eta_p(t,\infty)\big]=0.
\ee

The following study extends
\cite[Proposition~2.1]{LiYu2025}
to the setting of SDEs with jumps and $p\ges2$.
Unlike the method of continuation used therein,
our proof is based on direct $L^p$ estimates.

\begin{theorem}\label{thm:infinite-horizon-sde}\sl
Suppose that {\bf (A1)--(A3)} hold. Then, for every
$
(t,\xi)\in[0,\infty)\times
L_{\mathcal F_t}^{p}(\Omega;\mathbb R^n),
$
SDE \eqref{eq:infinite-horizon-sde} admits a unique strong solution
$
X\in
{\mathcal S^p}(t,\infty;\mathbb R^n).
$
Moreover, there exists a constant $C_p>0$, depending only on $p$
and the integrability bounds in {\bf (A2)}, such that
\begin{equation}\label{eq:state-estimate}
\mathbb E_t\Big[
\sup_{s\ges t}|X_s|^p
\Big]
\les
C_p\Big\{
|\xi|^p
+
\mathbb E_t[\eta_p(t,\infty)]
\Big\},
\qquad \dbP\mbox{-a.s.}
\end{equation}
Furthermore, let $\bar X$ denote the solution of
\eqref{eq:infinite-horizon-sde} corresponding to another initial value
$
\bar\xi\in
L_{\mathcal F_t}^{p}(\Omega;\mathbb R^n)
$
and another coefficient triple
$(\bar b,\bar\sigma,\bar\gamma)$ satisfying
{\bf (A1)--(A3)}. Then
\begin{equation}\label{eq:stability-estimate}
\begin{aligned}
&\mathbb E_t\[
\sup_{s\ges t}
|X_s-\bar X_s|^p
\]\!
\les\! C_p\Big\{
 |\xi-\bar\xi|^p
+
\mathbb E_t\[
\(\!
\int_t^\infty\!|(b -\bar b)(r,\bar X_r)|\,\mathrm d r\!
\)^p+
\(
\int_t^\infty\!|(\sigma -\bar \sigma)(r,\bar X_r)|^2\,\mathrm d r
\)^{\frac p2} \\
&\qq\qq+
\(
\int_t^\infty\!\!\int_E
|(\gamma -\bar \gamma)(r, e,\bar X_{r-} )|^2
\lambda(\mathrm d e)\,\mathrm d r
\)^{\frac p2}
+
\int_t^\infty\!\!\int_E
|(\gamma -\bar \gamma)(r, e,\bar X_{r-})|^p
\lambda(\mathrm d e)\,\mathrm d r
\]
\Big\},\quad \dbP\mbox{-a.s.}
\end{aligned}
\end{equation}
\end{theorem}

\begin{proof}

\noindent{\it Step 1. Finite-horizon solvability.}
Although finite-horizon well-posedness results for
SDEs with jumps are well established
(see, e.g., \cite[Chapter 3]{Situ2005} and
\cite[Chapter 6]{Applebaum2009}),
we provide a proof under {\bf (A1)}--{\bf (A3)}
to obtain a deterministic time partition based on
$\Theta_p$, which will be used in the subsequent analysis.

  Fix $T>t$, and let
$
t\les\tau<\varrho\les T.
$
For any
$
\eta\in L_{\mathcal F_\tau}^p(\Omega;\mathbb R^n)
$
and
$
x_\cd\in{\mathcal S^p}(\tau,\varrho;\mathbb R^n),
$
define
$$
(\Gamma x)_s
:=
\eta
+\int_\tau^s b(r,x_r)\,\mathrm d r
+\int_\tau^s \sigma(r,x_r)\,\mathrm d B_r +
\int_\tau^s\int_E
\gamma(r,e,x_{r-})\,
\widetilde\mu(\mathrm d r,\mathrm d e),
\qquad s\in[\tau,\varrho].
$$
By {\bf (A1)--(A3)},
$
\Gamma
$
maps
$
{\mathcal S^p}(\tau,\varrho;\mathbb R^n)
$
into itself.
Moreover, for any
$
x_\cd,\bar x_\cd\in
{\mathcal S^p}(\tau,\varrho;\mathbb R^n),
$
the Burkholder--Davis--Gundy (BDG) inequalities and {\bf (A2)} yield
$$
\mathbb E\Big[
\sup_{s\in[\tau,\varrho]}
|(\Gamma x)_s-(\Gamma\bar x)_s|^p
\Big]
\les
c_p\Theta_p(\tau,\varrho)
\mathbb E\Big[
\sup_{s\in[\tau,\varrho]}|x_s-\bar x_s|^p
\Big].
$$

Note that, for
 $t\les\tau<\varrho\les\infty$,
$$
\Theta_p(\tau,\varrho)
\les
\Big(\int_\tau^\varrho q_p(r)\,\mathrm d r\Big)^p
+
2\Big(\int_\tau^\varrho q_p(r)\,\mathrm d r\Big)^{\frac p2}
+
\int_\tau^\varrho q_p(r)\,\mathrm d r.
$$
Since $q_p\in L^1(0,\infty)$, choose $\varepsilon>0$ sufficiently
small such that
$\ds
c_p\big(
\varepsilon^p+2\varepsilon^{\frac p2}+\varepsilon
\big)
\les\frac12.
$
 We can then choose a finite deterministic partition
$
t=s_0<s_1<\cdots<s_N=\infty
$
such that
$$
\int_{s_{j-1}}^{s_j}q_p(r)\,\mathrm d r
\les\varepsilon,
\quad j=1,\ldots,N,
\qquad
N\les 1+\varepsilon^{-1}\int_0^\infty q_p(r)\,\mathrm d r.
$$
Consequently,
$ \ds
c_p\Theta_p(s_{j-1},s_j)\les\frac12,$
$  j=1,\ldots,N.
$

This partition is independent of $T$ and will also be used below.
For the fixed $T>t$, we truncate it at $T$ and apply
Banach's fixed-point theorem successively on the intervals
$[s_{j-1},s_j\wedge T]$ with $s_{j-1}<T$,
starting from $\xi$ and using the terminal value of each
local solution as the initial value on the next interval.
Concatenating these local solutions yields a unique
strong solution
$
X^T\in{\cS^p}(t,T;\mathbb R^n)
$
of \eqref{eq:infinite-horizon-sde} on $[t,T]$.

\ms

 {\it Step 2. Infinite-horizon solvability.}
Fix $T>t$, and let $X^T$ be the solution constructed in Step 1.
The conditional
BDG inequalities and {\bf (A2)} yield
\begin{equation} \label{supX-T}
\mathbb E_\tau\Big[
\sup_{s\in[\tau,\varrho]}|X_s^T|^p
\Big]
 \les c_p\Big\{|X_\tau^T|^p
+ \Theta_p(\tau,\varrho)
\mathbb E_\tau\Big[
\sup_{s\in[\tau,\varrho]}|X_s^T|^p
\Big]
+ \mathbb E_\tau[\eta_p(\tau,\varrho)]\Big\},
\end{equation}
with $c_p$ chosen as in Step 1.

Let $t=s_0<s_1<\cdots<s_N=\infty$ be the partition
constructed in Step 1, and let $m$ be determined by
$s_{m-1}<T\les s_m$.
On each interval $[s_{j-1},s_j\wedge T]$, $j=1,\ldots,m$,
we have
$
c_p\Theta_p(s_{j-1},s_j\wedge T)\les\frac12.
$
Absorbing the supremum term in the preceding estimate \eqref{supX-T}
and iterating with the tower property, we obtain
\begin{equation}\label{X-T-2}
\mathbb E_t\Big[
\sup_{s\in[t,T]}|X_s^T|^p
\Big]
 \les
(1+2c_p)^m
\Big\{
|\xi|^p+
\mathbb E_t\Big[
\sum_{j=1}^m\eta_p(s_{j-1},s_j\wedge T)
\Big]
\Big\} \les
C_p\Big\{
|\xi|^p+\mathbb E_t[\eta_p(t,T)]
\Big\}.
\end{equation}
Here we used $m\les N$.
By the bound on $N$ in Step 1, the constant $C_p>0$
depends only on $p$ and the integrability bounds
in {\bf (A2)}, and is independent of $t$ and $T$.

Set $T_n:=t+n$, $n\ges1$.
By finite-horizon pathwise uniqueness, the solutions
$X^{T_n}$ are consistent on overlapping intervals,
outside a common null set.
Thus,
$
X_s:=X_s^{T_n}
$
for $s\in[t,T_n]$ defines an
$\mathbb F$-adapted c\`adl\`ag process on $[t,\infty)$
satisfying \eqref{eq:infinite-horizon-sde}.
Applying the preceding estimate \eqref{X-T-2}  with $T=T_n$
and letting $n\to\infty$ by conditional monotone
convergence yields \eqref{eq:state-estimate}.
Taking expectations and using {\bf (A3)},
we obtain
$
X\in{\cS^p}(t,\infty;\mathbb R^n).
$

Finally, any two solutions are indistinguishable
on each $[t,T_n]$ by Step 1.
A countable intersection of these full-probability
events gives indistinguishability on $[t,\infty)$,
and hence uniqueness.

\ms%

{\it Step 3. {Continuous dependence}.}
Let $X$ and $\bar X$ be the solutions of
\eqref{eq:infinite-horizon-sde} corresponding to
$(\xi,b,\sigma,\gamma)$ and
$(\bar\xi,\bar b,\bar\sigma,\bar\gamma)$, respectively.
Set
$
\widehat X:=X-\bar X
$
and
$
\widehat\xi:=\xi-\bar\xi.
$
Then
\begin{equation}\label{eq:sde-difference}
\widehat X_s
 =\widehat\xi
+\int_t^s b_1(r,\widehat X_r)\,\mathrm d r
+\int_t^s \sigma_1(r,\widehat X_r)\,\mathrm d B_r
+\int_t^s\int_E
\gamma_1(r,e,\widehat X_{r-})
\,\widetilde\mu(\mathrm d r,\mathrm d e),
\qquad s\ges t,
\end{equation}
where
$
b_1(r,x)
 :=b(r,x+\bar X_r)-\bar b(r,\bar X_r),$ $
\sigma_1(r,x)
 :=\sigma(r,x+\bar X_r)-\bar\sigma(r,\bar X_r),$ $
\gamma_1(r,e,x)
 :=\gamma(r,e,x+\bar X_{r-})
-\bar\gamma(r,e,\bar X_{r-}).
$
The coefficients $(b_1,\sigma_1,\gamma_1)$ satisfy
{\bf (A1)--(A2)} with the same Lipschitz moduli as
$(b,\sigma,\gamma)$.
Moreover, {\bf (A2)--(A3)} for both coefficient triples,
together with
$
\bar X\in{\cS^p}(t,\infty;\mathbb R^n),
$
ensure the integrability required  in {\bf (A3)} by $b_1$, $\sigma_1$ and $\gamma_1$.
Applying \eqref{eq:state-estimate} to \eqref{eq:sde-difference} yields \eqref{eq:stability-estimate}.
\end{proof}

%
%

Unlike the dissipative framework of \cite{LuoLiWei2026}, the present setting does not require the state process to vanish at infinity. Under {\bf (A1)--(A3)}, it converges to a possibly nonzero random variable, as shown in the following lemma.

 \begin{lemma}\label{lem:state-limit}\sl
Suppose that  {\bf (A1)--(A3)} hold.
Fix $t\ges0$ and
$\xi\in L^p_{\mathcal F_t}(\Omega;\mathbb R^n)$,
and let $X$ be the unique solution of
\eqref{eq:infinite-horizon-sde}.
Then there exists
$
X_\infty\in L^p_{\mathcal F_\infty}(\Omega;\mathbb R^n)
$
such that
$
X_T\to X_\infty
$
in $L^p(\Omega;\mathbb R^n)$ and $\mathbb P$-a.s.,
as $T\to\infty$.
Moreover,
\begin{equation}\label{eq:state-tail}
\lim_{T\to\infty}
\mathbb E\Big[
\sup_{s\ges T}|X_s-X_T|^p
\Big]=0.
\end{equation}
\end{lemma}

 \begin{proof}
By Theorem~\ref{thm:infinite-horizon-sde},
$
X\in{\cS^p}(t,\infty;\mathbb R^n).
$
Since $p\ges2$, {\bf (A2)--(A3)} imply
$$
\mathbb E\Big[
\Big(\int_t^\infty|b(r,X_r)|\,\mathrm d r\Big)^2
+\int_t^\infty
\big(
|\sigma(r,X_r)|^2
+|\gamma(r,\cdot,X_{r-})|_{\lambda,2}^2
\big)\,\mathrm d r
\Big]<\infty.
$$
The martingale convergence theorem
(see \cite{JacodShiryaev2003}),
together with the absolute convergence of the drift integral,
therefore yields an $\mathcal F_\infty$-measurable random variable
$X_\infty$ such that
$
X_T\to X_\infty
$
$\mathbb P$-a.s., as $T\to\infty$. Moreover,
$\ds
|X_\infty|\les\sup_{s\ges t}|X_s|,
$
so that
$
X_\infty\in L^p_{\mathcal F_\infty}(\Omega;\mathbb R^n).
$

The almost sure convergence of $X$ implies
$\ds
\sup_{s\ges T}|X_s-X_T|\to0
$
$\mathbb P$-a.s.\ as $T\to\infty$.
Since $\ds
\sup_{s\ges T}|X_s-X_T|^p
\les 2^p\sup_{s\ges t}|X_s|^p
\in L^1(\Omega),
$
the dominated convergence theorem yields
both \eqref{eq:state-tail} and
$X_T\to X_\infty$ in $L^p(\Omega;\mathbb R^n)$.
\end{proof}

We next establish increment estimates for $X$
over both short intervals and the infinite horizon.

\begin{corollary}\label{cor:short-time-state-increment}\sl
Suppose that {\bf(A1)--(A3)} hold.
For each $(t,x)\in[0,\infty)\times\dbR^n$,
let $X^{t,x}$ denote the solution of
\eqref{eq:infinite-horizon-sde} with
initial value $X_t^{t,x}=x$.
Then the following estimates hold.

 (i)
There exists a deterministic
$\delta_0\in(0,1]$, depending only on
$p$ and $q_p$, such that, for every
$t\ges0$, $x\in\dbR^n$, and
$0<\delta\les\delta_0$,
\begin{equation}\label{eq:short-time-state-increment}
 \mathbb E_t\Big[
\sup_{s\in[t,t+\delta]}|X_s^{t,x}-x|^p
\Big] \les C_p\Big\{
|x|^p\Theta_p(t,t+\delta)
+\mathbb E_t[\eta_p(t,t+\delta)]
\Big\},
\qquad \dbP\mbox{-a.s.}
\end{equation}
Moreover, for every $t\ges0$ and
$x\in\dbR^n$,
\begin{equation}\label{eq:global-state-increment}
 \mathbb E_t\Big[
\sup_{s\ges t}|X_s^{t,x}-x|^p
\Big] \les C_p\Big\{
|x|^p\Theta_p(t,\infty)
+\mathbb E_t[\eta_p(t,\infty)]
\Big\},
\qquad \dbP\mbox{-a.s.}
\end{equation}

 (ii)
Fix $t\ges0$. Suppose additionally that
the following condition holds:

\noindent{\bf(A4)}
There exists a constant $L_*>0$ such that
$
q_p(r)+q_p^0(r)\les L_*
$
for a.e.\ $r\in[t,t+\delta_0]$.

Then, for every $x\in\dbR^n$ and
$0<\delta\les\delta_0$,
\begin{equation}\label{eq:short-time-state-increment-rate}
 \mathbb E_t\Big[
\sup_{s\in[t,t+\delta]}|X_s^{t,x}-x|^p
\Big] \les C_{p,L_*}(1+|x|^p)
\bigl(\delta^p+\delta^{\frac p2}+\delta\bigr),
\qquad \dbP\mbox{-a.s.}
\end{equation}
Moreover, if $\gamma\equiv0$ or $p=2$, then
\begin{equation}
\label{eq:short-time-state-increment-continuous-rate}
 \mathbb E_t\Big[
\sup_{s\in[t,t+\delta]}|X_s^{t,x}-x|^p
\Big] \les C_{p,L_*}(1+|x|^p)
\delta^{\frac p2},
\qquad \dbP\mbox{-a.s.}
\end{equation}
\end{corollary}

\begin{proof}
Fix $t\ges0$ and $x\in\dbR^n$,
and write $X:=X^{t,x}$.

{\rm(i)}
The conditional BDG inequalities,
together with {\bf(A2)}, yield
\begin{equation}\label{X-x}
\begin{aligned}
 &\mathbb E_t\Big[
\sup_{s\in[t,t+\delta]}|X_s-x|^p
\Big]  \les c_p\mathbb E_t\Big[
\Big(\int_t^{t+\delta}|b(r,X_r)|\,\mathrm d r\Big)^p
+\Big(\int_t^{t+\delta}|\sigma(r,X_r)|^2
\,\mathrm d r\Big)^{\frac p2}\\
&\qquad\qquad\qquad\qquad\qquad\qquad
+\Big(\int_t^{t+\delta}
|\gamma(r,\cdot,X_{r-})|_{\lambda,2}^2
\,\mathrm d r\Big)^{\frac p2}
+\int_t^{t+\delta}
|\gamma(r,\cdot,X_{r-})|_{\lambda,p}^p
\,\mathrm d r
\Big]\\
& \les c_p\mathbb E_t[\eta_p(t,t+\delta)]
+c_p\Theta_p(t,t+\delta)
\Big\{
|x|^p+
\mathbb E_t\Big[
\sup_{s\in[t,t+\delta]}|X_s-x|^p
\Big]
\Big\},
\end{aligned}
\end{equation}
where $c_p>0$ depends only on $p$.

Since $q_p\in L^1(0,\infty)$,
the estimate for $\Theta_p$ in Step~1
of the proof of
Theorem~\ref{thm:infinite-horizon-sde}
and the absolute continuity of
the Lebesgue integral imply
$$
\lim_{\delta\downarrow0}
\sup_{t\ges0}\Theta_p(t,t+\delta)=0.
$$
We may therefore choose a deterministic
$\delta_0\in(0,1]$, depending only on
$p$ and $q_p$, such that
$
c_p\Theta_p(t,t+\delta)\les\frac12
$
for all $t\ges0$ and
$0<\delta\les\delta_0$.
Absorbing the last term on the right-hand
side of \eqref{X-x} into the left-hand
side proves
\eqref{eq:short-time-state-increment}.

To establish \eqref{eq:global-state-increment},
apply the stability estimate in
Theorem~\ref{thm:infinite-horizon-sde}
to $X^{t,x}$ and the constant process
$\bar X_s\equiv x$, corresponding to
the zero coefficient triple
$(\bar b,\bar\sigma,\bar\gamma)=(0,0,0)$.
By the Lipschitz bounds in {\bf(A2)}
and the definitions of $\Theta_p$ and
$\eta_p$, we obtain
$$
  \mathbb E_t\Big[
\sup_{s\ges t}|X_s^{t,x}-x|^p
\Big] \les C_p\Big\{
|x|^p\Theta_p(t,\infty)
+\mathbb E_t[\eta_p(t,\infty)]
\Big\},
$$
which proves \eqref{eq:global-state-increment}.

{\rm (ii)} Under the additional condition {\bf(A4)},
\bel{Theta-eta-p}
\Theta_p(t,t+\delta)
+\mathbb E_t[\eta_p(t,t+\delta)]
\les C_{p,L_*}
\bigl(\delta^p+\delta^{\frac p2}+\delta\bigr),
\qquad 0<\delta\les\delta_0.
\ee
Combining this with \eqref{eq:short-time-state-increment}
proves \eqref{eq:short-time-state-increment-rate}.

If $\gamma\equiv0$, the two jump terms in the initial
estimate vanish, and the same argument gives
$$
\mathbb E_t\Big[
\sup_{s\in[t,t+\delta]}|X_s-x|^p
\Big]
\les C_{p,L_*}(1+|x|^p)
\bigl(\delta^p+\delta^{\frac p2}\bigr)\les C_{p,L_*}(1+|x|^p)
 \delta^{\frac p2} .
$$
This proves
\eqref{eq:short-time-state-increment-continuous-rate}.
For $p=2$, the same conclusion follows directly from
\eqref{eq:short-time-state-increment-rate}.

\end{proof}

\begin{remark}\label{rem:short-time-condition}\sl
The local boundedness assumption {\bf (A4)}
is sufficient but not necessary for the
short-time estimate
\eqref{eq:short-time-state-increment-rate}.
Indeed, its proof only requires
\eqref{Theta-eta-p}, which can also be
obtained under weaker integral conditions.
For example, fix $t\ges0$ and suppose that
there exist constants $\delta_1\in(0,1]$
and $L>0$ such that, for all
$0<\delta\les\delta_1$,
$$
\begin{aligned}
&\Big(\int_t^{t+\delta}
\big(
\|l_b(r)\|_\infty+\|b(r,0)\|_\infty
\big)\,\mathrm dr\Big)^2
+\int_t^{t+\delta}
\big(
\|l_\sigma(r)\|_\infty^2
+\|\sigma(r,0)\|_\infty^2
\big)\,\mathrm dr
\les L\delta,\\
&\int_t^{t+\delta}\int_E
\big(
\|l_\gamma(r,e)\|_\infty^2
+\|\gamma(r,e,0)\|_\infty^2
\big)\lambda(\mathrm de)\,\mathrm dr
\les L\delta^{\frac2p},\\
&\int_t^{t+\delta}\int_E
\big(
\|l_\gamma(r,e)\|_\infty^p
+\|\gamma(r,e,0)\|_\infty^p
\big)\lambda(\mathrm de)\,\mathrm dr
\les L\delta.
\end{aligned}
$$
Then
$$
  \Theta_p(t,t+\delta)
+\mathbb E_t[\eta_p(t,t+\delta)] \les
 L^{\frac p2}\delta^{\frac p2}
+\big(L^{\frac p2}+L\big)\delta,
\qquad \mathbb P\text{-a.s.}
$$
Thus, \eqref{Theta-eta-p} remains valid,
with a possibly different constant, for
$0<\delta\les\delta_0\wedge\delta_1$.
These integral conditions are weaker than
{\bf (A4)}, allowing certain coefficient
bounds to be unbounded in time.
\end{remark}

\section{Infinite-horizon BSDEs with jumps}

With the known process $X$ satisfying \eqref{eq:infinite-horizon-sde}, consider the following infinite-horizon BSDE with jumps,
\begin{equation}\label{eq:infinite-horizon-bsde}
\left\{\2n
\begin{aligned}
 &
Y_s
=
Y_T
+
\int_s^T
f\Big(
r,X_r,Y_r,Z_r,
\int_E\psi(r,e,K_r(e))\lambda(\mathrm d e)
\Big)\mathrm d r
\\
&\qquad
-
\int_s^T Z_r\,\mathrm d B_r
-
\int_s^T\int_E
 K_r(e)\widetilde\mu(\mathrm d r,\mathrm d e),
\q t\les s\les T<\infty;
\\[1ex]
 &
 Y_T\to\Phi(X_\infty)
 \quad \mathbb P\text{-a.s. as }T\to\infty,
\end{aligned}
\right.
\end{equation}
where the coefficients $f:\Omega \times [0,\infty)  \times  \mathbb{R}^{n} \times  \mathbb{R} \times  \mathbb{R}^{d} \times  \mathbb{R}\rightarrow \mathbb{R}$,
$ \psi:\Omega  \times  [0,\infty)  \times E \times  \mathbb{R}\rightarrow \mathbb{R}$,
$\Phi:\Omega\times\mathbb{R}^{n}\rightarrow\mathbb{R}$  satisfy the following conditions.

 \begin{description}

 \item[(B1)]
For every
$
(x,y,z,k)\in
\mathbb R^n\times\mathbb R\times\mathbb R^d\times\mathbb R,
$
$f(\cdot,\cdot,x,y,z,k)$ is $\mathbb F$-progressively
measurable. For every $k\in\mathbb R$,
$\psi(\cdot,\cdot,\cdot,k)$ is
${\mathcal P_{0,\infty}\otimes\mathcal B(E)}$-measurable.
Moreover, for every
$x\in\mathbb R^n$,
$\Phi(\cdot,x)$ is $\mathcal F_\infty$-measurable.

 \item[(B2)]
There exist nonnegative
processes $
l_x,l_y\in L_{\mathbb F}^{1,\infty}(0,\infty),
$ $
l_z,l_k\in L_{\mathbb F}^{2,\infty}(0,\infty),
$ and a nonnegative deterministic measurable function
$
\rho:[0,\infty)\times E\to[0,\infty),
$
and a constant $C_\Phi\ges0$   such that $
\rho(r,\cdot)\in L_\lambda^2(E)
$
for a.e.\ $r\ges0$,
$\ds
\int_0^\infty
\|l_k(r)\|_\infty^2
|\rho(r,\cdot)|_{\lambda,2}^2
\,\mathrm dr<\infty,
$ and the following inequalities hold:
$
x,x'\in\mathbb R^n,
$
$
y,y',k,k'\in\mathbb R,
$
and
$
z,z'\in\mathbb R^d,
$
$$
\begin{aligned}
&
|f(s,x,y,z,k)-f(s,x',y',z',k')|
\\
&\quad\les
l_x(s)|x-x'|
+l_y(s)|y-y'|
+l_z(s)|z-z'|
+l_k(s)|k-k'|,
\qquad
\mathrm d\mathbb P\otimes\mathrm d s\text{-a.e.},
\\[1mm]
&
|\psi(s,e,k)-\psi(s,e,k')|
\les
\rho(s,e)|k-k'|,
\qquad
\mathrm d\mathbb P\otimes\mathrm d s\otimes\lambda(\mathrm d e)\text{-a.e.},
\\[1mm]
&
 |\Phi(x)-\Phi(x')|
 \les
 C_\Phi|x-x'|,
 \qquad
 \mathbb P\text{-a.s.}
\end{aligned}
$$
\item[(B3)]
For the fixed $p\ges2$,
$
f(\cdot,0,0,0,0)
\in\mathcal H^{p,1}(t,\infty;\dbR),$
and
$\ds
\mathbb E\[
\(
\int_t^\infty
|\psi(r,\cdot,0)|_{\lambda,1}^2\mathrm d r
\)^{\frac p2}
\]
<\infty.
$

 \item[(B4)]
There exists a constant $B_{\Phi}>0$ such that
$
|\Phi(  0)|
\les
B_{\Phi},
$ $
\mathbb P\text{-a.s.}
$

 \end{description}

For notational convenience, for $p\ges2$,
$r\ges0$, $0\les t<s\les\infty$, and
$k\in L_\lambda^2(E)$, set
\begin{equation}\label{eq:bsde-data}
\begin{aligned}
&\Psi_r(k)
:=
\int_E\psi(r,e,k(e))\,\lambda(\mathrm de),\\
&\zeta_p(t,s)
 :=
\Big(\int_t^s
|\psi(r,\cdot,0)|_{\lambda,1}^2\,\mathrm dr
\Big)^{\frac p2}
+
\Big(\int_t^s
|f(r,0,0,0,0)|\,\mathrm dr
\Big)^p.
\end{aligned}
\end{equation}
We first establish an a priori $L^p$-estimate for
solutions of \eqref{eq:infinite-horizon-bsde}.

\begin{proposition}\label{prop:bsde-apriori}
\label{prop:stability-infinite-bsde}\sl
Suppose that {\bf (A1)--(A3)} and
{\bf (B1)--(B4)} hold.
{Let $(Y,Z,K)$ be a solution of
\eqref{eq:infinite-horizon-bsde} with}
$
 (Y,Z,K)\in
\mathbb B_\lambda^p(t,\infty).
$
Then the following estimates hold.

\begin{enumerate}
\item[(i)]
There exists a constant $C>0$, depending only
on $p$, $C_\Phi$, $B_\Phi$, and the integrability
bounds in {\bf (A2)} and {\bf (B2)}, such that
\begin{equation}\label{eq:bsde-apriori}
\begin{aligned}
&\mathbb E_t\Big[
\sup_{r\ges t}|Y_r|^p
+\Big(\int_t^\infty|Z_r|^2\,\mathrm dr
\Big)^{\frac p2} +
\Big(\int_t^\infty
|K_r(\cdot)|_{\lambda,2}^2\,\mathrm dr
\Big)^{\frac p2}
{+\int_t^\infty\int_E |K_r(e)|^p
\lambda(\mathrm d e)\,\mathrm dr}
\Big]\\
& \les C\Big\{
1+|\xi|^p
+\mathbb E_t\big[
\eta_p(t,\infty)+\zeta_p(t,\infty)
\big]
\Big\},
\qquad\mathbb P\text{-a.s.}
\end{aligned}
\end{equation}

\item[(ii)]
Let
$\bar X\in{\mathcal S^p}(t,\infty;\dbR^n)$
have an almost sure limit $\bar X_\infty$.
Suppose that another coefficient triple
$(\bar\Phi,\bar f,\bar\psi)$ satisfies
{\bf (B1)--(B4)} with common structural constants
and integrability bounds.
Let $(\bar Y,\bar Z,\bar K)\in
\mathbb B_\lambda^p(t,\infty)$ be the corresponding
solution of  BSDE \eqref{eq:infinite-horizon-bsde} associated with $(\bar X,\bar\Phi,\bar f,\bar\psi)$.
Then there exists a constant $C>0$, depending
only on $p$ and the common structural constants
and integrability bounds, such that
\be\label{eq:stability-infinite-bsde}
\begin{aligned}
&\mathbb E_t\Big[
 \sup_{r\in[t,\infty)}|Y_r-\bar Y_r|^p
 +
 \(\int_t^\infty|Z_r-\bar Z_r|^2\,\mathrm d r\)^{\frac p2}
 +
 \(\int_t^\infty
         |K_r-\bar K_r(\cdot)|_{\lambda,2}^2\,\mathrm d r
  \)^{\frac p2}
 \\
&\qquad {+
\int_t^\infty\int_E
|K_r(e)-\bar K_r(e)|^p\lambda(\mathrm d e)\,\mathrm d r
\Big]}\\
&\les
 C\mathbb E_t\Big[
 |(\Phi -\bar\Phi)(\bar X_\infty)|^p
 +
 \sup_{r\in[t,\infty)}|X_r-\bar X_r|^p
  +
  \(\int_t^\infty|(f-\bar f)\big(r,\bar X_r,\bar Y_r,\bar Z_r,
       \bar\Psi_r(\bar K_r)\big) |\,\mathrm d r\)^p
  \\
& \qq\q{+
  \(\int_t^\infty\|l_k(r)\|_\infty
  |(\psi-\bar\psi)(r,\cd,\bar K_r(\cd))|_{\lambda,1}
  \,\mathrm d r\)^p}\Big],
\end{aligned}\ee
Here $\bar\Psi$ is defined analogously to $\Psi$ using the barred coefficients.

\end{enumerate}
\end{proposition}

\begin{proof}
\emph{Step 1. Conditional finite-horizon $L^p$ estimates.}
Fix $T>t$. The restriction of $(Y,Z,K)$ to $[t,T]$
solves the finite-horizon BSDE with terminal value
$Y_T$ and generator
$
F_r(y,z,k):=f(r,X_r,y,z,\Psi_r(k)).
$ Moreover, {\bf (B2)} yields the uniform bound
$$
\int_t^T\Big[
\|l_y(r)\|_\infty+\|l_z(r)\|_\infty^2
+\|l_k(r)\|_\infty^2
|\rho(r,\cdot)|_{\lambda,2}^2
\Big]\,\mathrm d r
\les L.
$$
Applying It\^o's formula to $|Y|^2$,
we obtain, for $s\in[t,T]$,
\begin{equation}\label{eq:bsde-energy}
\begin{aligned}
&\q |Y_s|^2
+\int_s^T|Z_r|^2\,\mathrm dr
+\int_s^T\int_E|K_r(e)|^2
\,\mu(\mathrm dr,\mathrm de)\\
&  =
|Y_T|^2
+2\int_s^T
Y_rF_r(Y_r,Z_r,K_r)\,\mathrm dr-
2\int_s^T Y_{r-}Z_r\,\mathrm dB_r
-2\int_s^T\int_E
Y_{r-}K_r(e)\,
\widetilde\mu(\mathrm dr,\mathrm de).
\end{aligned}
\end{equation}
All subsequent estimates are obtained under
standard localization, which is removed
by conditional Fatou's lemma.

Taking conditional expectations in
\eqref{eq:bsde-energy}, and using
the conditional BDG and Young inequalities,
we obtain the standard finite-horizon $L^2$
estimate
\begin{equation}\label{eq:bsde-L2-finite} \mathbb E_s\Big[
\sup_{r\in[s,T]}|Y_r|^2
+\int_s^T|Z_r|^2\,\mathrm dr
+\int_s^T|K_r(\cdot)|_{\lambda,2}^2
\,\mathrm dr
\Big] \les C\mathbb E_s\Big[
|Y_T|^2
+\sup_{r\in[s,T]}|X_r|^2
+\zeta_2(s,T)
\Big].
\end{equation}
Then,
applying the conditional Doob inequality with exponent
${\frac{p}{2}}>1$ to \eqref{eq:bsde-L2-finite}, we obtain
\begin{equation}\label{eq:Y-p-finite}
\mathbb E_t\Big[
\sup_{r\in[t,T]}|Y_r|^p
\Big]
\les C\mathbb E_t\Big[
|Y_T|^p
+\sup_{r\in[t,T]}|X_r|^p
+\zeta_p(t,T)
\Big].
\end{equation}
For $p=2$, this follows directly from
\eqref{eq:bsde-L2-finite}.

Returning to \eqref{eq:bsde-energy} with $s=t$,
we raise the resulting energy inequality to the
power ${\frac{p}{2}}$. By the conditional BDG and Young
inequalities, for every $\varepsilon>0$,
$$
\begin{aligned}
&\mathbb E_t\Big[
\Big|\int_t^T Y_{r-}Z_r\,\mathrm dB_r
\Big|^{\frac p2}
\Big] \les
\varepsilon\mathbb E_t\Big[
\Big(\int_t^T|Z_r|^2\,\mathrm dr
\Big)^{\frac p2}
\Big]
+C_\varepsilon\mathbb E_t\Big[
\sup_{r\in[t,T]}|Y_r|^p
\Big],\\[1mm]
&\mathbb E_t\Big[
\Big|\int_t^T\int_E
Y_{r-}K_r(e)\,
\widetilde\mu(\mathrm dr,\mathrm de)
\Big|^{\frac p2}
\Big] \les
\varepsilon\mathbb E_t\Big[
\Big(\int_t^T\int_E|K_r(e)|^2
\,\mu(\mathrm dr,\mathrm de)
\Big)^{\frac p2}
\Big]
+C_\varepsilon\mathbb E_t\Big[
\sup_{r\in[t,T]}|Y_r|^p
\Big].
\end{aligned}
$$
Combining \eqref{eq:poisson-compensator-moment}
with {\bf (B2)}
and \eqref{eq:bsde-energy}, and choosing
$\varepsilon>0$ sufficiently small,
we obtain, after the standard localization
and absorption arguments,
\begin{equation}\label{eq:ZK-p-finite}
\begin{aligned}
& \mathbb E_t\Big[
\Big(\int_t^T|Z_r|^2\,\mathrm dr\Big)^{\frac p2}
+\Big(\int_t^T
|K_r(\cdot)|_{\lambda,2}^2\,\mathrm dr
\Big)^{\frac p2} +
\Big(\int_t^T\int_E|K_r(e)|^2
\,\mu(\mathrm dr,\mathrm de)
\Big)^{\frac p2}
\Big]\\
& \les C\mathbb E_t\Big[
\sup_{r\in[t,T]}|Y_r|^p
+\sup_{r\in[t,T]}|X_r|^p
+\zeta_p(t,T)
\Big].
\end{aligned}
\end{equation}
The constant $C>0$ depends only on $p$
and the integrability bounds in {\bf (B2)}
and is independent of $T$.

\medskip

\emph{Step 2. Passage to the infinite horizon.}
The terminal condition in \eqref{eq:infinite-horizon-bsde}
and
$
Y\in\cS ^p(t,\infty;\mathbb R)
$
imply that
$
Y_T\to\Phi(X_\infty)
$
in $L^p(\Omega)$ and $\mathbb P$-a.s.\ as $T\to\infty$.
Consequently,
$
|Y_T|^p\to|\Phi(X_\infty)|^p
$
in $L^1(\Omega)$.
Combining \eqref{eq:Y-p-finite} and
\eqref{eq:ZK-p-finite} with \eqref{eq:poisson-p-moment},
then letting $T\to\infty$
by conditional dominated and monotone convergence,
we obtain
\begin{equation}\label{YZK-p}
\begin{aligned}
&\mathbb E_t\Big[
\sup_{r\in[t,\infty)}|Y_r|^p
+\Big(\int_t^\infty|Z_r|^2\,\mathrm d r\Big)^{\frac p2}
+\Big(\int_t^\infty
|K_r(\cdot)|_{\lambda,2}^2\,\mathrm d r\Big)^{\frac p2}
{+\int_t^\infty\int_E|K_r(e)|^p
\lambda(\mathrm d e)\,\mathrm d r}
\Big]\\
&\les C\mathbb E_t\Big[
|\Phi(X_\infty)|^p
+\sup_{r\in[t,\infty)}|X_r|^p
+\zeta_p(t,\infty)
\Big] \les C\Big\{
1+|\xi|^p
+\mathbb E_t\big[
\eta_p(t,\infty)
+\zeta_p(t,\infty)
\big]\Big\},
\end{aligned}
\end{equation}
where the last inequality follows from
{\bf (B2)}, {\bf (B4)}, and the state estimate
\eqref{eq:state-estimate}.
This proves \eqref{eq:bsde-apriori} for  $p\ges 2$.

 \emph{Step 3. {Continuous dependence}.}
Set
$
\widehat X:=X-\bar X,$ $
\widehat Y:=Y-\bar Y,$ $
\widehat Z:=Z-\bar Z,$ $
\widehat K:=K-\bar K.
$
Define
$$
\widehat f(r,x,y,z,v)
:=
f\bigl(
r,x+\bar X_r,y+\bar Y_r,z+\bar Z_r,
v+\bar\Psi_r(\bar K_r)
\bigr) -
\bar f\bigl(
r,\bar X_r,\bar Y_r,\bar Z_r,
\bar\Psi_r(\bar K_r)
\bigr),
$$
 and
$
\widehat\psi(r,e,k)
:=
\psi\bigl(r,e,k+\bar K_r(e)\bigr)
-\bar\psi\bigl(r,e,\bar K_r(e)\bigr),
$ $
\widehat\Psi_r(k)
:=
\int_E\widehat\psi(r,e,k(e))\,\lambda(\mathrm d e).
$
The coefficients $\widehat f$ and $\widehat\psi$
satisfy the Lipschitz conditions in {\bf (B2)}
with the same moduli as $f$ and $\psi$, respectively.
For convenience,  set
$$
\begin{aligned}
&\widehat f_r^0
 :=\widehat f(r,0,0,0,0)
=(f-\bar f)\bigl(
r,\bar X_r,\bar Y_r,\bar Z_r,
\bar\Psi_r(\bar K_r)
\bigr),\\
&\widehat\psi_r^0(e)
 :=\widehat\psi(r,e,0)
=(\psi-\bar\psi)\bigl(r,e,\bar K_r(e)\bigr).
\end{aligned}
$$
Note that
$
\widehat\Psi_r(0)
=\int_E\widehat\psi_r^0(e)\,\lambda(\mathrm d e),
$
so that
$
|\widehat f(r,0,0,0,\widehat\Psi_r(0))|
\les
|\widehat f_r^0|
+l_k(r)|\widehat\psi_r^0|_{\lambda,1}.
$
Then, $(\widehat Y,\widehat Z,\widehat K)$
satisfies
\begin{equation}\label{eq:bsde-difference}
\widehat Y_s
= \widehat Y_T
+\int_s^T
\widehat f\bigl(
r,\widehat X_r,\widehat Y_r,\widehat Z_r,
\widehat\Psi_r(\widehat K_r)
\bigr)\,\mathrm d r -\int_s^T\widehat Z_r\,\mathrm d B_r
-\int_s^T\int_E
 \widehat K_r(e)\,\widetilde\mu(\mathrm d r,\mathrm d e),
\quad t\les s\les T<\infty,
\end{equation}
with $\widehat Y_T\longrightarrow
\Phi(X_\infty)-\bar\Phi(\bar X_\infty)$,
$\mathbb P$-a.s., as $T\to\infty$.

Applying the conditional finite-horizon
$L^p$ estimates from Step 1 to
\eqref{eq:bsde-difference}, and then passing
to the infinite horizon as in Step 2,
we obtain
\begin{equation}\label{eq:preliminary-stability-estimate}
\begin{aligned}
&\mathbb E_t\Big[
\sup_{r\ges t}|\widehat Y_r|^p
+\Big(\int_t^\infty|\widehat Z_r|^2
\,\mathrm dr\Big)^{\frac p2}+
\Big(\int_t^\infty
|\widehat K_r(\cdot)|_{\lambda,2}^2
\,\mathrm dr\Big)^{\frac p2}
{+\int_t^\infty\int_E|\widehat K_r(e)|^p
\lambda(\mathrm d e)\,\mathrm dr}
\Big]\\
&  \les C\mathbb E_t\Big[
|\Phi(X_\infty)-\bar\Phi(\bar X_\infty)|^p
+\sup_{r\ges t}|\widehat X_r|^p +
\Big(\int_t^\infty
|\widehat f_r^0|\,\mathrm dr\Big)^p +
\Big(\int_t^\infty
\|l_k(r)\|_\infty
|\widehat\psi_r^0|_{\lambda,1}
\,\mathrm dr\Big)^p
\Big].
\end{aligned}
\end{equation}
Finally, {using the Lipschitz continuity of
$\Phi$ to bound the terminal difference in
\eqref{eq:preliminary-stability-estimate},} we obtain
(ii)
from \eqref{eq:preliminary-stability-estimate}.
\end{proof}

%
%

Based on the preceding a priori
$L^p$ estimates, we establish
the well-posedness of the
infinite-horizon BSDE with jumps
\eqref{eq:infinite-horizon-bsde}.

\begin{theorem}
\label{thm:existence-infinite-bsde}\sl
Suppose that {\bf (A1)--(A3)}
and {\bf (B1)--(B4)} hold.
Then, for every $t\ges0$ and
$\xi\in L_{\mathcal F_t}^p(\Omega;\mathbb R^n)$,
the infinite-horizon BSDE with jumps
\eqref{eq:infinite-horizon-bsde}
admits a unique solution
$(Y,Z,K)\in\mathbb B_\lambda^p(t,\infty).$
{Moreover,} the solution satisfies
$
Y_T\to\Phi(X_\infty)
$
in $L^p(\Omega;\mathbb R)$
as $T\to\infty$.
\end{theorem}

\begin{proof}
{Uniqueness follows directly from}
Proposition~\ref{prop:stability-infinite-bsde}.
{We now prove existence.}

For $s\ges t$, set
$\eta_s:=\mathbb E_s[\Phi(X_\infty)]$.
By Lemma 3.1 and {\bf(B2)}, {\bf(B4)},
$\Phi(X_\infty)\in L^p_{\mathcal F_\infty}(\Omega)$.
Thus $\eta\in{\mathcal S^p}(t,\infty;\dbR)$
is a uniformly integrable {c\`adl\`ag} martingale.
Moreover, the martingale convergence theorem yields
$\eta_s\to\Phi(X_\infty)$ in $L^p(\Omega)$
and $\mathbb P$-a.s. as $s\to\infty$.
The martingale representation theorem then yields
square-integrable predictable processes
$(\zeta,\kappa)$ such that
$$
\eta_s=\eta_t
+\int_t^s\zeta_r\,\mathrm dB_r
+\int_t^s\int_E
\kappa_r(e)\,\widetilde\mu(\mathrm dr,\mathrm de),
\qquad s\ges t.
$$
We further show that
$
(\zeta,\kappa)\in
\mathcal H^{p,2}(t,\infty;\dbR^d)
\times
\bigl(
\mathcal K_{\lambda,2}^p(t,\infty;\dbR)
\cap
\mathcal K_{\lambda,p}^p(t,\infty;\dbR)
\bigr).
$
To this end, set $M_s:=\eta_s-\eta_t$ for $s\ges t$.
By Doob's maximal inequality and the
BDG inequality,
\begin{equation}\label{eq:martingale-jump-p-moment}
\begin{aligned}
&\mathbb E\[
\(
\int_t^\infty|\zeta_r|^2\,\mathrm dr
+\int_t^\infty\int_E
|\kappa_r(e)|^2\,\mu(\mathrm dr,\mathrm de)
\)^{{\frac{p}{2}}}
\]\\
&
=\mathbb E\bigl[[M]_\infty^{{\frac{p}{2}}}\bigr]
\les C_p\mathbb E\Big[\sup_{s\ges t}|M_s|^p\Big]
\les C_p\mathbb E|\Phi(X_\infty)|^p
<\infty.
\end{aligned}
\end{equation}
Applying Lemma~\ref{lem:poisson-moment-estimates}
with $K=\kappa$ and using
\eqref{eq:martingale-jump-p-moment}, we obtain
$$
 \mathbb E\(
\int_t^\infty|\kappa_r(\cdot)|_{\lambda,2}^2
\,\mathrm dr\)^{{\frac{p}{2}}} +
\mathbb E\int_t^\infty
|\kappa_r(\cdot)|_{\lambda,p}^p\,\mathrm dr
<\infty.
$$
Together with \eqref{eq:martingale-jump-p-moment},
this establishes the claimed $L^p$-integrability
of $(\zeta,\kappa)$.

Define
$
F(r,y,z,k):=f(r,X_r,y,z,\Psi_r(k)),
$
for $k\in L_\lambda^2(E)$.
For each $T>t$, consider the finite-horizon
BSDE with jumps
\begin{equation}\label{eq:thm31-truncated}
Y_s^T =\eta_T
+\int_s^T F(r,Y_r^T,Z_r^T,K_r^T)\,\mathrm d r -\int_s^T Z_r^T\,\mathrm d B_r
-\int_s^T\int_E K_r^T(e)\,
\widetilde\mu(\mathrm d r,\mathrm d e),
\quad t\les s\les T.
\end{equation}
By {\bf (A1)--(A3)} and {\bf (B1)--(B3)},
the generator $F$ satisfies the time-dependent
Lipschitz conditions with deterministic integrable
moduli. Moreover,
$\ds
\mathbb E\Big[
\Big(\int_t^\infty
|F(r,0,0,0)|\,\mathrm dr\Big)^p
\Big]<\infty.
$
The finite-horizon existence and uniqueness result  for BSDEs with jumps
(see
\cite[Chapter 1, Lemma 6 and Theorem 7]{Situ2000})
yields a unique solution
$
(Y^T,Z^T,K^T)\in
{\mathcal S^2}(t,T;\dbR)
\times\mathcal H^{2,2}(t,T;\dbR^d)
\times\mathcal K_{\lambda,2}^2(t,T;\dbR).
$  Moreover,
$$
|Y_s^T|^2
\les C\mathbb E_s\Big[
|\eta_T|^2
+\Big(\int_t^T|F(r,0,0,0)|\,\mathrm dr\Big)^2
\Big],
\qquad s\in[t,T],
$$
where $C$ is independent of $T$.

\medskip

 \textit{Step 1. Uniform estimates.}
Arguing as in the proof of
Proposition~\ref{prop:bsde-apriori},
and using {\bf (B2)} together with the
uniform $L^p$-bounds on $\eta_T$ and
$F(\cdot,0,0,0)$, we obtain
\begin{equation}\label{eq:thm31-uniform}
\sup_{T>t}\mathbb E\Big[
 \sup_{s\in[t,T]}|Y_s^T|^p
+\Big(\int_t^T|Z_s^T|^2\,\mathrm ds\Big)^{\frac p2}+\Big(\int_t^T
|K_s^T(\cdot)|_{\lambda,2}^2\,\mathrm ds
\Big)^{\frac p2} +\Big(\int_t^T\int_E
|K_s^T(e)|^2\,\mu(\mathrm ds,\mathrm de)
\Big)^{\frac p2}
\Big]<\infty.
\end{equation}

\medskip
 \textit{Step 2. Tail estimate.}
Set
$ \ds
\varepsilon_m:=\mathbb E\Big[
\Big(\int_m^\infty|F(r,\eta_r,\zeta_r,\kappa_r)|\,\mathrm dr\Big)^p
\Big]
$, for $m>t$.
By \textbf{(B2)}, H\"older's inequality, and the preceding estimates, we have
$$
\mathbb E\Big[
\Big(\int_t^\infty|F(r,\eta_r,\zeta_r,\kappa_r)|\,\mathrm dr\Big)^p
\Big]<\infty,
$$
which  yields
$\varepsilon_m\to0$, as $m\to\infty$.

Fix $T>m$. On $[m,T]$, set
$
\widehat Y_s:=Y_s^T-\eta_s,
 $ $
\widehat Z_s:=Z_s^T-\zeta_s,
 $ $
\widehat K_s:=K_s^T-\kappa_s.
$
Then $(\widehat Y,\widehat Z,\widehat K)$ satisfies
the following BSDE on $[m,T]$:
\begin{equation}\label{eq:thm31-shifted}
\widehat Y_s= \int_s^T
\widehat F(r,\widehat Y_r,\widehat Z_r,\widehat K_r)
\,\mathrm dr -\int_s^T\widehat Z_r\,\mathrm dB_r
-\int_s^T\int_E\widehat K_r(e)\,
\widetilde\mu(\mathrm dr,\mathrm de),
\qquad m\les s\les T.
\end{equation}
where
$
\widehat F(r,y,z,k)
:=F(r,\eta_r+y,\zeta_r+z,\kappa_r+k).
$
The generator $\widehat F$ inherits the Lipschitz moduli of $F$,
whose integral bound is independent of $m$ and $T$.
{The finite-horizon $L^p$ estimates
\eqref{eq:Y-p-finite} and \eqref{eq:ZK-p-finite} from
Proposition~\ref{prop:bsde-apriori}, together with
\eqref{eq:poisson-p-moment}, give for \eqref{eq:thm31-shifted}}
\begin{equation}\label{eq:thm31-tail}
\begin{aligned}
&\mathbb E\Big[
 \sup_{s\in[m,T]}|Y_s^T-\eta_s|^p
 +\Big(\int_m^T|Z_s^T-\zeta_s|^2\,\mathrm ds\Big)^{\frac p2} +\Big(\int_m^T
|K_s^T(\cdot)-\kappa_s(\cdot)|_{\lambda,2}^2
\,\mathrm ds\Big)^{\frac p2}\\
&\quad+{\int_m^T\int_E
 |K_s^T(e)-\kappa_s(e)|^p\,\lambda(\mathrm de)\,\mathrm ds}
\Big] \les C\mathbb E\Big[
\Big(\int_m^T|F(r,\eta_r,\zeta_r,\kappa_r)|\,\mathrm dr\Big)^p
\Big] \les C\varepsilon_m,
\end{aligned}
\end{equation}
where $C$ is independent of $m$ and $T$.

\medskip

 \emph{Step 3. The Cauchy property.}
For each $T>t$, extend the finite-horizon solution
to $[t,\infty)$ by setting
$
\widetilde Y_s^T:=Y_s^T\mathbf1_{\{s\les T\}}
+\eta_s\mathbf1_{\{s>T\}} $
and
$
\widetilde Z_s^T
 :=Z_s^T\mathbf1_{\{s\les T\}}
+\zeta_s\mathbf1_{\{s>T\}},$ $
\widetilde K_s^T
 :=K_s^T\mathbf1_{\{s\les T\}}
+\kappa_s\mathbf1_{\{s>T\}}.
$
Then $\widetilde Y^T$ is $\dbF$-adapted and c\`adl\`ag.
{Moreover, $\widetilde K^T\in
\mathcal K_{\lambda,p}^p(t,\infty)$: on $[t,T]$ this follows
from \eqref{eq:thm31-uniform} and
\eqref{eq:poisson-p-moment}; after $T$, it equals $\kappa$,
whose integrability follows from
\eqref{eq:martingale-jump-p-moment} and
Lemma~\ref{lem:poisson-moment-estimates}.}

Fix $m>t$. For every $T>m$, \eqref{eq:thm31-tail} gives
\begin{equation}\label{eq:thm31-extended-tail}
\begin{aligned}
&\mathbb E\Big[
\sup_{s\ges m}|\widetilde Y_s^T-\eta_s|^p
+\Big(\int_m^\infty
|\widetilde Z_s^T-\zeta_s|^2
\,\mathrm d s\Big)^{\frac p2} +\Big(\int_m^\infty
|\widetilde K_s^T(\cdot)-\kappa_s(\cdot)|_{\lambda,2}^2
\,\mathrm d s\Big)^{\frac p2}\\
&\quad+{\int_m^\infty\int_E
  |\widetilde K_s^T(e)-\kappa_s(e)|^p\,
  \lambda(\mathrm de)\,\mathrm ds}
\Big]
\les C\varepsilon_m,
\end{aligned}
\end{equation}
where $C$ is independent of $m$ and $T$.

Let $T,U>m$. On $[t,m]$, the solutions
$(Y^T,Z^T,K^T)$ and
$(Y^U,Z^U,K^U)$
satisfy the same BSDE, with generator $F$
and terminal values $Y_m^T$ and $Y_m^U$,
respectively.
Set
$$
\delta Y:=Y^T-Y^U,\quad
\delta Z:=Z^T-Z^U,\quad
\delta K:=K^T-K^U
$$
and
$
\delta F_r
:=F(r,Y_r^T,Z_r^T,K_r^T)
-F(r,Y_r^U,Z_r^U,K_r^U).
$
Then,
$$
\delta Y_s
= \delta Y_m+\int_s^m\delta F_r\,\mathrm d r
-\int_s^m\delta Z_r\,\mathrm d B_r -\int_s^m\int_E
\delta K_r(e)\,\widetilde\mu(\mathrm d r,\mathrm d e),\q s\in[t,m].
$$
Repeating the   $L^p$ arguments
in the proof of Proposition~\ref{prop:bsde-apriori},
we obtain
\begin{equation}\label{eq:thm31-local-stability}
 \begin{aligned}
 &\mathbb E_t\Big[
\sup_{s\in[t,m]}|\delta Y_s|^p
+\Big(\int_t^m|\delta Z_r|^2\,\mathrm d r\Big)^{\frac p2} +
\Big(\int_t^m
|\delta K_r(\cdot)|_{\lambda,2}^2
\,\mathrm d r\Big)^{\frac p2}\\
&\qq +\Big(\int_t^m\int_E
|\delta K_s(e)|^2\,\mu(\mathrm ds,\mathrm de)\Big)^{\frac p2}
\Big]
\les C\mathbb E_t[|\delta Y_m|^p],
\end{aligned}\end{equation}
where $C>0$ is independent of $m$, $T$, and $U$.

Taking expectations in \eqref{eq:thm31-local-stability} and using
{\eqref{eq:poisson-p-moment} for $\delta K$ and}
\eqref{eq:thm31-tail}, we obtain
\begin{equation}\label{eq:thm31-local-stability-1}
\begin{aligned}
&\mathbb E\Big[
 \sup_{s\in[t,m]}|\delta Y_s|^p
+\Big(\int_t^m|\delta Z_s|^2\,\mathrm ds\Big)^{\frac p2}
+\Big(\int_t^m
 |\delta K_s(\cdot)|_{\lambda,2}^2\,\mathrm ds\Big)^{\frac p2}
 +{\int_t^m\int_E
 |\delta K_s(e)|^p\,\lambda(\mathrm de)\,\mathrm ds}
\Big]\\
&
\les C\mathbb E[|\delta Y_m|^p]\les 2^{p-1}C\mathbb E\Big[
|Y_m^T-\eta_m|^p+|Y_m^U-\eta_m|^p
\Big]
\les C\varepsilon_m,
\end{aligned}
\end{equation}
where $C$ is independent of $m,T,U$.

On $[m,\infty)$, the triangle inequality and
\eqref{eq:thm31-extended-tail} yield
$$
\begin{aligned}
&\mathbb E\Big[
 \sup_{s\ges m}
|\widetilde Y_s^T-\widetilde Y_s^U|^p
+\Big(\int_m^\infty
|\widetilde Z_s^T-\widetilde Z_s^U|^2\,\mathrm d s\Big)^{\frac p2} +\Big(\int_m^\infty
|\widetilde K_s^T(\cdot)-\widetilde K_s^U(\cdot)|_{\lambda,2}^2
\,\mathrm d s\Big)^{\frac p2}\\
&\quad+{\int_m^\infty\int_E
  |\widetilde K_s^T(e)-\widetilde K_s^U(e)|^p\,
  \lambda(\mathrm de)\,\mathrm ds}
\Big]
\les C\varepsilon_m.
\end{aligned}
$$
Combining the estimates on $[t,m]$ and $[m,\infty)$,
we conclude that, for $T,U>m$,
\begin{equation}\label{eq:thm31-cauchy}
\begin{aligned}
&\mathbb E\Big[
\sup_{s\ges t}
|\widetilde Y_s^T-\widetilde Y_s^U|^p
+\Big(\int_t^\infty
|\widetilde Z_s^T-\widetilde Z_s^U|^2
\,\mathrm d s\Big)^{\frac p2}+\Big(\int_t^\infty
|\widetilde K_s^T(\cdot)-\widetilde K_s^U(\cdot)|_{\lambda,2}^2
\,\mathrm d s\Big)^{\frac p2}\\
&\quad+{\int_t^\infty\int_E
  |\widetilde K_s^T(e)-\widetilde K_s^U(e)|^p\,
  \lambda(\mathrm de)\,\mathrm ds}
\Big]
\les C\varepsilon_m.
\end{aligned}
\end{equation}

Due to  $\varepsilon_m\to0$,
$(\widetilde Y^T,\widetilde Z^T,\widetilde K^T)$ forms a Cauchy family
in $
\mathbb B_\lambda^p(t,\infty).
$
By completeness, there exists $(Y,Z,K)$ in this
product space such that
\begin{equation}\label{eq:thm31-convergence}
\begin{aligned}
 &\mathbb E\Big[
 \sup_{s\in[t,\infty)}
|\widetilde Y_s^T-Y_s|^p
+\Big(\int_t^\infty
|\widetilde Z_s^T-Z_s|^2\,\mathrm ds
\Big)^{\frac p2}\\
 &\quad+\Big(\int_t^\infty
|\widetilde K_s^T(\cdot)-K_s(\cdot)|_{\lambda,2}^2
\,\mathrm ds\Big)^{\frac p2}
+{\int_t^\infty\int_E
 |\widetilde K_s^T(e)-K_s(e)|^p\,\lambda(\mathrm de)\,\mathrm ds}
\Big]\longrightarrow0\quad\text{as }T\to\infty.
\end{aligned}
\end{equation}
Furthermore, a subsequence
$\widetilde Y^{T_n}$ converges to $Y$ uniformly on
$[t,\infty)$ almost surely.

\medskip
 \emph{Step 4. Passage to the limit.}
Fix $t<R<\infty$. For $T>R$,
\eqref{eq:thm31-truncated} gives
$$
Y_s^T= Y_R^T
+\int_s^R
f\bigl(r,X_r,Y_r^T,Z_r^T,\Psi_r(K_r^T)\bigr)
\,\mathrm dr -\int_s^R Z_r^T\,\mathrm dB_r
-\int_s^R\int_E K_r^T(e)\,
\widetilde\mu(\mathrm dr,\mathrm de),
\quad s\in[t,R].
$$
By {\bf (B2)},  the It\^o isometry,
Doob's inequality and \eqref{eq:thm31-convergence},
taking $T\to\i$ in
${\mathcal S^2}(t,R;\mathbb R)$, we obtain
\begin{equation}\label{eq:thm31-limit-equation}
Y_s=  Y_R+\int_s^R
f\bigl(r,X_r,Y_r,Z_r,\Psi_r(K_r)\bigr)
\,\mathrm dr
 -\int_s^R Z_r\,\mathrm dB_r
-\int_s^R\int_E K_r(e)\,
\widetilde\mu(\mathrm dr,\mathrm de),
\quad t\les s\les R<\infty.
\end{equation}

\medskip
 \emph{Step 5. The terminal condition.}
Letting $T\to\infty$ in \eqref{eq:thm31-extended-tail}
and using \eqref{eq:thm31-convergence} and Fatou's lemma,
we obtain
\begin{equation}\label{eq:thm31-terminal-tail}
\mathbb E\Big[
\sup_{s\ges m}|Y_s-\eta_s|^p
\Big]\les C\varepsilon_m.
\end{equation}
Then the nonincreasing random
variables
$\ds
 \sup_{s\ges m}|Y_s-\eta_s|^p
$
converge to zero almost surely, and  $Y_s-\eta_s\to0$ in $L^p(\Omega)$.
Together with the convergence of $\eta_s$ to
$\Phi(X_\infty)$ in $L^p(\Omega)$ and almost surely,
this yields
$$
Y_s\longrightarrow\Phi(X_\infty)
\quad\text{in }L^p(\Omega)\text{ and }\mathbb P\text{-a.s.}
\quad\text{as }s\to\infty.
$$
  \end{proof}

By the conditional stability estimate
\eqref{eq:stability-infinite-bsde}
and the state estimate
\eqref{eq:state-estimate},
we obtain the following conditional
$L^p$ Lipschitz estimate with respect
to the initial state.

\begin{corollary}[Lipschitz dependence on the state]
\label{prop:bsde-Y-lipschitz}\sl
Suppose that {\bf (A1)--(A3)}
and {\bf (B1)--(B4)} hold.
For $t\ges0$ and
$\xi_i\in L_{\mathcal F_t}^p(\Omega;\dbR^n)$,
$i=1,2$, let $X^i$ be the solution of
\eqref{eq:infinite-horizon-sde} with initial value
$\xi_i$, and let $(Y^i,Z^i,K^i)$ be the
corresponding solution of
\eqref{eq:infinite-horizon-bsde}.
Then there exists a constant $C>0$, depending
only on $p$, $C_\Phi$, and the integrability bounds
in {\bf (A2)} and {\bf (B2)}, such that
\begin{equation}\label{eq:bsde-Y-initial-Lipschitz}
\mathbb E_t\Big[
\sup_{s\ges t}|Y_s^1-Y_s^2|^p
\Big]
\les C|\xi_1-\xi_2|^p,
\qquad\mathbb P\text{-a.s.}
\end{equation}
\end{corollary}

To establish the comparison theorem for
\eqref{eq:infinite-horizon-bsde},
we impose the following additional condition.
\begin{description}

\item[(B5)]
 For every
$(x,y,z)\in\dbR^n\times\dbR\times\dbR^d$,
the mapping
$
v\longmapsto f(s,x,y,z,v)
$
is nondecreasing,
$\mathrm d\mathbb P\otimes\mathrm ds$-a.e.
Moreover, for all $k,k'\in\dbR$ with $k\ne k'$,
$$
\|l_k(s)\|_\infty
\frac{\psi(s,e,k)-\psi(s,e,k')}{k-k'}
\ges-1,\q \mathrm d\mathbb P\otimes\mathrm ds
\otimes\lambda(\mathrm de)\mbox{-a.s.}
$$

\end{description}

\begin{theorem}[Comparison theorem]
\label{thm:comparison-infinite-bsde}\sl
Suppose that {\bf (A1)--(A3)} hold.
For $i=1,2$,
and the same state process $X$, let
$(Y^i,Z^i,K^i)$ be the unique solution of
\eqref{eq:infinite-horizon-bsde}
corresponding to $(\Phi_i,f_i,\psi_i)$ satisfying
{\bf (B1)--(B4)}. Suppose additionally that
$(f_1,\psi_1)$ satisfies {\bf (B5)}.
{Write $l_y^1,l_z^1,l_k^1$ and $\rho_1$ for the
Lipschitz moduli of the first data set $(f_1,\psi_1)$.}
Assume that
$
\Phi_1(X_\infty)\ges\Phi_2(X_\infty)
$
$\mathbb P$-a.s., and
$
f_1(r,{X_r},Y_r^2,Z_r^2,\Psi^1(K_r^2))
\ges f_2(r,{X_r},Y_r^2,Z_r^2,\Psi^2(K_r^2)),
$ $
\mathrm d\mathbb P\otimes\mathrm dr\text{-a.e.}
$
Then $
Y_s^1\ges Y_s^2,$
$s\ges t$, $\mathbb P\text{-a.s.}$

\end{theorem}

\begin{proof}
Set
$
\Delta Y:=Y^1-Y^2,$ $
\Delta Z:=Z^1-Z^2,$ $
\Delta K:=K^1-K^2,
$
and
$$
h_r:=
f_1(r,X_r,Y_r^2,Z_r^2,\Psi^1(K_r^2))
-f_2(r,X_r,Y_r^2,Z_r^2,\Psi^2(K_r^2))
\ges0.
$$
By the standard difference-quotient linearization,
there exist progressively measurable processes
$\alpha,\beta$ and a predictable function $\theta$
such that, for $t\les s\les T<\infty$,
\begin{equation}\label{eq:comparison-linearized}
\begin{aligned}
\Delta Y_s
={}&\Delta Y_T
+\int_s^T\Big(
\alpha_r\Delta Y_r
+\langle\beta_r,\Delta Z_r\rangle +
\int_E\theta_r(e)\Delta K_r(e)
\,\lambda(\mathrm de)+h_r
\Big)\,\mathrm dr\\
&-\int_s^T\Delta Z_r\,\mathrm dB_r
-\int_s^T\int_E
\Delta K_r(e)\,\widetilde\mu(\mathrm dr,\mathrm de).
\end{aligned}
\end{equation}
Moreover, by {\bf (B2)} and {\bf (B5)}, these coefficients satisfy
$
|\alpha_r| \les l_y^1(r),$ $
|\beta_r| \les l_z^1(r),$   $\theta_r(e)=c_r\nu_r(e) $ with  $0\les c_r\les l_k^1(r)$ and
$$
\nu_r(e):=
\begin{cases}
\displaystyle
\frac{\psi_1(r,e,K_r^1(e))-\psi_1(r,e,K_r^2(e))}
     {K_r^1(e)-K_r^2(e)},
& K_r^1(e)\ne K_r^2(e),\\[2ex]
0,
& K_r^1(e)=K_r^2(e).
\end{cases}
$$
The monotonicity of $f_1$ and {\bf (B2), (B5)} imply
$
 \theta_r(e)\ges -1,
$ $
|\theta_r(e)|
\les {\|l_k^1(r)\|_\infty\rho_1(r,e)},
$
and hence
$$
\int_t^\infty
\Big(
\|\alpha_r\|_\infty
+\|\beta_r\|_\infty^2
+\|\theta_r(\cdot)\|_{\lambda,2}^2
\Big)\,\mathrm dr
\les L,
$$
where $L$ is a deterministic constant
depending only on the integrability bounds
in {\bf (B2)}.

For each $s\ges t$, let $\Gamma_{s,r}$, solve
\begin{equation}\label{eq:comparison-dual}
\mathrm d\Gamma_{s,r}
=\Gamma_{s,r-}\Big(
\alpha_r\,\mathrm dr
+\langle\beta_r,\mathrm dB_r\rangle
+\int_E\theta_r(e)\,
\widetilde\mu(\mathrm dr,\mathrm de)
\Big),\q r\ges s,\q
\Gamma_{s,s}=1.
\end{equation}
Since $\theta_r(e)\ges-1$, the stochastic
exponential representation implies
$\Gamma_{s,r}\ges0$.
{Moreover, {\bf (B2)} and the standard
$L^2$ estimate for linear SDEs with jumps yield}
\begin{equation}\label{eq:comparison-dual-L2}
\sup_{T\ges s}
\mathbb E_s\Big[
\sup_{r\in[s,T]}|\Gamma_{s,r}|^2
\Big]\les C,
\end{equation}
where $C$ is independent of $s$ and $T$.

Applying It\^o's product formula to
$\Gamma_{s,r}\Delta Y_r$, and using the
standard localization argument, we obtain
\begin{equation}\label{eq:comparison-dual-representation}
\Delta Y_s
=
\mathbb E_s\Big[
\Gamma_{s,T}\Delta Y_T
+\int_s^T\Gamma_{s,r}h_r\,\mathrm dr
\Big]\ges
-\mathbb E_s\big[
\Gamma_{s,T}(\Delta Y_T)^-
\big].
\end{equation}
On the other hand, by
Theorem~\ref{thm:existence-infinite-bsde},
$$
\Delta Y_T\to
\Delta\Phi_\infty
:=
\Phi_1(X_\infty)-\Phi_2(X_\infty)
\quad\text{in }L^p(\Omega),
$$
where $\Delta\Phi_\infty\ges0$.
Since $p\ges2$, it follows that
$\ds
\mathbb E\big[
|(\Delta Y_T)^-|^2
\big]\to0.
$
Consequently, by
\eqref{eq:comparison-dual-L2},
$$
0\les\mathbb E[(\Delta Y_s)^-]
 \les
\mathbb E\big[
\Gamma_{s,T}(\Delta Y_T)^-
\big] \les
C\big(
\mathbb E[|(\Delta Y_T)^-|^2]
\big)^{{\frac{1}{2}}}
\to0.
$$
Thus $\Delta Y_s\ges0$, $\mathbb P$-a.s.,
for each fixed  $s\ges t$. By the c\`adl\`ag property, the inequality
holds   for all $s\ges t$.
\end{proof}

\section{Infinite-horizon non-autonomous stochastic recursive differential games
with jumps}
\label{sec:stochastic-game}

We now apply the preceding SDE and BSDE results
to study an infinite-horizon two-player zero-sum
stochastic recursive differential game with jumps
and establish the connection between its value
functions and the associated HJBI equations.
 For this purpose, it suffices to work
in the square-integrable framework, and we henceforth
fix $p=2$.

Let the following deterministic Borel measurable
coefficients be given:
$$
\begin{aligned}
b&:[0,\infty)\times\mathbb R^n\times U\times V
\longrightarrow\mathbb R^n,\q
\sigma:[0,\infty)\times\mathbb R^n\times U\times V
\longrightarrow\mathbb R^{n\times d},\\
\gamma&:[0,\infty)\times E\times\mathbb R^n\times U\times V
\longrightarrow\mathbb R^n,
\end{aligned}
$$
and
$$
f :[0,\infty)\times\mathbb R^n\times\mathbb R
\times\mathbb R^d\times\mathbb R\times U\times V
\longrightarrow\mathbb R,\q
\rho :[0,\infty)\times E\longrightarrow[0,\infty),\q
\Phi:\mathbb R^n\longrightarrow\mathbb R,
$$
where
$U\subset\dbR^{m_1}$ and $V\subset\dbR^{m_2}$
are  nonempty compact sets.
We impose the following assumptions.

\begin{description}

\item[(G1)] The functions $b(r,x,u,v)$ and
$\sigma(r,x,u,v)$ are continuous in $(r,x,u,v)$.
For every $e\in E$,   $\gamma(r,e,x,u,v)$
is continuous in $(r,x,u,v)$.
Moreover, there exist   nonnegative measurable functions
$\ell_b,\ell_\sigma$ on $[0,\infty)$ and
$\ell_\gamma$ on $[0,\infty)\times E$ satisfying
$$
{\int_0^\infty\Big[
\ell_b(r)+\ell_\sigma(r)^2
+|\ell_\gamma(r,\cdot)|_{\lambda,2}^2
\Big]\,\mathrm dr<\infty},
$$
such that, for all $r\ges0$, $x,x'\in\mathbb R^n$, $e\in E$
and $(u,v)\in U\times V$,
$$
\begin{aligned}
&|b(r,x,u,v)-b(r,x',u,v)|
\les\ell_b(r)|x-x'|,\q
 |\sigma(r,x,u,v)-\sigma(r,x',u,v)|
\les\ell_\sigma(r)|x-x'|,\\
&|\gamma(r,e,x,u,v)-\gamma(r,e,x',u,v)|
\les\ell_\gamma(r,e)|x-x'|.
\end{aligned}
$$

\item[(G2)]
The function $f$ is jointly continuous in
$(r,x,y,z,k,u,v)$, and
$\rho:[0,\infty)\times E\to[0,\infty)$
is Borel measurable, with
$\rho(r,\cdot)\in L_\lambda^2(E)$
for a.e.\ $r\ges0$.
There exist   nonnegative measurable functions
$\ell_x,\ell_y,\ell_z,\ell_k$ on $[0,\infty)$, with
$$
\int_0^\infty\Big[
\ell_x(r)+\ell_y(r)+\ell_z(r)^2
+\ell_k(r)^2+\ell_k(r)^2
 |\rho(r,\cdot)|_{\lambda,2}^2
\Big]\,\mathrm dr<\infty,
$$
 such that, for all $r\ges0$, $x,x'\in\mathbb R^n$,
$y,y',k,k'\in\mathbb R$, $z,z'\in\mathbb R^d$,
and $(u,v)\in U\times V$,
$$
 |f(r,x,y,z,k,u,v)-f(r,x',y',z',k',u,v)| \les
\ell_x(r)|x-x'|+\ell_y(r)|y-y'|
+\ell_z(r)|z-z'|+\ell_k(r)|k-k'|.
$$

\item[(G3)]
The terminal function $\Phi$ is Lipschitz continuous in $x$
with constant $C_\Phi$.

\item[(G4)]
For every $(r,x,y,z,u,v)$, the mapping
$k\mapsto f(r,x,y,z,k,u,v)$ is nondecreasing.

\item[(G5)] The coefficients at the origin satisfy
$ \ds
{\int_0^\infty
\big[
b_0(r)+\sigma_0(r)^2
+\gamma_0(r)^2
+f_0(r)
\big]\,\mathrm dr<\infty},
$
where, for $r\ges0$,
$$
\begin{aligned}
&b_0(r)
:=\sup_{(u,v)\in U\times V}|b(r,0,u,v)|,\q
\sigma_0(r)
 :=\sup_{(u,v)\in U\times V}|\sigma(r,0,u,v)|,\\
&{\gamma_0(r)
 :=\sup_{(u,v)\in U\times V}
 |\gamma(r,\cdot,0,u,v)|_{\lambda,2}},\q f_0(r)
 :=\sup_{(u,v)\in U\times V}
|f(r,0,0,0,0,u,v)|.
\end{aligned}
$$

\end{description}

For notational convenience, for $r\ges0$ and
$0\les t<s\les\infty$, set
\begingroup
$$
 \bar q(r)
:=\ell_b(r)+\ell_\sigma(r)^2
+|\ell_\gamma(r,\cdot)|_{\lambda,2}^2,\qq
{} \bar q^0(r)
:= b_0(r)+\sigma_0(r)^2
+\gamma_0(r)^2,
$$
\endgroup
and
\begingroup
$$
\begin{aligned}
\bar\Theta(t,s)
:={}&\Big(\int_t^s\ell_b(r)\,\mathrm dr\Big)^2
+\int_t^s\ell_\sigma(r)^2\,\mathrm dr
+\int_t^s|\ell_\gamma(r,\cdot)|_{\lambda,2}^2\,\mathrm dr,\\
\bar\eta(t,s)
:={}&\Big(\int_t^s b_0(r)\,\mathrm dr\Big)^2
+\int_t^s\sigma_0(r)^2\,\mathrm dr
+\int_t^s\gamma_0(r)^2\,\mathrm dr.
\end{aligned}
$$
\endgroup
By {\bf (G1)} and {\bf (G5)},
${\bar q,\bar q^0},f_0\in L^1(0,\infty)$,
and ${\bar\Theta(t,s)}$ and ${\bar\eta(t,s)}$
are finite for all $0\les t<s\les\infty$.

 \medskip

For $0\les t<T\les\infty$, let $\mathcal U_{t,T}$ and
$\mathcal V_{t,T}$ denote the sets of all $U$-valued and
$V$-valued $\mathbb F$-predictable processes on $[t,T]$,
respectively.
Their elements are called admissible controls for
Player~I and Player~II, respectively.
For $(t,x)\in[0,\infty)\times\mathbb R^n$ and
$(u,v)\in\mathcal U_{t,\infty}\times\mathcal V_{t,\infty}$,
consider the controlled state equation
\begin{equation}\label{eq:game-controlled-state}
\begin{cases}
\begin{aligned}
\mathrm dX_s^{t,x;u,v}
={}&b(s,X_s^{t,x;u,v},u_s,v_s)\,\mathrm ds+\sigma(s,X_s^{t,x;u,v},u_s,v_s)\,\mathrm dB_s\\
&+\int_E\gamma(s,e,X_{s-}^{t,x;u,v},u_s,v_s)
\,\widetilde\mu(\mathrm ds,\mathrm de),
\qquad s\ges t,
\end{aligned}\\[1ex]
X_t^{t,x;u,v}=x,
\end{cases}
\end{equation}
and the following infinite-horizon BSDE with jumps:
\begin{equation}\label{eq:game-controlled-bsde}
\begin{cases}
\ds Y_s^{t,x;u,v}=Y_T^{t,x;u,v}
+\int_s^T f\bigl(r,X_r^{t,x;u,v},Y_r^{t,x;u,v},
Z_r^{t,x;u,v},\Psi_r(K_r^{t,x;u,v}),u_r,v_r\bigr)
\,\mathrm dr\\[1ex]
\displaystyle\qquad
-\int_s^T Z_r^{t,x;u,v}\,\mathrm dB_r
-\int_s^T\int_E K_r^{t,x;u,v}(e)
\,\widetilde\mu(\mathrm dr,\mathrm de),
\quad t\les s\les T<\infty,\\[1ex]
\displaystyle
Y_T^{t,x;u,v}\to\Phi(X_\infty^{t,x;u,v})
\quad\text{in }L^2(\Omega)\text{ and }\mathbb P\text{-a.s.},
\quad\text{as }T\to\infty,
\end{cases}
\end{equation}
with
$ \ds
\Psi_r(k):=\int_E\rho(r,e)k(e)\,\lambda(\mathrm d e),
$ $k\in L_\l^2(E)$.

\ms

Applying Theorem~\ref{thm:infinite-horizon-sde},
Corollary~\ref{cor:short-time-state-increment},
and Lemma~\ref{lem:state-limit},
we obtain the following properties of the
controlled state equation
\eqref{eq:game-controlled-state}.

\begin{lemma}\label{lemma:game-state-estimates}\sl
Assume {\bf(G1)} and {\bf(G5)}.
For every $t\ges0$, $x\in\mathbb R^n$, and
$(u,v)\in\mathcal U_{t,\infty}\times\mathcal V_{t,\infty}$,
equation~\eqref{eq:game-controlled-state} admits a unique
$\mathbb F$-adapted c\`adl\`ag strong solution
$X^{t,x;u,v}\in{{\mathcal S^2}}
(t,\infty;\mathbb R^n)$.
Moreover, there exists
$X_\infty^{t,x;u,v}\in
{L_{\mathcal F_\infty}^2}(\Omega;\mathbb R^n)$
such that
$
X_s^{t,x;u,v}\to X_\infty^{t,x;u,v}
$ in {$L^2(\Omega;\mathbb R^n)$}
 and $\mathbb P\text{-a.s.}
$ as $s\to\infty.
$
Furthermore, there exists a constant $C>0$,
independent of $t,x,u,v$, such that
\begin{equation}\label{eq:game-state-growth}
\begin{aligned}
&\mathbb E_t\Big[
\sup_{s\ges t}|X_s^{t,x;u,v}|^{2}
\Big]
\les C\bigl(|x|^{2}+{\bar\eta(t,\infty)}\bigr),\q  \mathbb P\mbox{-a.s.}\\
 &\dbE_t\Big[
\sup_{s\ges t}|X_s^{t,x;u,v}-x|^2
\Big] \les C\big[
|x|^2\bar\Theta(t,\infty)
+\bar\eta(t,\infty)
\big],
\qquad \dbP\mbox{-a.s.}
\end{aligned}\end{equation}
and, for every $x,x'\in\mathbb R^n$,
\begin{equation}\label{eq:game-state-stability}
\mathbb E_t\Big[
\sup_{s\ges t}
|X_s^{t,x;u,v}-X_s^{t,x';u,v}|^{2}
\Big]
\les C|x-x'|^{2}, \q  \mathbb P\mbox{-a.s.}
\end{equation}
There also exists a constant
$\delta_0\in(0,1]$, depending only on
the Lipschitz profiles in {\bf (G1)},
such that, for every $0<\delta\les\delta_0$,
\begin{equation}\label{eq:game-state-short-time}
\mathbb E_t\Big[
\sup_{s\in[t,t+\delta]}|X_s^{t,x;u,v}-x|^{2}
\Big]
\les C\Big[
|x|^{2}{\bar\Theta(t,t+\delta)}
+{\bar\eta(t,t+\delta)}
\Big],\q  \mathbb P\mbox{-a.s.}
\end{equation}
\end{lemma}

Similarly, applying
Theorem~\ref{thm:existence-infinite-bsde},
Proposition~\ref{prop:bsde-apriori},
Corollary~\ref{prop:bsde-Y-lipschitz},
and Lemma~\ref{lemma:game-state-estimates},
we obtain the following properties of
the controlled BSDE
\eqref{eq:game-controlled-bsde}.

\begin{lemma}\label{lemma:game-bsde-estimates}\sl
Assume {\bf(G1)--(G3)} and  {\bf(G5)}.
For every $t\ges0$, $x\in\mathbb R^n$, and
$(u,v)\in\mathcal U_{t,\infty}\times\mathcal V_{t,\infty}$,
equation~\eqref{eq:game-controlled-bsde} admits a
unique solution
$
 (Y^{t,x;u,v},Z^{t,x;u,v},
 K^{t,x;u,v})  \in
\mathbb B_\lambda^2(t,\infty).
$
Moreover,   there exists a constant $C>0$,
depending only on the structural
and integrability bounds in
{\bf (G1)--(G3)} and {\bf (G5)}, such that
\begingroup
\begin{equation}\label{eq:game-bsde-growth}
\begin{aligned}
\mathbb E_t\Big[
&\sup_{s\ges t}|Y_s^{t,x;u,v}|^2
+\int_t^\infty|Z_r^{t,x;u,v}|^2\,\mathrm d r
+\int_t^\infty
|K_r^{t,x;u,v}(\cdot)|_{\lambda,2}^2
\,\mathrm d r
\Big]
\les C(1+|x|^2),\q  \mathbb P\mbox{-a.s.}
\end{aligned}
\end{equation}
\endgroup
For every $x,x'\in\mathbb R^n$, we also have
\begingroup
\begin{equation}\label{eq:game-bsde-state-stability}
\begin{aligned}
&\mathbb E_t\Big[
\sup_{s\ges t}
|Y_s^{t,x;u,v}-Y_s^{t,x';u,v}|^2+\int_t^\infty
|Z_r^{t,x;u,v}-Z_r^{t,x';u,v}|^2
\,\mathrm d r\\
&\qq+\int_t^\infty
|K_r^{t,x;u,v}(\cdot)
-K_r^{t,x';u,v}(\cdot)|_{\lambda,2}^2
\,\mathrm d r
\Big] \les C|x-x'|^2,\q  \mathbb P\mbox{-a.s.}
\end{aligned}
\end{equation}
\endgroup

\end{lemma}

For $(t,x)\in[0,\infty)\times\mathbb R^n$ and
$(u,v)\in\mathcal U_{t,\infty}\times\mathcal V_{t,\infty}$,
define the recursive cost functional  by
\begin{equation}\label{eq:game-payoff}
J(t,x;u,v):=Y_t^{t,x;u,v}.
\end{equation}
It  can be extended to any random initial state
$\xi\in{L_{\mathcal F_t}^2}(\Omega;\mathbb R^n)$, i.e.,
$ \ds
J(t,\xi;u,v)=Y_t^{t,\xi;u,v},
$ $\mathbb P\text{-a.s.}$ See \cite{BL2008, BuckdahnHuLi2011}.

For any stopping time $S$ taking values in $[t,\infty]$,
set
$
\llbracket t,S\rrbracket
:=\bigl\{(\omega,r)\in\Omega\times[t,\infty):
t\les r\les S(\omega)\bigr\}.
$
Following \cite{BL2008}, we introduce the
following nonanticipative strategies.

\begin{definition}\label{def:game-strategies}\sl
Let $0\les t<T\les\infty$.
A nonanticipative strategy for Player~I on $[t,T]$ is a mapping
$\alpha:\mathcal V_{t,T}\to\mathcal U_{t,T}$ such that,
for every stopping time $S$ with values in $[t,T]$ and every
$v_1,v_2\in\mathcal V_{t,T}$,
$
v_1=v_2
$ $\mathrm d\mathbb P\otimes\mathrm d r$-a.e. on $\llbracket t,S\rrbracket
$
implies
$
\alpha(v_1)=\alpha(v_2)
\quad\mathrm d\mathbb P\otimes\mathrm d r\text{-a.e. on }\llbracket t,S\rrbracket.
$
The set of such strategies is denoted by $\mathcal A_{t,T}$.
Nonanticipative strategies
$\beta:\mathcal U_{t,T}\to\mathcal V_{t,T}$
for Player~II are defined analogously, and their set is denoted
by $\mathcal B_{t,T}$.
\end{definition}

Since Player~I seeks to maximize the payoff
while Player~II seeks to minimize it,
we define the lower and upper value functions
through Elliott--Kalton nonanticipative strategies as follows,
\begin{equation}\label{eq:game-lower-value}
W(t,x)
:=
\essinf_{\beta\in\mathcal B_{t,\infty}}
\esssup_{u\in\mathcal U_{t,\infty}}
J(t,x;u,\beta(u))
\end{equation}
and
\begin{equation}\label{eq:game-upper-value}
U(t,x)
:=
\esssup_{\alpha\in\mathcal A_{t,\infty}}
\essinf_{v\in\mathcal V_{t,\infty}}
J(t,x;\alpha(v),v).
\end{equation}
These essential extrema are a priori
$\mathcal F_t$-measurable random variables.
However, the following proposition shows that they are
in fact deterministic.

\begin{proposition}\label{prop:game-determinacy}\sl
Assume {\bf(G1)--(G3)} and  {\bf(G5)}.
For every $(t,x)\in[0,\infty)\times\mathbb R^n$,
we have
$$
W(t,x)=\mathbb E[W(t,x)],
\qquad
U(t,x)=\mathbb E[U(t,x)],
\quad\mathbb P\text{-a.s.}
$$
\end{proposition}

The conclusion follows by adapting the invariance argument
of {\cite[Lemmas~3.1--3.2]{BuckdahnHuLi2011}}.
Indeed, the measure-preserving transformations in
{\cite[(3.11)--(3.12)]{BuckdahnHuLi2011}}
preserve the Brownian and Poisson increments after time $t$.
In our infinite-horizon setting, the terminal
convergence in $L^2$ and almost surely is also preserved
under these transformations.
Thus, uniqueness of the controlled equations yields
the corresponding transformation identity for the payoff
and the invariance of both value functions.
Their determinacy then follows from
{\cite[Lemma~3.2]{BuckdahnHuLi2011}}.

By Lemma~\ref{lemma:game-bsde-estimates}
and the definitions of the lower and upper value functions,
we obtain the following properties.

\begin{proposition}\label{prop:game-value-regularity}\sl
Assume {\bf(G1)--(G3)} and  {\bf(G5)}.
There exists a constant $C>0$, depending only on  the   structural bounds, such that, for all $t\ges0$ and $x,x'\in\mathbb R^n$,
\begin{equation}\label{eq:game-value-growth}
\begin{aligned}
&|W(t,x)|+|U(t,x)|\les C(1+|x|) ,\\
&|W(t,x)-W(t,x')|+|U(t,x)-U(t,x')|
\les C|x-x'|.
\end{aligned}\end{equation}

\end{proposition}

In what follows, we focus on the lower value function $W$,
as the upper value function $U$ can be treated analogously.

\subsection{Dynamic programming principle}

In this subsection, we establish the dynamic programming
principle (DPP, for short) for the stochastic differential
game. To this end, following the backward stochastic
semigroup approach introduced by Peng~\cite{Peng1997},
we {recall its extension to BSDEs with jumps
from \cite{BuckdahnHuLi2011}.}

\begin{definition}
\label{def:game-backward-semigroup}\sl
Fix $0\les t<T<\infty$, $x\in\dbR^n$, and
$(u,v)\in\mathcal U_{t,T}\times\mathcal V_{t,T}$.
For any $\zeta\in L_{\mathcal F_T}^2(\Omega;\dbR)$,
the backward stochastic semigroup is defined by
\begin{equation}\label{eq:game-backward-semigroup}
G_{s,T}^{t,x;u,v}[\zeta]
:=\widetilde Y_s,
\qquad s\in[t,T],
\end{equation}
where $(\widetilde Y,\widetilde Z,\widetilde K)$
is the solution of the finite-horizon BSDE with jumps
\begin{equation}\label{eq:game-finite-bsde}
\begin{aligned}
\widetilde Y_s
={}&\zeta
+\int_s^T
f\bigl(
r,X_r^{t,x;u,v},
\widetilde Y_r,\widetilde Z_r,
\Psi_r(\widetilde K_r),u_r,v_r
\bigr)\,\mathrm dr\\
&-\int_s^T\widetilde Z_r\,\mathrm dB_r
-\int_s^T\int_E\widetilde K_r(e)\,
\widetilde\mu(\mathrm dr,\mathrm de),
\qquad s\in[t,T],
\end{aligned}
\end{equation}
and $X^{t,x;u,v}$ denotes the solution of
\eqref{eq:game-controlled-state} on $[t,T]$.
\end{definition}

Under {\bf (G1)--(G2)} and {\bf (G5)},
the finite-horizon well-posedness result
for BSDEs with jumps
(see \cite[Chapter~1, Lemma~6 and Theorem~7]{Situ2000})
ensures that \eqref{eq:game-finite-bsde}
admits a unique solution
$
(\widetilde Y,\widetilde Z,\widetilde K)
\in\mathbb B_\lambda^2(t,T).
$
Hence, the backward stochastic semigroup
is well-defined.

By the uniqueness of the  BSDE,
the backward   semigroup satisfies
the flow property
$$
G_{s,T}^{t,x;u,v}[\zeta]
=
G_{s,r}^{t,x;u,v}
\bigl[G_{r,T}^{t,x;u,v}[\zeta]\bigr],
\qquad t\les s<r<T.
$$
In particular, for
$(u,v)\in\mathcal U_{t,\infty}\times\mathcal V_{t,\infty}$
and $t\les s<\infty$,
\begin{equation}\label{eq:game-payoff-flow}
J(t,x;u,v)
 =G_{t,s}^{t,x;u,v}
\bigl[Y_s^{t,x;u,v}\bigr] =G_{t,s}^{t,x;u,v}
\bigl[J(s,X_s^{t,x;u,v};u,v)\bigr].
\end{equation}

\begin{proposition}[Dynamic programming principle]
\label{pro:4.1}\sl
Assume {\bf (G1)--(G5)}.
Then, for every $0\les t<s<\infty$
and $x\in\dbR^n$, the lower value
function satisfies
\begin{equation}\label{eq:game-dpp-lower}
W(t,x)
=
\essinf_{\beta\in\mathcal B_{t,s}}
\esssup_{u\in\mathcal U_{t,s}}
G_{t,s}^{t,x;u,\beta(u)}
\bigl[W(s,X_s^{t,x;u,\beta(u)})\bigr].
\end{equation}

\end{proposition}

The proof follows the standard backward-semigroup
and nonanticipative strategy-pasting argument;
see {\cite[Theorem~3.1 and Appendix~A.2]{BuckdahnHuLi2011}}.
The required comparison, conditional $L^2$ stability,
and payoff flow properties have already been established
under {\bf (G1)--(G5)}, also for random initial states.
Moreover, the control-uniform integrability bounds
in {\bf (G5)} ensure that the required pasting operations
preserve admissibility on the infinite-horizon
continuation interval.
Hence the same argument applies to the present setting,
and we omit the details.

Using the dynamic programming principle and the preceding
spatial Lipschitz estimates, we establish the continuity
of the value functions with respect to time.

\begin{theorem}\label{thm:4.1}\sl
Assume {\bf (G1)--(G5)}.
Then, the lower value function  $W(\cdot,x)$
is continuous on $[0,\infty)$.
Consequently,
$
W \in C([0,\infty)\times\dbR^n).
$
\end{theorem}

 \begin{proof}

Fix $0\les t<T<\infty$ with $\delta:=T-t\les\delta_0$,
where $\delta_0$ is given in
Lemma~\ref{lemma:game-state-estimates},
and fix $x\in\mathbb R^n$.
Throughout, $C>0$ denotes a generic constant independent
of $t,T,x,u,v$.

For any $(u,v)\in\mathcal U_{t,T}\times\mathcal V_{t,T}$,
write $X:=X^{t,x;u,v}$. We have the following decomposition
\begin{equation}\label{eq:game-time-decomposition}
G_{t,T}^{t,x;u,v}[W(T,X_T)]-W(T,x)
=I_\delta^1+I_\delta^2,
\end{equation}
where
$
I_\delta^1
:=G_{t,T}^{t,x;u,v}[W(T,X_T)]
-G_{t,T}^{t,x;u,v}[W(T,x)]
$
and
$
I_\delta^2
:=G_{t,T}^{t,x;u,v}[W(T,x)]-W(T,x).
$

For convenience, set
$ \ds
{d(t,T;x)
:=\big(
|x|^2\bar\Theta(t,T)
+\bar\eta(t,T)
\big)^{{\frac{1}{2}}}}.
$
We first estimate $I_\delta^1$.
Since the two BSDEs have the same generator,
the conditional $L^2$ stability estimate yields
$$
|I_\delta^1|^2
\les C\mathbb E_t\Big[
|W(T,X_T)-W(T,x)|^2
\Big]
\les C\mathbb E_t[|X_T-x|^2]
\les C{d(t,T;x)^2}.
$$
To estimate $I_\delta^2$, let
$(\bar Y,\bar Z,\bar K)$ be the solution defining
$G_{t,T}^{t,x;u,v}[W(T,x)]$.
The constant triple $(W(T,x),0,0)$ solves the BSDE with zero
generator and terminal value $W(T,x)$. Applying It\^o's formula to
$|\bar Y_r-W(T,x)|^2$ and the conditional $L^2$ estimate, we obtain
\begin{equation}\label{eq:game-time-I2-estimate}
\begin{aligned}
&|I_\delta^2|^2
\les
C\mathbb E_t\Big[ \Big(\int_t^T
|f(r,X_r,W(T,x),0,0,u_r,v_r)|
\,\mathrm dr\Big)^2
\Big]\\
&\les C(1+|x|^2)
\Big(\int_t^T(f_0(r)+\ell_x(r)+\ell_y(r))
\,\mathrm dr\Big)^2.
\end{aligned}
\end{equation}

Combining the estimates for $I_\delta^1$ and $I_\delta^2$
in \eqref{eq:game-time-decomposition} yields a deterministic
bound for
$|G_{t,T}^{t,x;u,v}[W(T,X_T)]-W(T,x)|$,
uniformly over $(u,v)$.
Taking the essential supremum and infimum in the DPP
\eqref{eq:game-dpp-lower} therefore gives
\begin{equation}\label{eq:game-value-time-modulus}
\begin{aligned}
|W(t,x)-W(T,x)|
\les C\Big\{
&{d(t,T;x)}
+(1+|x|)\int_t^T
\bigl(f_0(r)+\ell_x(r)+\ell_y(r)\bigr)
\,\mathrm dr
\Big\}.
\end{aligned}
\end{equation}
By {\bf (G1)}, {\bf (G2)}, and {\bf (G5)},
the right-hand side of
\eqref{eq:game-value-time-modulus}
tends to zero as $T-t\to0$,
uniformly in $t,T$ for each fixed $x$.
Thus $W(\cdot,x)$ is continuous on
$[0,\infty)$.
Together with the uniform spatial Lipschitz
continuity in
Proposition~\ref{prop:game-value-regularity},
this yields
$
W\in C([0,\infty)\times\dbR^n).
$

\end{proof}

\begin{proposition}
\label{prop:game-terminal-uniform}\sl
Assume {\bf (G1)--(G3)} and {\bf (G5)}. Then, for every $R>0$,
\begin{equation}\label{eq:game-terminal-uniform}
\lim_{t\to\infty}\sup_{|x|\les R}
\bigl(
|W(t,x)-\Phi(x)|
+|U(t,x)-\Phi(x)|
\bigr)=0.
\end{equation}
\end{proposition}

\begin{proof}
Fix $R>0$. For $|x|\les R$ and
$(u,v)\in\mathcal U_{t,\infty}
\times\mathcal V_{t,\infty}$,
write $X:=X^{t,x;u,v}$ and set
$ d_R(t)
:=\big(R^2\bar\Theta(t,\infty)
+\bar\eta(t,\infty)\big)^\frac 12.
$
Obviously, $d_R(t)\to0$ as $t\to\infty$.
Moreover, by \eqref{eq:game-state-growth}, we have
\begin{equation}\label{eq:terminal-state-uniform}
\mathbb E_t\Big[
\sup_{r\ges t}|X_r-x|^{2}
\Big]
\les C{d_R(t)^2}.
\end{equation}

Next, consider the auxiliary state process
$\bar X_r\equiv x$ and the constant BSDE solution
$
\bar Y_r\equiv\Phi(x),$ $
\bar Z_r=0,$ $
\bar K_r=0,
$
corresponding to the zero generator and
the constant terminal function
$\bar\Phi(\cdot)\equiv\Phi(x)$.
Applying the infinite-horizon BSDE stability
estimate \eqref{eq:stability-infinite-bsde},
we obtain
$$
\begin{aligned}
&|J(t,x;u,v)-\Phi(x)|^{2} \les C\mathbb E_t\Big[
|\Phi(X_\infty)-\Phi(x)|^2
+\Big(\int_t^\infty
|f(r,X_r,\Phi(x),0,0,u_r,v_r)|
\,\mathrm dr\Big)^2
\Big]\\
&\les C\mathbb E_t\Big[
\sup_{r\ges t}|X_r-x|^{2} +
\Big(\int_t^\infty
|f(r,{X_r},\Phi(x),0,0,u_r,v_r)|
\,\mathrm dr\Big)^{2}
\Big]\\
&\les C\mathbb E_t\Big[
\sup_{r\ges t}|X_r-x|^{2}
+\Big(\int_t^\infty
\bigl(f_0(r)+\ell_x(r)|X_r|
+\ell_y(r)|\Phi(x)|\bigr)
\,\mathrm dr\Big)^{2}
\Big]\\
 & \les C_R\Big\{
{d_R(t)^2}+
\Big(\int_t^\infty
\bigl(f_0(r)+\ell_x(r)+\ell_y(r)\bigr)
\,\mathrm dr\Big)^{2}
\Big\}.
\end{aligned}
$$
The right-hand side is deterministic,
independent of $x,u,v$ for $|x|\les R$,
and tends to zero as $t\to\infty$.
Taking the corresponding essential extrema,
we obtain
$$
 \sup_{|x|\les R}
\bigl(
|W(t,x)-\Phi(x)|
+|U(t,x)-\Phi(x)|
\bigr) \les C_R\Big\{
d_R(t)+
\int_t^\infty
\bigl(f_0(r)+\ell_x(r)+\ell_y(r)\bigr)
\,\mathrm dr
\Big\}
\to0
$$
as $t\to\infty$, proving
\eqref{eq:game-terminal-uniform}.
\end{proof}

\subsection{Viscosity characterization of the lower value function}

In this subsection, we show that the lower value function
$W$ is a viscosity solution of the corresponding
non-autonomous HJBI integro-differential equation
with a terminal condition at infinity.
In addition to {\bf (G1)--(G5)}, we impose the following
regularity and jump-size condition.

\begin{description}

\item[(G6)]
The mapping
$t\mapsto\rho(t,\cdot)$
is continuous with values in $L_\lambda^2(E)$.
Moreover, for every $L>0$, there exists
a constant $C_L>0$ such that
$$
|\gamma(t,e,x,u,v)|
\les C_L(1+|x|)(1\wedge|e|),
$$
for all
$(t,x,u,v)\in[0,L]\times\mathbb R^n\times U\times V$
and $\lambda$-a.e.\ $e\in E$.
\end{description}
Note that, under {\bf (G1)} and {\bf (G6)}, the mapping
$
(t,x,u,v)\longmapsto
\gamma(t,\cdot,x,u,v)
 $
is continuous with values in
$L_\lambda^2(E;\mathbb R^n)$.

We first introduce the relevant differential
and nonlocal operators.
For
$(t,x,u,v)\in[0,\infty)\times\mathbb R^n\times U\times V$
and a sufficiently smooth function $\varphi$, define
\begin{equation}\label{eq:operator-C-full}
\begin{aligned}
&
\mathfrak L^{u,v}\varphi(t,x)
:=D\varphi(t,x)\cdot b(t,x,u,v)
+\frac12\operatorname{tr}\big[
\sigma\sigma^\top(t,x,u,v)D^2\varphi(t,x)
\big],\\
&\mathfrak B^{u,v}\varphi(t,x)
:=
\int_E\big[
 \varphi\bigl(t,x+\gamma(t,e,x,u,v)\bigr)
-\varphi(t,x) -D\varphi(t,x)\cdot\gamma(t,e,x,u,v)
\big]\lambda(\mathrm de),\\
&\mathfrak C^{u,v}\varphi(t,x)
:=
\int_E\big[
 \varphi\bigl(t,x+\gamma(t,e,x,u,v)\bigr)
-\varphi(t,x)
\big]\rho(t,e)\lambda(\mathrm de),\\
&H^{u,v}[\varphi](t,x)
:=
\mathfrak L^{u,v}\varphi(t,x)
+\mathfrak B^{u,v}\varphi(t,x) +
f\big(
t,x,\varphi(t,x),
D\varphi(t,x)\cdot\sigma(t,x,u,v),
\mathfrak C^{u,v}\varphi(t,x),
u,v
\big).
\end{aligned}
\end{equation}
Consider the following non-autonomous
HJBI integro-differential equation:
\begin{equation}\label{equ:5.1}
\begin{cases}
-\partial_tW(t,x)-H^-[W](t,x)=0,
& (t,x)\in[0,\infty)\times\mathbb R^n,\\[1ex]
\displaystyle
 \lim_{t\to\infty}W(t,x)=\Phi(x),
& \text{locally uniformly in }x\in\dbR^n.
\end{cases}
\end{equation}
where
the lower Isaacs Hamiltonian is given  by
\begin{equation}\label{eq:Hminus}
H^-[\varphi](t,x)
:=
\sup_{u\in U}\inf_{v\in V}
H^{u,v}[\varphi](t,x).
\end{equation}

Let $C_{l,b}^3([0,\infty)\times\dbR^n)$
denote the class of functions $\varphi$ with
continuous derivatives
$\partial_s\varphi$, $D_x^j\varphi$ $(j=1,2,3)$,
and $D_x\partial_s\varphi$, such that, for every $L>0$,
\begin{equation}\label{eq:Isaacs-test-class}
\sup_{\substack{0\les s\les L\\x\in\dbR^n}}
\bigg[ \frac{|\varphi(s,x)|+|\partial_s\varphi(s,x)|}
{1+|x|^2} +\frac{|D_x\varphi(s,x)|
+|D_x\partial_s\varphi(s,x)|}{1+|x|} +|D_x^2\varphi(s,x)|+|D_x^3\varphi(s,x)|
\bigg]<\infty.
\end{equation}
To formulate the notion of viscosity solution
for \eqref{equ:5.1}, we introduce $\delta$-split
nonlocal operators to account for the jump terms.
For $\delta>0$, set
$
E_\delta:=\{e\in E:|e|<\delta\}
$
and
$
E_\delta^c:=E\setminus E_\delta.
$
For a continuous function $w$ with at most
linear growth and
$\varphi\in C_{l,b}^3([0,\infty)\times\mathbb R^n)$,
define
\begin{equation}\label{eq:Bdelta}
\begin{aligned}
&\mathfrak B_\delta^{u,v}(w,\varphi)(t,x)
:=
\int_{E_\delta}
\big[
\varphi(t,x+\gamma)-\varphi(t,x)
-D\varphi(t,x)\cdot\gamma
\big]\lambda(\mathrm de)\\
&\hskip 3 cm+
\int_{E_\delta^c}
\big[
w(t,x+\gamma)-w(t,x)
-D\varphi(t,x)\cdot\gamma
\big]\lambda(\mathrm de),\\
&\mathfrak C_\delta^{u,v}(w,\varphi)(t,x)
:=
\int_{E_\delta}
\big[
\varphi(t,x+\gamma)-\varphi(t,x)
\big]\rho(t,e)\lambda(\mathrm de) +
\int_{E_\delta^c}
\big[
w(t,x+\gamma)-w(t,x)
\big]\rho(t,e)\lambda(\mathrm de),
\end{aligned}
\end{equation}
where, for brevity,
$
\gamma=\gamma(t,e,x,u,v).
$
Accordingly, the $\delta$-split lower Isaacs
Hamiltonian is defined by
\begin{equation}\label{eq:Hdelta-minus}
\begin{aligned}
&\mathcal H_\delta^-(t,x,y,\varphi;w) :=
\sup_{u\in U}\inf_{v\in V}
\Big\{
\mathfrak L^{u,v}\varphi(t,x)
+\mathfrak B_\delta^{u,v}(w,\varphi)(t,x)\\
&\hskip 3.5cm
+f\bigl(
t,x,y,
D\varphi(t,x)\cdot\sigma(t,x,u,v),
\mathfrak C_\delta^{u,v}(w,\varphi)(t,x),
u,v
\bigr)
\Big\}.
\end{aligned}
\end{equation}

We now introduce the notion of viscosity solutions
for HJBI integro-differential equation \eqref{equ:5.1}.

\begin{definition}\label{def:viscosity-Isaacs-equation}\sl
A continuous function
$w:[0,\infty)\times\dbR^n\to\dbR$
with at most linear growth, uniformly in time,
is called

\noindent
(i)
 a viscosity
subsolution of \eqref{equ:5.1} if
$
\ds
\limsup_{\substack{t\to\infty,y\to x}}w(t,y)
\les\Phi(x),
$
for every $x\in\dbR^n$, and, for every
$\varphi\in C_{l,b}^3([0,\infty)\times\dbR^n)$,
whenever $w-\varphi$ attains a local maximum
at $(t,x)\in[0,\infty)\times\dbR^n$,
one has, for all sufficiently small $\delta>0$,
\begin{equation}\label{eq:viscosity-subsolution}
-\partial_t\varphi(t,x)
-\mathcal H_\delta^-
(t,x,w(t,x),\varphi;w)
\les0.
\end{equation}

\noindent
(ii)  a viscosity
supersolution of \eqref{equ:5.1} if
$
\ds
\liminf_{\substack{t\to\infty,y\to x}}w(t,y)
\ges\Phi(x),
$
for every $x\in\dbR^n$, and, for every
$\varphi\in C_{l,b}^3([0,\infty)\times\dbR^n)$,
whenever $w-\varphi$ attains a local minimum
at $(t,x)\in[0,\infty)\times\dbR^n$,
one has, for all sufficiently small $\delta>0$,
\begin{equation}\label{eq:viscosity-supersolution}
-\partial_t\varphi(t,x)
-\mathcal H_\delta^-
(t,x,w(t,x),\varphi;w)
\ges0.
\end{equation}

\noindent
(iii)   a viscosity
solution of \eqref{equ:5.1} if it is both a
viscosity subsolution and a viscosity supersolution.
\end{definition}

Following \cite[Remark~4.3 and Lemma~4.1]{BuckdahnHuLi2011},
we record two useful equivalences in the viscosity solution
framework, concerning local versus global contact points
and $\delta$-split versus full nonlocal operators.

\begin{remark}\label{rem:viscosity-equivalence}\sl
Assume {\bf (G1)--(G6)}, and let
$w\in C([0,\infty)\times\dbR^n)$ have at most
linear growth, uniformly in time.

{\rm (i)}
In Definition~\ref{def:viscosity-Isaacs-equation},
local maxima (resp., local minima) of $w-\varphi$
may equivalently be replaced by global maxima
(resp., global minima) on
$[0,\infty)\times\dbR^n$.

{\rm (ii)}
For every $\varphi\in C_{l,b}^3([0,\infty)\times\dbR^n)$,
the full nonlocal terms
$\mathfrak B^{u,v}\varphi(t,x)$ and
$\mathfrak C^{u,v}\varphi(t,x)$
are well defined and finite.
Moreover, at global contact points,
the viscosity inequalities for subsolutions
(resp., supersolutions) may equivalently
be formulated by replacing
$$
\mathfrak B_\delta^{u,v}(w,\varphi)(t,x)
\quad\text{with}\quad
\mathfrak B^{u,v}\varphi(t,x),\q
\mathfrak C_\delta^{u,v}(w,\varphi)(t,x)
\quad\text{with}\quad
\mathfrak C^{u,v}\varphi(t,x),
$$
where $\varphi$ touches $w$ from above
(resp., below).

\end{remark}

\begin{proof}

Assertion (i) follows from the localization argument in
\cite[Remark 4.3]{BuckdahnHuLi2011}, which applies under
{\bf (G6)} and the uniform linear growth of $w$.
The argument is local in time and preserves the
test-function class used here.

For (ii), by Taylor's formula, the
Cauchy--Schwarz inequality, and {\bf(G6)},
the full nonlocal terms
$\mathfrak B^{u,v}\varphi(t,x)$ and
$\mathfrak C^{u,v}\varphi(t,x)$
are well defined and finite.

At global contact points, the implication
from the split viscosity inequalities to
the full nonlocal inequalities follows
directly from the contact condition,
$\rho\ges0$, and the monotonicity of $f$
in its jump argument.
The converse implication follows from
the one-sided smooth approximation
argument in \cite[Lemma~4.1]{BuckdahnHuLi2011}.

 Consider the subsolution case and fix a
global maximum point $(t_0,x_0)$ of
$w-\varphi$, normalized so that
$w(t_0,x_0)=\varphi(t_0,x_0)$.
Following the one-sided smooth approximation
in \cite[Proof of Lemma~4.1, properties
(i)--(vi)]{BuckdahnHuLi2011},
we construct
$\varphi_\alpha\in C_{l,b}^3(
[0,\infty)\times\dbR^n)$
such that
$
w\les\varphi_\alpha\les\varphi,
$
with $\varphi_\alpha=\varphi$ in a
neighborhood of $(t_0,x_0)$
and $\varphi_\alpha\to w$ locally uniformly
as $\alpha\downarrow0$.
The construction can be localized in time
and space, preserving the test-function
class used here.
Then, for every $(u,v)\in U\times V$,
we have
$$
\begin{aligned}
\mathfrak B^{u,v}\varphi_\alpha(t_0,x_0)
&\les
\mathfrak B_\delta^{u,v}
(\varphi_\alpha,\varphi)(t_0,x_0),\\
\mathfrak C^{u,v}\varphi_\alpha(t_0,x_0)
&\les
\mathfrak C_\delta^{u,v}
(\varphi_\alpha,\varphi)(t_0,x_0).
\end{aligned}
$$

Applying the full nonlocal subsolution
inequality to $\varphi_\alpha$, and using
$\rho\ges0$ and the monotonicity of $f$
in its jump argument, we obtain
$$
0
\les
\partial_t\varphi_\alpha(t_0,x_0)
+H^-[\varphi_\alpha](t_0,x_0)
\les
\partial_t\varphi(t_0,x_0)
+\mathcal H_\delta^-
(t_0,x_0,w(t_0,x_0),\varphi;\varphi_\alpha).
$$

To pass to the limit as $\alpha\downarrow0$,
by the continuity of $f$ and the compactness
of $U\times V$, it suffices to show that,
for every fixed $\delta>0$,
$$
\begin{aligned}
&\sup_{u,v}
\big|
\mathfrak B_\delta^{u,v}
(\varphi_\alpha,\varphi)(t_0,x_0)
-\mathfrak B_\delta^{u,v}
(w,\varphi)(t_0,x_0)
\big|\longrightarrow0,\\
&\sup_{u,v}
\big|
\mathfrak C_\delta^{u,v}
(\varphi_\alpha,\varphi)(t_0,x_0)
-\mathfrak C_\delta^{u,v}
(w,\varphi)(t_0,x_0)
\big|\longrightarrow0.
\end{aligned}
$$
Indeed, by  {\bf (G6)},
all jump landing points
$
x_0+\gamma(t_0,e,x_0,u,v)
$
lie in a common compact set, uniformly over
$(u,v)\in U\times V$ and $e\in E$.
The locally uniform convergence of
$\varphi_\alpha$ to $w$ thus yields
$\varepsilon_\alpha\downarrow0$ such that
$$
\begin{aligned}
&\sup_{u,v}
\big|
\mathfrak B_\delta^{u,v}
(\varphi_\alpha,\varphi)(t_0,x_0)
-\mathfrak B_\delta^{u,v}
(w,\varphi)(t_0,x_0)
\big|
\les
\varepsilon_\alpha\lambda(E_\delta^c)
\longrightarrow0.
\end{aligned}
$$
Similarly, for the weighted nonlocal operator,
the Cauchy--Schwarz inequality gives
$$
\sup_{u,v}
\big|
\mathfrak C_\delta^{u,v}
(\varphi_\alpha,\varphi)(t_0,x_0)
-\mathfrak C_\delta^{u,v}
(w,\varphi)(t_0,x_0)
\big|
\les
\varepsilon_\alpha
\int_{E_\delta^c}
\rho(t_0,e)\lambda(\mathrm de)
\les
\varepsilon_\alpha
\lambda(E_\delta^c)^{\frac12}
|\rho(t_0,\cdot)|_{\lambda,2}
\longrightarrow0.
$$
Both convergences are uniform over
$(u,v)\in U\times V$.
Consequently, letting $\alpha\downarrow0$
in the preceding viscosity inequality yields
$$
0\les
\partial_t\varphi(t_0,x_0)
+\mathcal H_\delta^-
(t_0,x_0,w(t_0,x_0),\varphi;w),
$$
which is the desired split viscosity inequality.
The supersolution case follows analogously.
\end{proof}

Fix $\varphi\in C_{l,b}^3([0,\infty)\times\dbR^n)$.
 For $a\in\dbR^n$, $y\in\dbR$,
$z\in\dbR^d$, and $k\in L_\lambda^2(E)$,
define
\begin{equation}\label{eq:jump-frozen-F}
\begin{aligned}
&F^\varphi(r,a,y,z,k,u,v)
:= \partial_r\varphi(r,a)
+\mathfrak L^{u,v}\varphi(r,a)
+\mathfrak B^{u,v}\varphi(r,a)\\
&\qq\qq\q+f\bigl(r,a,y+\varphi(r,a),
z+\sigma(r,a,u,v)^\top D\varphi(r,a),\Psi_r(k)+
\mathfrak C^{u,v}\varphi(r,a),u,v\bigr).
\end{aligned}
\end{equation}
We next record some properties of $F^\varphi$
that will be used below.
For each $L>0$, define
$$
a_\varphi(r):=C_{L,\varphi}\Big[
 1+\bar q(r)+\bar q^0(r)+f_0(r)
+\ell_x(r)+\ell_y(r) +\ell_z(r)^2
+\ell_k(r)^2|\rho(r,\cdot)|_{\lambda,2}^2
\Big],\qquad r\in[0,L],
$$
where $C_{L,\varphi}>0$ is a sufficiently
large constant.
By {\bf(G1)}, {\bf(G2)}, and {\bf(G5)},
$a_\varphi\in L^1(0,L)$.

\begin{remark}\label{Re-F}\sl
{\rm (i)}
The mapping
$
(r,a,u,v)\longmapsto
F^\varphi(r,a,0,0,0,u,v)
$
is jointly continuous.
This follows from the joint continuity of
$b,\sigma,f$ in {\bf (G1)--(G2)},
the $L_\lambda^2$-continuity of
$\gamma,\rho$ in {\bf (G6)},
and the regularity of $\varphi$.

{\rm (ii)}
For every $L>0$, the generator $F^\varphi$
satisfies
\begin{equation}\label{eq:jump-F-Lipschitz}
\begin{aligned}
&|F^\varphi(r,a,y,z,k,u,v)
-F^\varphi(r,a',y',z',k',u,v)|\\
&\les
a_\varphi(r)(1+|a|^2+|a'|^2)|a-a'|
+\ell_y(r)|y-y'| +\ell_z(r)|z-z'|
+\ell_k(r)|\rho(r,\cdot)|_{\lambda,2}
|k-k'|_{\lambda,2},
\end{aligned}
\end{equation}
for $r\in[0,L]$, $a,a'\in\dbR^n$,
$y,y'\in\dbR$, $z,z'\in\dbR^d$,
$k,k'\in L_\lambda^2(E)$, and
$(u,v)\in U\times V$.

In particular,
\begin{equation}\label{eq:jump-F-difference}
\begin{aligned}
&|F^\varphi(r,a,y,z,k,u,v)
-F^\varphi(r,a',y,z,k,u,v)| \\
&\les
|F^\varphi(r,a,0,0,0,u,v)
-F^\varphi(r,a',0,0,0,u,v)| +2\ell_y(r)|y|
+2\ell_z(r)|z|
+2\ell_k(r)|\rho(r,\cdot)|_{\lambda,2}
|k|_{\lambda,2}.
\end{aligned}
\end{equation}
%

{\rm (iii)}
For every $L>0$, we have
\begin{equation}\label{eq:freeze-g-growth}
|F^\varphi(r,a,0,0,0,u,v)|
\les a_\varphi(r)(1+|a|^2),
\end{equation}
for all $r\in[0,L]$, $a\in\dbR^n$,
and $(u,v)\in U\times V$.

Indeed, the growth conditions in
{\bf (G1)} and {\bf (G5)},
together with Taylor's formula and
the Cauchy--Schwarz inequality, yield
$$
\begin{aligned}
&|\mathfrak L^{u,v}\varphi(r,a)|
+|\mathfrak B^{u,v}\varphi(r,a)|
\les C_{L,\varphi}
(\bar q(r)+\bar q^0(r))(1+|a|^2),\\
&|\mathfrak C^{u,v}\varphi(r,a)|
\les C_{L,\varphi}(1+|a|)
|\rho(r,\cdot)|_{\lambda,2}
|\gamma(r,\cdot,a,u,v)|_{\lambda,2}.
\end{aligned}
$$
Combining these estimates with
the weighted derivative bounds in
\eqref{eq:Isaacs-test-class},
the Lipschitz property of $f$,
and Young's inequality proves
\eqref{eq:freeze-g-growth}.
\end{remark}

Fix $(t,x)\in[0,\infty)\times\dbR^n$ and $0<h\les1$.
For $(u,v)\in\mathcal U_{t,t+h}\times\mathcal V_{t,t+h}$,
write $X:=X^{t,x;u,v}$.
Consider the following BSDEs with zero terminal values
on $[t,t+h]$:
\begin{equation}\label{eq:jump-Y1}
\begin{aligned}
Y_s^{1;u,v}
={}&\int_s^{t+h}
F^\varphi(r,X_r,Y_r^{1;u,v},Z_r^{1;u,v},
K_r^{1;u,v},u_r,v_r)\,\mathrm dr\\
&-\int_s^{t+h}Z_r^{1;u,v}\,\mathrm dB_r
-\int_s^{t+h}\int_E K_r^{1;u,v}(e)
\,\widetilde\mu(\mathrm dr,\mathrm de),
\end{aligned}
\end{equation}
and
\begin{equation}\label{eq:jump-Y2}
\begin{aligned}
Y_s^{2;u,v}
={}&\int_s^{t+h}
F^\varphi(r,x,Y_r^{2;u,v},Z_r^{2;u,v},
K_r^{2;u,v},u_r,v_r)\,\mathrm dr\\
&-\int_s^{t+h}Z_r^{2;u,v}\,\mathrm dB_r
-\int_s^{t+h}\int_E K_r^{2;u,v}(e)
\,\widetilde\mu(\mathrm dr,\mathrm de).
\end{aligned}
\end{equation}
Both BSDEs are well posed under {\bf (G1)--(G6)}.
By It\^o's formula and BSDE uniqueness, we have
\begin{equation}\label{eq:jump-semigroup-identity}
Y_s^{1;u,v}
=G_{s,t+h}^{t,x;u,v}
\bigl[\varphi(t+h,X_{t+h})\bigr]
-\varphi(s,X_s),\qquad s\in[t,t+h].
\end{equation}

The following lemma provides a uniform short-time
estimate for the difference between the solutions
of \eqref{eq:jump-Y1} and \eqref{eq:jump-Y2}.

\begin{lemma}
\label{lem:jump-short-freeze}\sl
Assume {\bf (G1)--(G6)}. For every $M,R>0$,
there exists a deterministic function
$\varepsilon_{M,R,\varphi}:(0,1]\to[0,\infty)$ with
$\varepsilon_{M,R,\varphi}(h)\to0$ as $h\downarrow0$
such that, for $t\in[0,M]$, $|x|\les R$, and all
sufficiently small $h>0$,
\begin{equation}\label{eq:jump-short-freeze}
|Y_t^{1;u,v}-Y_t^{2;u,v}|
\les h\varepsilon_{M,R,\varphi}(h),
\qquad\mathbb P\text{-a.s.},
\end{equation}
uniformly over all admissible control pairs $(u,v)$.
\end{lemma}

The proof involves several technical estimates.
To keep the main argument focused,
we defer it to Appendix~\ref{app:proof-jump-short-freeze}.

\begin{lemma}
\label{lem:jump-ode}\sl
Assume {\bf (G1)--(G6)}.
Let $y^h$ be the unique deterministic
solution of
\begin{equation}\label{eq:jump-frozen-ode}
-\frac{\mathrm d}{\mathrm dr}y_r^h
=F_0^\varphi(r,x,y_r^h),
\qquad y_{t+h}^h=0.
\end{equation}
where $\ds
F_0^\varphi(r,x,y)
:=\sup_{u\in U}\inf_{v\in V}
F^\varphi(r,x,y,0,0,u,v).
$
Then
\begin{equation}\label{eq:jump-frozen-value}
\essinf_{\beta\in\mathcal B_{t,t+h}}
\esssup_{u\in\mathcal U_{t,t+h}}
Y_t^{2;u,\beta(u)}
=y_t^h,\qquad\mathbb P\text{-a.s.}
\end{equation}

\end{lemma}

\begin{proof}
By the properties of $F^\varphi$
in  Remark \ref{Re-F}, together with
{\bf(G2)} and the compactness of $U\times V$,
$F_0^\varphi(r,x,y)$ is jointly continuous
in $(r,y)$, and
\bel{F0}
\begin{aligned}
&|F_0^\varphi(r,x,y)
-F_0^\varphi(r,x,y')|
\les \ell_y(r)|y-y'|,\\
&|F_0^\varphi(r,x,0)|
\les a_\varphi(r)(1+|x|^2).
\end{aligned}
\ee
Since $\ell_y,a_\varphi\in L^1(0,t+h)$,   \eqref{eq:jump-frozen-ode} admits
a unique deterministic solution
$y^h\in C^1([t,t+h])$.

By the continuity of $F^\varphi$ and the
compactness of $U$ and $V$, there exist
measurable selectors
$
u^*:[t,t+h]\to U
$
and
$
v^*:[t,t+h]\times U\to V
$
such that
$$
\begin{aligned}
&
\inf_{v\in V}F^\varphi(r,x,y_r^h,0,0,u^*(r),v)
=
\sup_{u\in U}\inf_{v\in V}F^\varphi(r,x,y_r^h,0,0,u,v)
=F_0^\varphi(r,x,y_r^h),\\
& F^\varphi(r,x,y_r^h,0,0,u,v^*(r,u))
=
\inf_{v\in V}F^\varphi(r,x,y_r^h,0,0,u,v).
\end{aligned}$$

For any $\beta\in\mathcal B_{t,t+h}$,
$$
F^\varphi(r,x,y_r^h,0,0,
u^*(r),\beta(u^*)_r)
\ges F_0^\varphi(r,x,y_r^h).
$$
Since $(y^h,0,0)$ solves the BSDE
with generator $F_0^\varphi$,
the comparison theorem yields
$
Y_t^{2;u^*,\beta(u^*)}\ges y_t^h.
$
Consequently,
$$
\essinf_{\beta\in\mathcal B_{t,t+h}}
\esssup_{u\in\mathcal U_{t,t+h}}
Y_t^{2;u,\beta(u)}
\ges y_t^h.
$$

For the reverse inequality, define
$
\beta^*(u)_r:=v^*(r,u_r).
$
The measurability of $v^*$ ensures that
$\beta^*(u)$ is predictable.
Moreover, since $\beta^*(u)_r$ depends
only on the current value $u_r$,
$\beta^*\in\mathcal B_{t,t+h}$.
For every $u\in\mathcal U_{t,t+h}$,
we have
$$
 F^\varphi(r,x,y_r^h,0,0,
u_r,\beta^*(u)_r) =
\inf_{v\in V}
F^\varphi(r,x,y_r^h,0,0,u_r,v)
\les F_0^\varphi(r,x,y_r^h).
$$
Applying the comparison theorem again gives
$
Y_t^{2;u,\beta^*(u)}\les y_t^h.
$
Hence,
$$
\essinf_{\beta\in\mathcal B_{t,t+h}}
\esssup_{u\in\mathcal U_{t,t+h}}
Y_t^{2;u,\beta(u)}
\les
\esssup_{u\in\mathcal U_{t,t+h}}
Y_t^{2;u,\beta^*(u)}
\les y_t^h.
$$
Combining the two inequalities proves
\eqref{eq:jump-frozen-value}.
\end{proof}

\begin{theorem}
\label{thm:jump-viscosity-existence}\sl
Assume {\bf (G1)--(G6)}.
Then the lower value function $W$
is a viscosity solution of the
HJBI integro-differential equation
\eqref{equ:5.1}, in the sense of
Definition~\ref{def:viscosity-Isaacs-equation}.
\end{theorem}

\begin{proof}
The continuity and linear growth of $W$ have already been
established. {Its terminal condition in
Definition~\ref{def:viscosity-Isaacs-equation} follows directly
from Proposition~\ref{prop:game-terminal-uniform} and the
continuity of $\Phi$.} By Remark~\ref{rem:viscosity-equivalence},
it suffices to verify the full-test viscosity
inequalities at global contact points.

Let $\varphi\in C_{l,b}^3([0,\infty)\times\dbR^n)$
and suppose that $W-\varphi$ attains a global
maximum zero at $(t,x)$. Then $W\les\varphi$
and $W(t,x)=\varphi(t,x)$.
The DPP, the   comparison theorem, and
\eqref{eq:jump-semigroup-identity} yield
$$
0
 \les\essinf_{\beta\in\mathcal B_{t,t+h}}
\esssup_{u\in\mathcal U_{t,t+h}}
\Big\{
G_{t,t+h}^{t,x;u,\beta(u)}
[\varphi(t+h,X_{t+h}^{t,x;u,\beta(u)})]
-\varphi(t,x)\Big\} =\essinf_{\beta\in\mathcal B_{t,t+h}}
\esssup_{u\in\mathcal U_{t,t+h}}
Y_t^{1;u,\beta(u)}.
$$
By Lemma~\ref{lem:jump-short-freeze}
and \eqref{eq:jump-frozen-value},
\bel{yth+he}
0\les y_t^h+h\varepsilon_{M,R,\varphi}(h),
\ee
where $M>t$, $R>|x|$, and
$\varepsilon_{M,R,\varphi}(h)\to0$ as $h\downarrow0$.

We next evaluate the short-time limit of $y_t^h$.
By \eqref{F0}, the ODE
\eqref{eq:jump-frozen-ode} and Gronwall's
inequality give
$$
\sup_{r\in[t,t+h]}|y_r^h|
\les C_{M,R,\varphi}h.
$$
Consequently,
$$
\begin{aligned}
&\left|
\frac{y_t^h}{h}
-F_0^\varphi(t,x,0)
\right| \les
\frac1h\int_t^{t+h}
\ell_y(r)|y_r^h|\,\mathrm dr +
\frac1h\int_t^{t+h}
|F_0^\varphi(r,x,0)
-F_0^\varphi(t,x,0)|\,\mathrm dr\\
& \les
C_{M,R,\varphi}
\int_t^{t+h}\ell_y(r)\,\mathrm dr +
\sup_{r\in[t,t+h]}
|F_0^\varphi(r,x,0)
-F_0^\varphi(t,x,0)|
\longrightarrow0,\q h\to 0,
\end{aligned}
$$
where the convergence follows from
$\ell_y\in L^1(0,\infty)$ and the continuity
of $F_0^\varphi(\cdot,x,0)$.
Hence,
$$
\lim_{h\downarrow0}\frac{y_t^h}{h}
=F_0^\varphi(t,x,0)
=\partial_t\varphi(t,x)+H^-[\varphi](t,x).
$$
Letting $h\downarrow0$ and using \eqref{yth+he}, we obtain
$
\partial_t\varphi(t,x)+H^-[\varphi](t,x)\ges0.
 $
Thus the full-test subsolution inequality holds.

If $W-\varphi$ attains a global minimum zero
at $(t,x)$, the DPP and BSDE comparison
inequalities are reversed. The same argument gives
$$
0\ges y_t^h-h\varepsilon_{M,R,\varphi}(h).
$$
Dividing by $h$ and using the same short-time
limit yields
$$
\partial_t\varphi(t,x)+H^-[\varphi](t,x)\les0.
$$
Therefore the full-test supersolution inequality
also holds.
Hence $W$ is a viscosity solution of
\eqref{equ:5.1}.
\end{proof}

\subsection{Comparison principle and uniqueness}
\label{subsec:jump-viscosity-comparison}

We next establish a comparison principle for continuous
viscosity sub- and supersolutions of \eqref{equ:5.1}
with uniform linear growth, provided that at least
one of them is uniformly Lipschitz continuous in
the spatial variable. Since the value function $W$
has this property, the comparison result will yield
uniqueness. The proof relies on the following
additional Lyapunov-type condition.

For a smooth function $\phi$ for which the nonlocal
terms are well defined, set
$$
\mathcal M[\phi](t,x)
:=\sup_{(u,v)\in U\times V}\Big\{
\mathfrak L^{u,v}\phi(t,x)
+\mathfrak B^{u,v}\phi(t,x)
+\ell_z(t)|\sigma(t,x,u,v)^\top D\phi(t,x)|
+\ell_k(t)(\mathfrak C^{u,v}\phi(t,x))^+
\Big\}.
$$

\begin{description}
\item[(G7)]
Let $\mathfrak p(x):=1+|x|^2$.
There exists a nonnegative function
$q\in C([0,\infty))\cap L^1(0,\infty)$ such that
$$
\ell_y(t)\les q(t),
\qquad
\mathcal M[\mathfrak p](t,x)
\les q(t)\mathfrak p(x),
\qquad
(t,x)\in[0,\infty)\times\dbR^n.
$$
\end{description}
Condition {\bf(G7)} does not require
dissipativity of the state dynamics.
Moreover, the spatial Lipschitz profiles
of the coefficients need not be locally bounded.
However, $\ell_y\les q$ and the continuity of $q$
imply that $\ell_y$ is locally bounded.
The following remark provides a sufficient
coefficient condition for {\bf(G7)}.

\begin{remark}\sl
Using the notation in {\bf(G5)}, set
$$
Q_\rho(t):=
\sup_{(u,v)\in U\times V}
\int_E\rho(t,e)|\gamma(t,e,0,u,v)|^2
\lambda(\mathrm de)
+\int_E\rho(t,e)\ell_\gamma(t,e)^2
\lambda(\mathrm de).
$$

{\rm (i)}
For a sufficiently large constant $C>0$, define
$$
\begin{aligned}
\mathcal R(t):=C\Big[
&\ell_y(t)+b_0(t)+\ell_b(t)
+\sigma_0(t)^2+\ell_\sigma(t)^2
+\gamma_0(t)^2
+|\ell_\gamma(t,\cdot)|_{\lambda,2}^2\\
&+\ell_z(t)[\sigma_0(t)+\ell_\sigma(t)]
+\ell_k(t)\Big\{
|\rho(t,\cdot)|_{\lambda,2}
[\gamma_0(t)+|\ell_\gamma(t,\cdot)|_{\lambda,2}]
+Q_\rho(t)\Big\}\Big].
\end{aligned}
$$
Then {\bf(G7)} holds if $\mathcal R$ admits
a nonnegative continuous integrable majorant.
Indeed,
$$
 \mathfrak L^{u,v}\mathfrak p
+\mathfrak B^{u,v}\mathfrak p
 =2x\cdot b+|\sigma|^2
+|\gamma|_{\lambda,2}^2,\qq \mathfrak C^{u,v}\mathfrak p
 =\int_E\rho(t,e)
(2x\cdot\gamma+|\gamma|^2)
\lambda(\mathrm de).
$$
Under {\bf(G6)},
$\mathfrak C^{u,v}\mathfrak p(t,x)$
is well defined and finite.
 Moreover, by {\bf(G1)}, {\bf(G2)}, and
{\bf(G5)}, together with $\rho\ges0$
and the Cauchy--Schwarz inequality,
$$
|\mathfrak C^{u,v}\mathfrak p(t,x)|
\les \Big\{
 |\rho(t,\cdot)|_{\lambda,2}
[\gamma_0(t)+|\ell_\gamma(t,\cdot)|_{\lambda,2}] +Q_\rho(t)
\Big\}\mathfrak p(x).
$$
Combining this with the growth estimates
for $b$, $\sigma$, and $\gamma$ yields
$$
\ell_y(t)\les\mathcal R(t),
\qquad
\mathcal M[\mathfrak p](t,x)
\les\mathcal R(t)\mathfrak p(x).
$$

{\rm (ii)}
Under {\bf(G1)--(G5)}, all terms in
$\mathcal R$ are integrable except possibly
$\ell_kQ_\rho$.
Thus $\ell_kQ_\rho\in L^1(0,\infty)$
is the only additional integrability requirement
in this sufficient criterion.
The existence of a nonnegative continuous
integrable majorant of $\mathcal R$ must
also be required separately.
In particular, {\bf(G7)} does not follow
automatically from {\bf(G1)--(G5)}.
\end{remark}

To establish the comparison principle,
we introduce a time transformation to obtain
a strictly decreasing generator and a
Lyapunov barrier to control the spatial
penalization. Let $q$ be the function
specified in {\bf(G7)}. Define
$$
\Lambda(t):=\int_0^t[q(r)+e^{-r}]\,\mathrm dr,
\qquad
\Gamma(t,x):=\exp\Big(
\int_0^t[e^{-r}-q(r)]\,\mathrm dr
\Big)\mathfrak p(x).
$$
The following lemma establishes the
properties needed for the comparison argument.

\begin{lemma}
\label{lem:jump-integrable-barrier}\sl
Assume {\bf(G7)}.
Then $\Lambda\in C^1([0,\infty))$ is bounded,
$\Gamma\in C^{1,\infty}([0,\infty)\times\dbR^n)$,
and there exist constants $0<c<C$ such that
$
c\mathfrak p(x)\les\Gamma(t,x)
\les C\mathfrak p(x),
$ $ (t,x)\in[0,\infty)\times\dbR^n.
$
Moreover,
\begin{equation}\label{eq:jump-integrable-barrier}
-\partial_t\Gamma(t,x)
+e^{-t}\Gamma(t,x)
-\mathcal M[\Gamma](t,x)\ges0.
\end{equation}
\end{lemma}

\begin{proof}
Set
$
Q:=\int_0^\infty q(r)\,\mathrm dr<\infty
$
and
$
a(t):=\exp\big(
\int_0^t[e^{-r}-q(r)]\,\mathrm dr
\big).
$
Since $q$ is continuous and integrable,
$\Lambda\in C^1([0,\infty))$ is bounded, and
$\Gamma\in C^{1,\infty}([0,\infty)\times\dbR^n)$.
Moreover,
$$
\int_0^t[e^{-r}-q(r)]\,\mathrm dr
=1-e^{-t}-\int_0^t q(r)\,\mathrm dr
\in[-Q,1].
$$
Hence,
$
c\mathfrak p(x)\les\Gamma(t,x)\les C\mathfrak p(x),
$
where $c=e^{-Q}$ and $C=e$.

It remains to verify \eqref{eq:jump-integrable-barrier}.
Since $\Gamma(t,x)=a(t)\mathfrak p(x)$ with $a(t)>0$,
the positive homogeneity of $\mathcal M$
and {\bf (G7)} give
$$
\mathcal M[\Gamma](t,x)
=a(t)\mathcal M[\mathfrak p](t,x)
\les q(t)\Gamma(t,x).
$$
Together with ${\partial_t\Gamma=(e^{-t}-q(t))\Gamma}$,
this yields
$$
 -\partial_t\Gamma(t,x)
+e^{-t}\Gamma(t,x)
-\mathcal M[\Gamma](t,x) =q(t)\Gamma(t,x)-\mathcal M[\Gamma](t,x)
\ges0.
$$
\end{proof}

Fix $T,R,K>0$. For $0<\varepsilon\les1$, set
$
\mathcal D_\varepsilon
:=\big\{(a,a',u,v):\,
 |a|,|a'|\les R,\
|a-a'|\les K\varepsilon,\
 (u,v)\in U\times V\big\}.
 $
For $(a,a',u,v)\in\mathcal D_\varepsilon$, write
$
\Delta\gamma_t(e)
:=\gamma(t,e,a,u,v)-\gamma(t,e,a',u,v).
$
Define the following nonnegative functions on $[0,T]$:
$$
\begin{aligned}
&B_\varepsilon(t)
 :=\sup_{\mathcal D_\varepsilon}
|b(t,a,u,v)-b(t,a',u,v)|,\\
&S_\varepsilon(t)
 :=\sup_{\mathcal D_\varepsilon}
|\sigma(t,a,u,v)-\sigma(t,a',u,v)|,\\
&D_\varepsilon(t)
 :=\frac1\varepsilon\sup_{\mathcal D_\varepsilon}
\Big[
|\sigma(t,a,u,v)-\sigma(t,a',u,v)|^2
+|\Delta\gamma_t|_{\lambda,2}^2
\Big],\\
&K_\varepsilon(t)
:=K\sup_{\mathcal D_\varepsilon}
\Big|\int_E\rho(t,e)\Delta\gamma_t(e)
\lambda(\mathrm de)\Big| +\frac1{2\varepsilon}
\sup_{\mathcal D_\varepsilon}
\int_E\rho(t,e)|\Delta\gamma_t(e)|^2
\lambda(\mathrm de) .
\end{aligned}
$$
We also define the transformed generator
$$
\widetilde f(t,a,y,z,k,u,v)
:=  e^{\Lambda(t)}
f\big(t,a,e^{-\Lambda(t)}y,
e^{-\Lambda(t)}z,
e^{-\Lambda(t)}k,u,v\big) -[q(t)+e^{-t}]y.
$$

The following lemma provides a time-integrable
error envelope for the spatial discrepancies
of the non-autonomous coefficients arising
in the doubling-of-variables argument,
without requiring these errors to vanish
at the maximizing time.

\begin{lemma}
\label{lem:jump-finite-scale-envelope}\sl
Assume {\bf(G1)--(G7)}.
The functions
$B_\varepsilon,S_\varepsilon,D_\varepsilon,K_\varepsilon$
are continuous on $[0,T]$ and satisfy
\begin{equation}\label{eq:jump-finite-scale-vanishing}
 \sup_{t\in[0,T]}
[B_\varepsilon(t)+S_\varepsilon(t)]
\longrightarrow0,\q\int_0^T
[D_\varepsilon(t)+\ell_k(t)K_\varepsilon(t)]
\,\mathrm dt
\longrightarrow0,
\qquad \varepsilon\downarrow0.
\end{equation}
Moreover, for any $B>0$ and $C\ges0$,
there exists a nonnegative function
$e_\varepsilon\in C([0,T])$ such that
\begin{equation}\label{eq:jump-time-envelope}
\int_0^T e_\varepsilon(t)\,\mathrm dt
\longrightarrow0
\qquad(\varepsilon\downarrow0),
\end{equation}
and
$$
 |\widetilde f(t,a,y,z+\zeta,k+\kappa,u,v)
-\widetilde f(t,a',y,z,k,u,v)|
+K B_\varepsilon(t)+C D_\varepsilon(t)
\les e_\varepsilon(t)
$$
for all $t\in[0,T]$,
$(a,a',u,v)\in\mathcal D_\varepsilon$,
$|y|+|z|+|k|\les B$,
$|\zeta|\les K S_\varepsilon(t)$,
and $|\kappa|\les K_\varepsilon(t)$.
The function $e_\varepsilon$ may depend on $B$ and $C$.
\end{lemma}

\begin{proof}

By the joint continuity of $b,\sigma$ and the
$L_\lambda^2$-continuity of $\gamma$,
the functions defining
$B_\varepsilon,S_\varepsilon,D_\varepsilon$
are jointly continuous in $(t,a,a',u,v)$.
Since $\mathcal D_\varepsilon$ is compact,
these three functions are continuous on $[0,T]$.

For $K_\varepsilon$, take
$(t,a,a',u,v)$ and
$(s,\bar a,\bar a',\bar u,\bar v)$
in $[0,T]\times\cD_\e$, and set
$$
\Delta\gamma_t(e)
 :=\gamma(t,e,a,u,v)-\gamma(t,e,a',u,v),\q
\Delta\gamma_s'(e)
 :=\gamma(s,e,\bar a,\bar u,\bar v)
-\gamma(s,e,\bar a',\bar u,\bar v).
$$
By {\bf(G6)}, $\Delta\gamma_t$ is jointly
continuous in $(t,a,a',u,v)$ as an
$L_\lambda^2$-valued function, and
$t\mapsto\rho(t,\cdot)$ is continuous
in $L_\lambda^2$.
Moreover, {\bf(G6)} yields
$$
|\Delta\gamma_t(e)|+|\Delta\gamma_s'(e)|
\les M_{T,R},\qquad \lambda\text{-a.e. }e\in E,
$$
and hence
$$
\big\||\Delta\gamma_t|^2
-|\Delta\gamma_s'|^2\big\|_{\lambda,2}
\les M_{T,R}
|\Delta\gamma_t-\Delta\gamma_s'|_{\lambda,2}.
$$
Consequently, by Cauchy--Schwarz, both
$$
\int_E\rho(t,e)\Delta\gamma_t(e)\lambda(\mathrm de),
\qquad
\int_E\rho(t,e)|\Delta\gamma_t(e)|^2
\lambda(\mathrm de)
$$
are jointly continuous in $(t,a,a',u,v)$.
Taking suprema over the fixed compact set
$\cD_\e$ gives $K_\varepsilon\in C([0,T])$.

We next establish the limits in
\eqref{eq:jump-finite-scale-vanishing}.
By the joint continuity of $b,\sigma$
and the compactness of $[0,T]\times U\times V
\times\overline B_R$, we have
$$
\sup_{t\in[0,T]}
[B_\varepsilon(t)+S_\varepsilon(t)]
\longrightarrow0.
$$
Moreover, {\bf(G1)} and
$|a-a'|\les K\varepsilon$ yield
$
D_\varepsilon(t)
\les K^2\varepsilon
\big[
\ell_\sigma(t)^2
+|\ell_\gamma(t,\cdot)|_{\lambda,2}^2
\big],
$
and hence
$\ds
\int_0^T D_\varepsilon(t)\,\mathrm dt
\longrightarrow0.
$

It remains to consider $K_\varepsilon$.
By {\bf(G1)} and {\bf(G6)},
$
|\Delta\gamma_t|_{\lambda,2}
 \les K\varepsilon
|\ell_\gamma(t,\cdot)|_{\lambda,2},$ and
$|\Delta\gamma_t|_{\lambda,4}^2
 \les M_{T,R}
|\Delta\gamma_t|_{\lambda,2}.
 $
Consequently, by Cauchy--Schwarz,
$$
K_\varepsilon(t)
\les C_{T,R,K}
|\rho(t,\cdot)|_{\lambda,2}
|\ell_\gamma(t,\cdot)|_{\lambda,2}.
$$
To prove pointwise convergence, fix $H>0$,
 $E=\{\rho(t,\cdot)\les H\} \cup\{\rho(t,\cdot)>H\}$.
By the preceding estimates,
$$
\begin{aligned}
K_\varepsilon(t)
\les{}&
K^2\varepsilon
|\rho(t,\cdot)|_{\lambda,2}
|\ell_\gamma(t,\cdot)|_{\lambda,2} +\frac{K^2H\varepsilon}{2}
|\ell_\gamma(t,\cdot)|_{\lambda,2}^2 +C_{T,R,K}
|\ell_\gamma(t,\cdot)|_{\lambda,2}
\(
\int_{\{\rho(t,e)>H\}}
\rho(t,e)^2\lambda(\mathrm de)
\)^{{\frac{1}{2}}}.
\end{aligned}
$$
For almost every $t\in[0,T]$,
letting first $\varepsilon\downarrow0$
and then $H\to\infty$, we obtain
$
K_\varepsilon(t)\longrightarrow0.
$
Finally, noting    $\ds
\int_0^T
\ell_k(t)|\rho(t,\cdot)|_{\lambda,2}
|\ell_\gamma(t,\cdot)|_{\lambda,2}
\,\mathrm dt<\infty,
$
the dominated convergence theorem therefore gives
$$
\int_0^T\ell_k(t)K_\varepsilon(t)
\,\mathrm dt\longrightarrow0.
$$
This proves
\eqref{eq:jump-finite-scale-vanishing}.

\medskip

Fix $B>0$ and $C\ges0$.
Define
$$
\begin{aligned}
\omega_\varepsilon(t):=
\sup\Big\{&
\big|
\widetilde f(t,a,y,
z+K S_\varepsilon(t)\xi,
k+K_\varepsilon(t)\chi,u,v) -\widetilde f(t,a',y,z,k,u,v)
\big|:
(a,a',u,v)\in\mathcal D_\varepsilon,\\
&\qq|y|+|z|+|k|\les B,\
|\xi|,|\chi|\les1
\Big\}.
\end{aligned}
$$
{By joint continuity of the functions involved,}
the supremum over the fixed compact parameter
set defines a nonnegative continuous function
$\omega_\varepsilon$ on $[0,T]$.

By {\bf(G2)} and the definition of
$\widetilde f$,
$$
\omega_\varepsilon(t)
\les
K\varepsilon e^{\Lambda(t)}\ell_x(t)
+K\ell_z(t)S_\varepsilon(t) +\ell_k(t)K_\varepsilon(t).
$$
Since $\Lambda$ is bounded and
$S_\varepsilon(t)\les K\varepsilon\ell_\sigma(t)$,
we obtain
$$
\int_0^T\omega_\varepsilon(t)\,\mathrm dt
\les
C_TK\varepsilon\int_0^T\ell_x(t)\,\mathrm dt +K^2\varepsilon
\int_0^T\ell_z(t)\ell_\sigma(t)\,\mathrm dt +\int_0^T\ell_k(t)K_\varepsilon(t)\,\mathrm dt
\longrightarrow0.
$$
Finally, set
$
e_\varepsilon(t):=
\omega_\varepsilon(t)
+K B_\varepsilon(t)+C D_\varepsilon(t).
$
Then $e_\varepsilon\in C([0,T])$ is nonnegative,
and the required pointwise estimate follows
directly from its definition.
Moreover, by the preceding convergence results,
$
\int_0^T e_\varepsilon(t)\,\mathrm dt
\longrightarrow0,
$
which proves \eqref{eq:jump-time-envelope}
and completes the proof.

\end{proof}

\br\label{Re-upp-loe-Con}\sl
Since $\Phi$ is continuous, the terminal
conditions in
Definition~\ref{def:viscosity-Isaacs-equation},
namely,
$$
\limsup_{\substack{t\to\infty\\y\to x}}w(t,y)
\les\Phi(x),
\qquad
\liminf_{\substack{t\to\infty\\y\to x}}w(t,y)
\ges\Phi(x),
$$
for subsolutions and supersolutions,
respectively, are equivalent to the
following one-sided locally uniform limits
on every compact $K\subset\dbR^n$:
$$
\lim_{t\to\infty}
\sup_{x\in K}(w(t,x)-\Phi(x))^+=0,
\qquad
\lim_{t\to\infty}
\sup_{x\in K}(\Phi(x)-w(t,x))^+=0.
$$
For a continuous viscosity solution satisfying both terminal
conditions, these reduce to
$$
\lim_{t\to\infty}\sup_{x\in K}|w(t,x)-\Phi(x)|=0
\qquad\text{for every compact }K\subset\dbR^n.
$$

\er

\begin{theorem}
\label{thm:jump-one-sided-comparison}\sl
Assume {\bf(G1)--(G7)}.
Let $w_1$ and $w_2$ be, respectively,
a continuous viscosity subsolution and
a continuous viscosity supersolution of
\eqref{equ:5.1}, in the sense of
Definition~\ref{def:viscosity-Isaacs-equation},
both with uniform linear growth.
If at least one of $w_1,w_2$ is uniformly
Lipschitz continuous in $x\in\dbR^n$,
then
$$
w_1(t,x)\les w_2(t,x),
\qquad (t,x)\in[0,\infty)\times\dbR^n.
$$
\end{theorem}

\begin{proof}
We first assume that $w_2$ is uniformly
Lipschitz continuous in $x\in\dbR^n$.
Suppose, by contradiction, that
$w_1(t_*,x_*)>w_2(t_*,x_*)$
 for some $(t_*,x_*)\in(0,\infty)\times\dbR^n$.
Define
$$
\widetilde w_i(t,x):=e^{\Lambda(t)}w_i(t,x),
\qquad i=1,2.
$$
Since $\Lambda$ is bounded,
$\widetilde w_1$ and $\widetilde w_2$
retain uniform linear growth.
Moreover, $\widetilde w_2$ is uniformly
Lipschitz continuous in space, with
Lipschitz constant denoted by $L_*$.
A direct calculation shows that
$\widetilde w_1$ and $\widetilde w_2$
are, respectively, a viscosity subsolution
and supersolution of the transformed
HJBI equation
\begin{equation}\label{eq:jump-transformed-HJBI}
\begin{cases}
\begin{aligned}
-\partial_t\widetilde w(t,x)
-\sup_{u\in U}\inf_{v\in V}\Big\{
&\mathfrak L^{u,v}\widetilde w(t,x)
+\mathfrak B^{u,v}\widetilde w(t,x)\\
&+\widetilde f\big(
t,x,\widetilde w(t,x),
\sigma(t,x,u,v)^\top D\widetilde w(t,x),
\mathfrak C^{u,v}\widetilde w(t,x),
u,v\big)
\Big\}=0,
\end{aligned}\\[1ex]
\displaystyle
\lim_{t\to\infty}\widetilde w(t,x)
=\widetilde\Phi(x),
\quad\text{locally uniformly in }x\in\dbR^n,
\end{cases}
\end{equation}
where
$
\widetilde\Phi(x):=e^{\Lambda_\infty}\Phi(x)
$
and
$\ds
\Lambda_\infty:=\lim_{t\to\infty}\Lambda(t).
$
For $\widetilde w_1$ and $\widetilde w_2$,
the terminal condition is understood
in the corresponding one-sided sense:
$$
\limsup_{\substack{t\to\infty,y\to x}}
\widetilde w_1(t,y)
\les\widetilde\Phi(x)
\les\liminf_{\substack{t\to\infty,y\to x}}
\widetilde w_2(t,y),
\qquad x\in\dbR^n.
$$
Moreover, for $y_1\ges y_2$, the Lipschitz
property of $f$ in $y$ and $q\ges\ell_y$ yield
\begin{equation}\label{eq:jump-properness}
 \widetilde f(t,a,y_1,z,k,u,v)
-\widetilde f(t,a,y_2,z,k,u,v) \les
[\ell_y(t)-q(t)-e^{-t}](y_1-y_2)
\les-e^{-t}(y_1-y_2).
\end{equation}
Thus, $\widetilde f$ is strictly decreasing in $y$,
with a time-dependent monotonicity rate $e^{-t}>0$.

\no \emph{Step 1. Terminal-time localization.}
Choose $0<a_*<t_*$ and sufficiently small
$\alpha,\nu>0$ such that
$$
\widetilde w_1(t_*,x_*)-\widetilde w_2(t_*,x_*)
-\alpha\Gamma(t_*,x_*)
-\frac{\nu}{t_*-a_*}>0.
$$
We first claim that, for some sufficiently large
$T>t_*$,
\begin{equation}\label{eq:jump-terminal-localization}
{
\widetilde w_1(T,x)-\widetilde w_2(T,x)
\les\frac{\alpha}{4}\Gamma(T,x),
\qquad x\in\dbR^n.}
\end{equation}
Indeed, since $\widetilde w_1,\widetilde w_2$
have uniform linear growth and
$\Gamma(t,x)\ges c(1+|x|^2)$,
there exists $R_\alpha>0$ such that
$$
\widetilde w_1(t,x)-\widetilde w_2(t,x)
\les\frac{\alpha}{4}\Gamma(t,x),
\qquad |x|\ges R_\alpha,\quad t\ges0.
$$
By Remark \ref{Re-upp-loe-Con}, on the compact ball $\overline B_{R_\alpha}$, we have
$$
\lim_{t\to\infty}\sup_{|x|\les R_\alpha}
(\widetilde w_1(t,x)-\widetilde\Phi(x))^+=0,
\qquad
\lim_{t\to\infty}\sup_{|x|\les R_\alpha}
(\widetilde\Phi(x)-\widetilde w_2(t,x))^+=0.
$$
Then
$
\ds
\limsup_{t\to\infty}\sup_{|x|\les R_\alpha}
[\widetilde w_1(t,x)-\widetilde w_2(t,x)]\les0.
$
Hence, for all sufficiently large $t$,
$$
\widetilde w_1(t,x)-\widetilde w_2(t,x)
\les\frac{\alpha c}{4}
\les\frac{\alpha}{4}\Gamma(t,x),
\qquad |x|\les R_\alpha.
$$
Combining these estimates proves
\eqref{eq:jump-terminal-localization}.
Fix such a $T$.

\noindent
\emph{Step 2. Doubling of variables and localization of
the maximum.}
Following the doubling-of-variables method
in \cite{CIL1992,BarlesImbert2008},
for any nonnegative continuous function
$A$ on $[0,T]$, define
$$
\Xi_\varepsilon^A(t,x,y)
:= \widetilde w_1(t,x)-\widetilde w_2(t,y)
-\frac{|x-y|^2}{2\varepsilon} -\alpha\Gamma(t,x)-A(t)
-\frac{\nu}{t-a_*},
$$
for $(t,x,y)\in(a_*,T]\times\dbR^n\times\dbR^n$.

We claim the following properties,
where {\rm(i)--(ii)} hold for any nonnegative
$A\in C([0,T])$ and $0<\varepsilon\les1$,
while {\rm(iii)--(v)} concern
$\Xi_\varepsilon:=\Xi_\varepsilon^{A_\varepsilon}$
defined below, for sufficiently small
$\varepsilon>0$:
\bel{The-Pro}
\left\{\2n
\begin{aligned}
{\rm(i)}\quad&
\Xi_\varepsilon^A(t,x,y)>0
\ \Longrightarrow\ |x|,|y|\les R,\\[1ex]
{\rm(ii)}\quad&
\text{for any positive global maximum point }
(t_\varepsilon,x_\varepsilon,y_\varepsilon)
\text{ of }\Xi_\varepsilon^A,\
|x_\varepsilon-y_\varepsilon|\les K\varepsilon,\\[1ex]
{\rm(iii)}\quad&
\Xi_\varepsilon(T,x,y)<0,
\qquad x,y\in\dbR^n,\\[1ex]
{\rm(iv)}\quad&
\Xi_\varepsilon(t,x,y)<0,
\qquad a_*<t\les a_*+\theta,\quad
x,y\in\dbR^n,\\[1ex]
{\rm(v)}\quad&
\text{there exists a global maximum point }
(t_\varepsilon,x_\varepsilon,y_\varepsilon)
\text{ of }\Xi_\varepsilon\text{ such that}\\
&
m_\varepsilon
:=\Xi_\varepsilon(t_\varepsilon,x_\varepsilon,y_\varepsilon)
\ges m_*>0,\
a_*<t_\varepsilon<T,\
|x_\varepsilon|,|y_\varepsilon|\les R,\
|x_\varepsilon-y_\varepsilon|\les K\varepsilon.
\end{aligned}
\right.
\ee
Here $R,K,\theta,m_*>0$ are independent
of $\varepsilon$, and $R,K$ are also
independent of $A$.
We verify these properties below.

\ss
  (i).
By the uniform linear growth of
$\widetilde w_1,\widetilde w_2$,
the lower bound on $\Gamma$, and $A\ges0$,
we have
$$
\begin{aligned}
\Xi_\varepsilon^A(t,x,y)
&\les C(1+|x|)
-\alpha c(1+|x|^2)
+C|x-y|
-\frac{|x-y|^2}{2\varepsilon}\\
&\les C(1+|x|)
-\alpha c(1+|x|^2)+C\varepsilon.
\end{aligned}
$$
Thus, whenever
$\Xi_\varepsilon^A(t,x,y)>0$,
we have $|x|\les R_1$ for some
$R_1>0$ independent of $A$ and
$\varepsilon\in (0,1]$.
Moreover, under the same condition,
$$
0<\Xi_\varepsilon^A(t,x,y)
\les C(1+R_1)+C|x-y|
-\frac{|x-y|^2}{2},
$$
which implies $|x-y|\les R_2$
for some $R_2>0$ independent of
$A$ and $\varepsilon$.
Consequently, there exists $R>0$, independent of
$A$ and $\varepsilon$,
such that
$$
\Xi_\varepsilon^A(t,x,y)>0
\quad\Longrightarrow\quad
|x|,|y|\les R.
$$

 (ii)
At any positive global maximum
$(t_\varepsilon,x_\varepsilon,y_\varepsilon)$
of $\Xi_\varepsilon^A$, comparison with
$(t_\varepsilon,x_\varepsilon,x_\varepsilon)$
and the spatial Lipschitz continuity of
$\widetilde w_2$ give
$$
\frac{|x_\varepsilon-y_\varepsilon|^2}
{2\varepsilon}
\les
\widetilde w_2(t_\varepsilon,x_\varepsilon)
-\widetilde w_2(t_\varepsilon,y_\varepsilon)\\
 \les L_*|x_\varepsilon-y_\varepsilon|.
$$
Hence,
\begin{equation}\label{eq:jump-one-sided-distance}
|x_\varepsilon-y_\varepsilon|
\les K\varepsilon,
\end{equation}
where $K:=\max\{1,2L_*\}$.
In particular, $R$ and $K$ are independent
of the choice of $A$.

\medskip

Before proving {\rm(iii)--(v)}, we construct
the auxiliary function $A_\varepsilon$.
Fix $C_0>0$, independently of $\varepsilon$,
sufficiently large to control the
$D_\varepsilon$ terms in the estimates of
$I_{1,u,v}^{\varepsilon,\alpha}$ and
$I_{3,u,v}^{\varepsilon,\alpha,\delta}$
below.

On $[0,T]\times\mathcal D_\varepsilon$,
the quantities
$\widetilde w_2(t,y)$ and
$\sigma(t,y,u,v)^\top(x-y)/\varepsilon$
are uniformly bounded.
Moreover, the spatial Lipschitz continuity
of $\widetilde w_2$ gives
$$
|\mathfrak C^{u,v}\widetilde w_2(t,y)|
\les L_*|\rho(t,\cdot)|_{\lambda,2}
|\gamma(t,\cdot,y,u,v)|_{\lambda,2}.
$$
The corresponding split weighted jump argument
$\mathfrak C_\delta^{u,v}
(\widetilde w_2,\varphi_2)
(t_\varepsilon,y_\varepsilon)$
is also uniformly bounded for sufficiently
small $\delta$, chosen for each fixed
$\varepsilon$, as verified in Step~4 below,
following \eqref{eq:jump-weighted-difference}.
 Here $\varphi_2$ is  the lower test function
for $\widetilde w_2$ introduced in Step 3.
Consequently, we may fix $B>0$,
independently of $\varepsilon$, $(u,v)$,
and sufficiently small split radii,
such that the three reference arguments $(y,z,k)$
in Lemma~\ref{lem:jump-finite-scale-envelope}
satisfy
 $
|y |+|z |
+|k |\les B.
$
Applying Lemma~\ref{lem:jump-finite-scale-envelope}
with the fixed constants $R,K,B,C_0$,
we obtain a nonnegative function
$e_\varepsilon\in C([0,T])$ satisfying
$
\int_0^T e_\varepsilon(r)\,\mathrm dr
\longrightarrow0$,
$\varepsilon\downarrow0$.

Now set
$
A_\varepsilon(t)
:=\int_t^T e_\varepsilon(r)\,\mathrm dr,
$ and  $
\Xi_\varepsilon:=\Xi_\varepsilon^{A_\varepsilon}.
$
Consequently,
$$
 A_\varepsilon\in C^1([0,T]),\qquad
A_\varepsilon'(t)=-e_\varepsilon(t),\q\|A_\varepsilon\|_\infty
=\int_0^T e_\varepsilon(r)\,\mathrm dr
\longrightarrow0
\quad\text{as }\varepsilon\downarrow0.
$$
By the initial choice of $\alpha,\nu$,
$
\Xi_\varepsilon(t_*,x_*,x_*)>0
$
for all sufficiently small $\varepsilon$.

(iii)
{{On the other hand,
\eqref{eq:jump-terminal-localization} and the spatial Lipschitz
continuity of $\widetilde w_2$ imply}}
$$
{
\Xi_\varepsilon(T,x,y)
 \les-\frac{3\alpha}{4}\Gamma(T,x)
+L_*|x-y|
-\frac{|x-y|^2}{2\varepsilon} \les-\frac{3\alpha c}{4}
+\frac{L_*^2\varepsilon}{2}<0}
$$
for sufficiently small $\varepsilon$.

\ss

 (iv)
By (i), every positive value of
$\Xi_\varepsilon$ must occur in
$|x|,|y|\les R$.
On this region, the uniform linear growth
of $\widetilde w_1,\widetilde w_2$
and the nonnegativity of the remaining
penalty terms imply
$$
\Xi_\varepsilon(t,x,y)
\les C_R-\frac{\nu}{t-a_*}.
$$
Since the right-hand side tends to
$-\infty$ as $t\downarrow a_*$,
there exists $\theta>0$, independent
of $\varepsilon$, such that
$$
\Xi_\varepsilon(t,x,y)<0,
\qquad a_*<t\les a_*+\theta,
\quad x,y\in\dbR^n.
$$
We may choose $\theta<t_*-a_*$.

 \ss

{\rm(v)}
By the initial choice of $\alpha,\nu$,
$$
d_*:=
\widetilde w_1(t_*,x_*)
-\widetilde w_2(t_*,x_*)
-\alpha\Gamma(t_*,x_*)
-\frac{\nu}{t_*-a_*}>0.
$$
Since $\|A_\varepsilon\|_\infty\to0$,
for all sufficiently small $\varepsilon$,
$$
\Xi_\varepsilon(t_*,x_*,x_*)
=d_*-A_\varepsilon(t_*)
\ges\frac{d_*}{2}>0.
$$
By {\rm(i)}, {\rm(iii)}, and {\rm(iv)},
all positive values of $\Xi_\varepsilon$
are confined to the compact set
$
[a_*+\theta,T]\times
\overline B_R\times\overline B_R.
$
Since $\Xi_\varepsilon$ is continuous
and takes a positive value at
$(t_*,x_*,x_*)$, it attains a positive
global maximum $m_\varepsilon$ at some
$(t_\varepsilon,x_\varepsilon,y_\varepsilon)$.
By {\rm(iii)}, this maximum cannot occur
at $t=T$. Moreover, {\rm(ii)} gives
$$
\begin{gathered}
m_\varepsilon
=\Xi_\varepsilon(t_\varepsilon,
x_\varepsilon,y_\varepsilon)
\ges m_*:=\frac{d_*}{2}>0,\\
a_*<t_\varepsilon<T,\qquad
|x_\varepsilon|,|y_\varepsilon|\les R,\qquad
|x_\varepsilon-y_\varepsilon|
\les K\varepsilon.
\end{gathered}
$$
Hence, \eqref{The-Pro} is proved.

\medskip
\noindent
\emph{Step 3. Application of the parabolic nonlocal
Jensen--Ishii lemma.}

At the positive global maximum
$(t_\varepsilon,x_\varepsilon,y_\varepsilon)$
obtained in {\rm(v)}, set
$$
P_\varepsilon(t,x,y)
 :=\frac{|x-y|^2}{2\varepsilon}
+\alpha\Gamma(t,x)+A_\varepsilon(t)
+\frac{\nu}{t-a_*},\q
p_\varepsilon
 :=\frac{x_\varepsilon-y_\varepsilon}{\varepsilon}.
$$
{
If necessary, we first smooth $P_\varepsilon$
locally in time to meet the $C^2$ requirement
of the parabolic nonlocal theorem of sums.
The first time derivatives converge uniformly
on the localized cylinder, so the smoothing
error is absorbed into the auxiliary remainder.
}
By the parabolic theorem of sums
\cite[Theorem~8.3]{CIL1992}
and its nonlocal extension
\cite[Corollary~2]{BarlesImbert2008},
there exist $a_\varepsilon^1,a_\varepsilon^2
\in\dbR$ and symmetric matrices
$X_\varepsilon,Y_\varepsilon\in\mathbb S^n$
such that
$$
 \big(
a_\varepsilon^1,
p_\varepsilon+\alpha D\Gamma(t_\varepsilon,x_\varepsilon),
X_\varepsilon
\big)
\in\overline{\mathcal P}^{2,+}
\widetilde w_1(t_\varepsilon,x_\varepsilon),\q \big(
a_\varepsilon^2,p_\varepsilon,Y_\varepsilon
\big)
\in\overline{\mathcal P}^{2,-}
\widetilde w_2(t_\varepsilon,y_\varepsilon),
$$
with
$$
a_\varepsilon^1-a_\varepsilon^2
 =\alpha\partial_t\Gamma(t_\varepsilon,x_\varepsilon)
+A_\varepsilon'(t_\varepsilon)
-\frac{\nu}{(t_\varepsilon-a_*)^2} =\alpha\partial_t\Gamma(t_\varepsilon,x_\varepsilon)
-e_\varepsilon(t_\varepsilon)
-\frac{\nu}{(t_\varepsilon-a_*)^2}.
$$
Moreover, by choosing the auxiliary
matrix tolerance sufficiently small,
$$
\begin{pmatrix}
X_\varepsilon&0\\
0&-Y_\varepsilon
\end{pmatrix}
\les
\frac1\varepsilon
\begin{pmatrix}
I&-I\\
-I&I
\end{pmatrix} +\alpha
\begin{pmatrix}
D^2\Gamma(t_\varepsilon,x_\varepsilon)&0\\
0&0
\end{pmatrix}
+{o_\vartheta(1)}I_{2n},
$$
{
where $o_\vartheta(1)\to0$
as $\vartheta\downarrow0$ for each fixed
$\varepsilon>0$.
}

{
For the nonlocal terms, we use the mixed
semijet--test-function argument in
\cite[Definition~4 and Corollary~2]{BarlesImbert2008}.
The spatial sections of the original doubled
penalty at $t_\varepsilon$ are, up to
additive constants,
$$
\varphi_1(t_\varepsilon,z)
:=\frac{|z-y_\varepsilon|^2}{2\varepsilon}
+\alpha\Gamma(t_\varepsilon,z),\qquad
\varphi_2(t_\varepsilon,z)
:=-\frac{|x_\varepsilon-z|^2}{2\varepsilon}.
$$
We use these sections in the split nonlocal
operators below. The localized smooth-test
approximants converge to them in $C^2$ on
compact spatial sets for fixed $\varepsilon$
and $\delta$; their contributions to the
viscosity inequalities differ by
$o_\vartheta(1)$ as the auxiliary
regularization parameter
$\vartheta\downarrow0$.
The local second-order terms are still
evaluated through the semijet matrices
$X_\varepsilon,Y_\varepsilon$.
}
Their spatial gradients at the contact
points satisfy
$$
D\varphi_1(t_\varepsilon,x_\varepsilon)
 =p_\varepsilon+
\alpha D\Gamma(t_\varepsilon,x_\varepsilon),\q
D\varphi_2(t_\varepsilon,y_\varepsilon)
 =p_\varepsilon.
$$
The test functions can be localized
to preserve the class $C_{l,b}^3$,
while their time slopes are supplied
by the parabolic theorem of sums.

Applying the split viscosity inequalities
in Definition~\ref{def:viscosity-Isaacs-equation},
we obtain, for all sufficiently small
$\delta>0$,
{up to the auxiliary errors
$o_\vartheta(1)$ described above,}
\begin{equation}\label{eq:jump-doubling-sub}
\begin{aligned}
&0\les a_\varepsilon^1
+\sup_{u\in U}\inf_{v\in V}\Big\{
b(t_\varepsilon,x_\varepsilon,u,v)
\cdot D\varphi_1(t_\varepsilon,x_\varepsilon) +\frac12\operatorname{tr}\big[
\sigma\sigma^\top(t_\varepsilon,x_\varepsilon,u,v)
X_\varepsilon\big]\\
&+\mathfrak B_\delta^{u,v}
(\widetilde w_1,\varphi_1)
(t_\varepsilon,x_\varepsilon)
 +\widetilde f\big(
t_\varepsilon,x_\varepsilon,
\widetilde w_1(t_\varepsilon,x_\varepsilon),
\sigma(t_\varepsilon,x_\varepsilon,u,v)^\top
D\varphi_1(t_\varepsilon,x_\varepsilon),
\mathfrak C_\delta^{u,v}
(\widetilde w_1,\varphi_1)
(t_\varepsilon,x_\varepsilon),u,v
\big)
\Big\}{+o_\vartheta(1)},
\end{aligned}
\end{equation}
and
\begin{equation}\label{eq:jump-doubling-super}
\begin{aligned}
0\ges{}&a_\varepsilon^2
+\sup_{u\in U}\inf_{v\in V}\Big\{
b(t_\varepsilon,y_\varepsilon,u,v)
\cdot D\varphi_2(t_\varepsilon,y_\varepsilon) +\frac12\operatorname{tr}\big[
\sigma\sigma^\top(t_\varepsilon,y_\varepsilon,u,v)
Y_\varepsilon\big] +\mathfrak B_\delta^{u,v}
(\widetilde w_2,\varphi_2)
(t_\varepsilon,y_\varepsilon)\\
&+\widetilde f\big(
t_\varepsilon,y_\varepsilon,
\widetilde w_2(t_\varepsilon,y_\varepsilon),
\sigma(t_\varepsilon,y_\varepsilon,u,v)^\top
D\varphi_2(t_\varepsilon,y_\varepsilon),
\mathfrak C_\delta^{u,v}
(\widetilde w_2,\varphi_2)
(t_\varepsilon,y_\varepsilon),u,v
\big)
\Big\}{-o_\vartheta(1)},
\end{aligned}
\end{equation}
Here the split nonlocal operators are
defined in \eqref{eq:Bdelta}.
For each fixed $\varepsilon$ and $\delta$,
the $C^2$ convergence of the regularized
spatial sections and {\bf(G6)} control the
small-jump integrals uniformly over
$U\times V$. For the large-jump integrals,
the common compact jump-landing region,
$\lambda(E_\delta^c)<\infty$, and
$\rho(t,\cdot)\in L_\lambda^2(E)$
justify the same passage to the limit.
The resulting auxiliary errors are included
in $o_\vartheta(1)$; the standard
regularization and limiting-semijet
arguments are omitted.

 \emph{Step 4. Estimates of the local and nonlocal terms.}
For each $(u,v)\in U\times V$, define
\begin{equation}\label{eq:jump-four-differences}
\begin{aligned}
I_{1,u,v}^{\varepsilon,\alpha}
:={}&\frac12\operatorname{tr}\big[
\sigma\sigma^\top(t_\varepsilon,x_\varepsilon,u,v)
X_\varepsilon -
\sigma\sigma^\top(t_\varepsilon,y_\varepsilon,u,v)
Y_\varepsilon\big],\\
I_{2,u,v}^{\varepsilon,\alpha}
:={}&
b(t_\varepsilon,x_\varepsilon,u,v)
\cdot D\varphi_1(t_\varepsilon,x_\varepsilon) -
b(t_\varepsilon,y_\varepsilon,u,v)
\cdot D\varphi_2(t_\varepsilon,y_\varepsilon),\\
I_{3,u,v}^{\varepsilon,\alpha,\delta}
:={}&
\mathfrak B_\delta^{u,v}
(\widetilde w_1,\varphi_1)
(t_\varepsilon,x_\varepsilon) -
\mathfrak B_\delta^{u,v}
(\widetilde w_2,\varphi_2)
(t_\varepsilon,y_\varepsilon),\\
I_{4,u,v}^{\varepsilon,\alpha,\delta}
:={}&
\widetilde f\big(
t_\varepsilon,x_\varepsilon,
\widetilde w_1(t_\varepsilon,x_\varepsilon),
\sigma(t_\varepsilon,x_\varepsilon,u,v)^\top
D\varphi_1(t_\varepsilon,x_\varepsilon),
\mathfrak C_\delta^{u,v}
(\widetilde w_1,\varphi_1)
(t_\varepsilon,x_\varepsilon),u,v
\big)\\
&-
\widetilde f\big(
t_\varepsilon,y_\varepsilon,
\widetilde w_2(t_\varepsilon,y_\varepsilon),
\sigma(t_\varepsilon,y_\varepsilon,u,v)^\top
D\varphi_2(t_\varepsilon,y_\varepsilon),
\mathfrak C_\delta^{u,v}
(\widetilde w_2,\varphi_2)
(t_\varepsilon,y_\varepsilon),u,v
\big).
\end{aligned}
\end{equation}
Subtracting \eqref{eq:jump-doubling-super}
from \eqref{eq:jump-doubling-sub}, together with
the elementary minimax inequality,
we obtain
\begin{equation}\label{eq:jump-four-sum}
0\les
a_\varepsilon^1-a_\varepsilon^2 +\sup_{(u,v)\in U\times V}
\Big\{
I_{1,u,v}^{\varepsilon,\alpha}
+I_{2,u,v}^{\varepsilon,\alpha} +
I_{3,u,v}^{\varepsilon,\alpha,\delta}
+I_{4,u,v}^{\varepsilon,\alpha,\delta}
\Big\}{+o_\vartheta(1)}.
\end{equation}

We estimate the four terms separately.
By the matrix inequality from the
parabolic theorem of sums,
\begin{equation}\label{eq:jump-I1}
\begin{aligned}
I_{1,u,v}^{\varepsilon,\alpha}
\les{}&
\frac{\alpha}{2}\operatorname{tr}\big[
\sigma\sigma^\top(t_\varepsilon,x_\varepsilon,u,v)
D^2\Gamma(t_\varepsilon,x_\varepsilon)
\big] +\frac{C}{\varepsilon}
\big|
\sigma(t_\varepsilon,x_\varepsilon,u,v)
-\sigma(t_\varepsilon,y_\varepsilon,u,v)
\big|^2+o_\vartheta(1)\\
\les{}&
\frac{\alpha}{2}\operatorname{tr}\big[
\sigma\sigma^\top(t_\varepsilon,x_\varepsilon,u,v)
D^2\Gamma(t_\varepsilon,x_\varepsilon)
\big]+C D_\varepsilon(t_\varepsilon)+o_\vartheta(1).
\end{aligned}
\end{equation}
Here $o_\vartheta(1)\to0$ as $\vartheta\downarrow0$
for fixed $\varepsilon,\alpha,\delta$, uniformly
in $(u,v)$, since $\sigma$ is bounded on the
localized compact cylinder.

 For $I_{2,u,v}^{\varepsilon,\alpha}$, recalling the expressions for
$D\varphi_1$ and $D\varphi_2$, we have
\begin{equation}\label{eq:jump-I2}
\begin{aligned}
I_{2,u,v}^{\varepsilon,\alpha}
={}&\alpha b(t_\varepsilon,x_\varepsilon,u,v)
\cdot D\Gamma(t_\varepsilon,x_\varepsilon) +\frac{x_\varepsilon-y_\varepsilon}{\varepsilon}
\cdot\big[
b(t_\varepsilon,x_\varepsilon,u,v)
-b(t_\varepsilon,y_\varepsilon,u,v)
\big]\\
\les{}&
\alpha b(t_\varepsilon,x_\varepsilon,u,v)
\cdot D\Gamma(t_\varepsilon,x_\varepsilon)
+K B_\varepsilon(t_\varepsilon).
\end{aligned}
\end{equation}

By the global maximum property of
$\Xi_\varepsilon$, for every $e\in E$,
$$
\begin{aligned}
&\widetilde w_1\big(
t_\varepsilon,
x_\varepsilon+
\gamma(t_\varepsilon,e,x_\varepsilon,u,v)
\big)
-\widetilde w_1(t_\varepsilon,x_\varepsilon)-
\widetilde w_2\big(
t_\varepsilon,
y_\varepsilon+
\gamma(t_\varepsilon,e,y_\varepsilon,u,v)
\big)
+\widetilde w_2(t_\varepsilon,y_\varepsilon)\\
&\les
\frac{x_\varepsilon-y_\varepsilon}{\varepsilon}
\cdot\Big[
\gamma(t_\varepsilon,e,x_\varepsilon,u,v)
-\gamma(t_\varepsilon,e,y_\varepsilon,u,v)
\big] +
\frac1{2\varepsilon}
\big|
\gamma(t_\varepsilon,e,x_\varepsilon,u,v)
-\gamma(t_\varepsilon,e,y_\varepsilon,u,v)
\big|^2\\
&\quad+
\alpha\big[
\Gamma\big(
t_\varepsilon,
x_\varepsilon+
\gamma(t_\varepsilon,e,x_\varepsilon,u,v)
\big)
-\Gamma(t_\varepsilon,x_\varepsilon)
\big].
\end{aligned}
$$
Then in  $I_{3,u,v}^{\varepsilon,\alpha,\delta}$,
 integrating the large-jump difference on
$E_\delta^c$ is bounded above by
$$\begin{aligned}
&
 \alpha\int_{E_\delta^c}\big[
\Gamma(t_\varepsilon,x_\varepsilon+\gamma(t_\varepsilon,e,x_\varepsilon,u,v))-\Gamma(t_\varepsilon,x_\varepsilon)
-D\Gamma(t_\varepsilon,x_\varepsilon)\cdot\gamma(t_\varepsilon,e,x_\varepsilon,u,v)\big]\lambda(\mathrm de) \\
&+\frac1{2\varepsilon}
\int_{E_\delta^c}|\gamma(t_\varepsilon,e,x_\varepsilon,u,v)-\gamma(t_\varepsilon,e,y_\varepsilon,u,v)|^2
\lambda(\mathrm de).
\end{aligned}
$$

By {\bf(G6)}, since
$t_\varepsilon\in[0,T]$ and
$|x_\varepsilon|,|y_\varepsilon|\les R$,
for $0<\delta\les1$,
$$
 |\gamma(t_\varepsilon,e,x_\varepsilon,u,v)|
+|\gamma(t_\varepsilon,e,y_\varepsilon,u,v)| \les C_{T,R}|e|,
\qquad \lambda\text{-a.e. }e\in E_\delta.
$$
{
For the small jumps, the spatial sections
of the original doubled penalty specified
in Step~3 give the following exact identity:
}
$$
\begin{aligned}
&\int_{E_\delta}\Big[
\varphi_1\big(
t_\varepsilon,
x_\varepsilon+
\gamma(t_\varepsilon,e,x_\varepsilon,u,v)
\big)
-\varphi_1(t_\varepsilon,x_\varepsilon)
-D\varphi_1(t_\varepsilon,x_\varepsilon)
\cdot\gamma(t_\varepsilon,e,x_\varepsilon,u,v)
\Big]\lambda(\mathrm de)\\
&\quad-\int_{E_\delta}\Big[
\varphi_2\big(
t_\varepsilon,
y_\varepsilon+
\gamma(t_\varepsilon,e,y_\varepsilon,u,v)
\big)
-\varphi_2(t_\varepsilon,y_\varepsilon)
-D\varphi_2(t_\varepsilon,y_\varepsilon)
\cdot\gamma(t_\varepsilon,e,y_\varepsilon,u,v)
\Big]\lambda(\mathrm de)\\
& =
\alpha\int_{E_\delta}\Big[
\Gamma\big(t_\varepsilon,
x_\varepsilon+\gamma(t_\varepsilon,e,x_\varepsilon,u,v)\big)
-\Gamma(t_\varepsilon,x_\varepsilon)
-D\Gamma(t_\varepsilon,x_\varepsilon)
\cdot\gamma(t_\varepsilon,e,x_\varepsilon,u,v)
\Big]\lambda(\mathrm de)\\
&\qquad+\frac1{2\varepsilon}
\int_{E_\delta}\Big[
|\gamma(t_\varepsilon,e,x_\varepsilon,u,v)
-\gamma(t_\varepsilon,e,y_\varepsilon,u,v)|^2
+2\gamma(t_\varepsilon,e,x_\varepsilon,u,v)
\cdot\gamma(t_\varepsilon,e,y_\varepsilon,u,v)
\Big]\lambda(\mathrm de)\\
& \les
\alpha\int_{E_\delta}\Big[
\Gamma\big(
t_\varepsilon,
x_\varepsilon+
\gamma(t_\varepsilon,e,x_\varepsilon,u,v)
\big)
-\Gamma(t_\varepsilon,x_\varepsilon)
-D\Gamma(t_\varepsilon,x_\varepsilon)
\cdot\gamma(t_\varepsilon,e,x_\varepsilon,u,v)
\Big]\lambda(\mathrm de)\\
&\qquad+\frac1{2\varepsilon}
\int_{E_\delta}
\big|
\gamma(t_\varepsilon,e,x_\varepsilon,u,v)
-\gamma(t_\varepsilon,e,y_\varepsilon,u,v)
\big|^2
\lambda(\mathrm de)
+r_{\varepsilon,\delta}^{B},
\end{aligned}
$$
where the nonnegative remainder can be chosen
so that, for each fixed $\varepsilon>0$
and sufficiently small $\delta$,
$$
0\les r_{\varepsilon,\delta}^{B}
\les\frac{C_{T,R}}{\varepsilon}
\int_{E_\delta}|e|^2\lambda(\mathrm de)
\longrightarrow0
\qquad(\delta\downarrow0),
$$
uniformly over $(u,v)\in U\times V$.

Based on the above,
we obtain
\begin{equation}\label{eq:jump-I3}
I_{3,u,v}^{\varepsilon,\alpha,\delta}
\les
\alpha\mathfrak B^{u,v}\Gamma
(t_\varepsilon,x_\varepsilon)
+C D_\varepsilon(t_\varepsilon)
+r_{\varepsilon,\delta}^{B}.
\end{equation}

In order to apply Lemma  \ref{lem:jump-finite-scale-envelope} to $I_{4,u,v}^{\varepsilon,\alpha,\delta}$,
we estimate the weighted jump arguments.
By the global maximum inequality,
$\rho\ges0$, and the smooth-test
approximation discussed above, we obtain
\begin{equation}\label{eq:jump-weighted-difference}
\begin{aligned}
&\mathfrak C_\delta^{u,v}
(\widetilde w_1,\varphi_1)
(t_\varepsilon,x_\varepsilon)
-\mathfrak C_\delta^{u,v}
(\widetilde w_2,\varphi_2)
(t_\varepsilon,y_\varepsilon)\\
&\les
\alpha\mathfrak C^{u,v}\Gamma
(t_\varepsilon,x_\varepsilon)
+\frac{x_\varepsilon-y_\varepsilon}{\varepsilon}
\cdot\int_E\rho(t_\varepsilon,e)
\big[
\gamma(t_\varepsilon,e,x_\varepsilon,u,v)
-\gamma(t_\varepsilon,e,y_\varepsilon,u,v)
\big]\lambda(\mathrm de)\\
&\quad+\frac1{2\varepsilon}
\int_{E_\delta^c}\rho(t_\varepsilon,e)
\big|
\gamma(t_\varepsilon,e,x_\varepsilon,u,v)
-\gamma(t_\varepsilon,e,y_\varepsilon,u,v)
\big|^2\lambda(\mathrm de)\\
&\quad+\frac1{2\varepsilon}
\int_{E_\delta}\rho(t_\varepsilon,e)
\Big[
|\gamma(t_\varepsilon,e,x_\varepsilon,u,v)|^2
+|\gamma(t_\varepsilon,e,y_\varepsilon,u,v)|^2
\Big]\lambda(\mathrm de)\\
&=
\alpha\mathfrak C^{u,v}\Gamma
(t_\varepsilon,x_\varepsilon)
+\frac{x_\varepsilon-y_\varepsilon}{\varepsilon}
\cdot\int_E\rho(t_\varepsilon,e)
\big[
\gamma(t_\varepsilon,e,x_\varepsilon,u,v)
-\gamma(t_\varepsilon,e,y_\varepsilon,u,v)
\big]\lambda(\mathrm de)\\
&\quad+\frac1{2\varepsilon}
\int_E\rho(t_\varepsilon,e)
\big|
\gamma(t_\varepsilon,e,x_\varepsilon,u,v)
-\gamma(t_\varepsilon,e,y_\varepsilon,u,v)
\big|^2\lambda(\mathrm de)\\
&\quad+\frac1\varepsilon
\int_{E_\delta}\rho(t_\varepsilon,e)
\gamma(t_\varepsilon,e,x_\varepsilon,u,v)
\cdot\gamma(t_\varepsilon,e,y_\varepsilon,u,v)
\lambda(\mathrm de)\\
&\les
\alpha\mathfrak C^{u,v}\Gamma
(t_\varepsilon,x_\varepsilon)
+K_\varepsilon(t_\varepsilon)
+r_{\varepsilon,\delta}^{C},
\end{aligned}
\end{equation}
{
where
$$
\begin{aligned}
r_{\varepsilon,\delta}^{C}
&:=\frac{C_{T,R}\delta}{\varepsilon}
\sup_{0\les r\les T}
\int_{E_\delta}\rho(r,e)|e|
\lambda(\mathrm de)\\
&\les \frac{C_{T,R}\delta}{\varepsilon}
\sup_{0\les r\les T}|\rho(r,\cdot)|_{\lambda,2}
\Big(\int_{E_\delta}|e|^2
\lambda(\mathrm de)\Big)^{1/2}
\longrightarrow0
\qquad(\delta\downarrow0)
\end{aligned}
$$
for each fixed $\varepsilon>0$,
uniformly over $(u,v)\in U\times V$.
}

Moreover, by the spatial Lipschitz continuity
of $\widetilde w_2$ and {\bf(G6)},
$$
\big|
\mathfrak C^{u,v}\widetilde w_2
(t_\varepsilon,y_\varepsilon)
\big|
\les
L_*|\rho(t_\varepsilon,\cdot)|_{\lambda,2}
|\gamma(t_\varepsilon,\cdot,y_\varepsilon,u,v)|_{\lambda,2}
\les C_{T,R}.
$$
After removing the auxiliary regularization,
the quadratic structure of the doubled penalty
and \eqref{eq:jump-one-sided-distance} yield
$$
 \big|
\mathfrak C_\delta^{u,v}
(\widetilde w_2,\varphi_2)(t_\varepsilon,y_\varepsilon)
-\mathfrak C^{u,v}\widetilde w_2
(t_\varepsilon,y_\varepsilon)
\big| \les
C_{T,R}\(K+L_*+\frac{\delta}{\varepsilon}\)
\sup_{0\les r\les T}
\int_{E_\delta}\rho(r,e)|e|\lambda(\mathrm de)
\longrightarrow0
\qquad(\delta\downarrow0),
$$
for each fixed $\varepsilon>0$,
uniformly over $(u,v)\in U\times V$.
Consequently, choosing $\delta$ sufficiently
small for each fixed $\varepsilon$, the
reference weighted jump argument is uniformly
bounded, as required in
Lemma~\ref{lem:jump-finite-scale-envelope}.

Moreover, the two $z$-arguments satisfy
$$
\begin{aligned}
&\sigma(t_\varepsilon,x_\varepsilon,u,v)^\top
D\varphi_1(t_\varepsilon,x_\varepsilon)
-\sigma(t_\varepsilon,y_\varepsilon,u,v)^\top
D\varphi_2(t_\varepsilon,y_\varepsilon)\\
& =
\alpha\sigma(t_\varepsilon,x_\varepsilon,u,v)^\top
D\Gamma(t_\varepsilon,x_\varepsilon) +
\big[
\sigma(t_\varepsilon,x_\varepsilon,u,v)
-\sigma(t_\varepsilon,y_\varepsilon,u,v)
\big]^\top
\frac{x_\varepsilon-y_\varepsilon}{\varepsilon}.
\end{aligned}
$$
By \eqref{eq:jump-one-sided-distance},
$$
\Big|
\big[
\sigma(t_\varepsilon,x_\varepsilon,u,v)
-\sigma(t_\varepsilon,y_\varepsilon,u,v)
\big]^\top
\frac{x_\varepsilon-y_\varepsilon}{\varepsilon}
\Big|
\les K S_\varepsilon(t_\varepsilon).
$$

Since $m_\varepsilon>0$, we have
$
\widetilde w_1(t_\varepsilon,x_\varepsilon)
>\widetilde w_2(t_\varepsilon,y_\varepsilon).
$
Moreover, as verified above, the reference
arguments
$
\widetilde w_2(t_\varepsilon,y_\varepsilon),$ $
\sigma(t_\varepsilon,y_\varepsilon,u,v)^\top
D\varphi_2(t_\varepsilon,y_\varepsilon),$ $
\mathfrak C_\delta^{u,v}
(\widetilde w_2,\varphi_2)
(t_\varepsilon,y_\varepsilon)
$
satisfy the uniform bound $B$ required in
Lemma~\ref{lem:jump-finite-scale-envelope}
for sufficiently small $\delta$.
Together with \eqref{eq:jump-weighted-difference}
and the preceding estimate of the $z$-arguments,
this allows us to apply
Lemma~\ref{lem:jump-finite-scale-envelope}.
Using also the monotonicity of $\widetilde f$
in $k$, \eqref{eq:jump-properness}, and
its Lipschitz bounds in $z$ and $k$,
we obtain
\begin{equation}\label{eq:jump-I4}
\begin{aligned}
I_{4,u,v}^{\varepsilon,\alpha,\delta}
\les{}&
-{e^{-t_\varepsilon}}
\Big[
\widetilde w_1(t_\varepsilon,x_\varepsilon)
-\widetilde w_2(t_\varepsilon,y_\varepsilon)
\Big] +\alpha\ell_z(t_\varepsilon)
\big|
\sigma(t_\varepsilon,x_\varepsilon,u,v)^\top
D\Gamma(t_\varepsilon,x_\varepsilon)
\big|\\
&+\alpha\ell_k(t_\varepsilon)
\big(
\mathfrak C^{u,v}\Gamma
(t_\varepsilon,x_\varepsilon)
\big)^++e_\varepsilon(t_\varepsilon)
-K B_\varepsilon(t_\varepsilon)
-C_0 D_\varepsilon(t_\varepsilon)
+r_{\varepsilon,\delta}^{f}.
\end{aligned}
\end{equation}
Here
$$
r_{\varepsilon,\delta}^{f}
:=\ell_k(t_\varepsilon)r_{\varepsilon,\delta}^{C}
\longrightarrow0
\qquad(\delta\downarrow0)
$$
for each fixed $\varepsilon>0$,
uniformly over $(u,v)\in U\times V$.

Combining
\eqref{eq:jump-I1}--\eqref{eq:jump-I4},
and recalling that $C_0$ was chosen
in the construction of $e_\varepsilon$
to dominate the combined
$D_\varepsilon$-terms in
\eqref{eq:jump-I1} and \eqref{eq:jump-I3},
{with the auxiliary $o_\vartheta(1)$
errors absorbed into the remainder,}
we obtain
\begin{equation}\label{eq:jump-four-estimates}
\begin{aligned}
&\sup_{(u,v)\in U\times V}
\Big\{
I_{1,u,v}^{\varepsilon,\alpha}
+I_{2,u,v}^{\varepsilon,\alpha}
+I_{3,u,v}^{\varepsilon,\alpha,\delta}
+I_{4,u,v}^{\varepsilon,\alpha,\delta}
\Big\}\\
&\quad\les
\alpha\mathcal M[\Gamma]
(t_\varepsilon,x_\varepsilon) -{e^{-t_\varepsilon}}
\Big[
\widetilde w_1(t_\varepsilon,x_\varepsilon)
-\widetilde w_2(t_\varepsilon,y_\varepsilon)
\Big] +e_\varepsilon(t_\varepsilon)
+r_{\varepsilon,\delta},
\end{aligned}
\end{equation}

{
For each fixed $\varepsilon>0$, choose
$0<\delta=\delta(\varepsilon)\les\varepsilon$
sufficiently small so that the required
small-jump localization and reference bounds
hold, and
$$
r_{\varepsilon,\delta}^{B}
+\ell_k(t_\varepsilon)r_{\varepsilon,\delta}^{C}
\les\varepsilon.
$$
With this $\delta$ fixed, choose the
auxiliary smooth-test regularization and
matrix tolerance sufficiently small that
their combined $o_\vartheta(1)$ errors
in \eqref{eq:jump-four-sum} and
\eqref{eq:jump-I1} do not exceed
$\varepsilon$, uniformly over $U\times V$.
Including these errors in
$r_{\varepsilon,\delta}$, we obtain
$$
0\les r_{\varepsilon,\delta}\les2\varepsilon
\longrightarrow0
\qquad(\varepsilon\downarrow0),
$$
uniformly over $(u,v)\in U\times V$.
}

Combining \eqref{eq:jump-four-sum}
and \eqref{eq:jump-four-estimates}, we obtain
\begin{equation}\label{eq:jump-comparison-final}
\begin{aligned}
&0\les
a_\varepsilon^1-a_\varepsilon^2
+\alpha\mathcal M[\Gamma]
(t_\varepsilon,x_\varepsilon) -{e^{-t_\varepsilon}}
\big[
\widetilde w_1(t_\varepsilon,x_\varepsilon)
-\widetilde w_2(t_\varepsilon,y_\varepsilon)
\big] +e_\varepsilon(t_\varepsilon)
+r_{\varepsilon,\delta}\\
&=
\alpha\big[
\partial_t\Gamma+\mathcal M[\Gamma]
\big](t_\varepsilon,x_\varepsilon) -{e^{-t_\varepsilon}}
\big[
\widetilde w_1(t_\varepsilon,x_\varepsilon)
-\widetilde w_2(t_\varepsilon,y_\varepsilon)
\big] -\frac{\nu}{(t_\varepsilon-a_*)^2}
+r_{\varepsilon,\delta}.
\end{aligned}
\end{equation}
On the other hand, by the definition of
$m_\varepsilon$,
$$
 \widetilde w_1(t_\varepsilon,x_\varepsilon)
-\widetilde w_2(t_\varepsilon,y_\varepsilon) =
m_\varepsilon
+\frac{|x_\varepsilon-y_\varepsilon|^2}
{2\varepsilon}
+\alpha\Gamma(t_\varepsilon,x_\varepsilon) +
A_\varepsilon(t_\varepsilon)
+\frac{\nu}{t_\varepsilon-a_*} \ges
m_\varepsilon+\alpha\Gamma(t_\varepsilon,x_\varepsilon).
$$
Substituting this into
\eqref{eq:jump-comparison-final}, we get
$$
0\les
\alpha\big[
\partial_t\Gamma+\mathcal M[\Gamma]
-{e^{-t}\Gamma}
\big](t_\varepsilon,x_\varepsilon)
 -{e^{-t_\varepsilon}}m_\varepsilon
-\frac{\nu}{(t_\varepsilon-a_*)^2}
+r_{\varepsilon,\delta}.
$$
By the Lyapunov barrier inequality
\eqref{eq:jump-integrable-barrier},
\begin{equation}\label{eq:jump-comparison-contradiction}
0\les
-{e^{-t_\varepsilon}}m_\varepsilon
-\frac{\nu}{(t_\varepsilon-a_*)^2}
+r_{\varepsilon,\delta}.
\end{equation}

\emph{Step 5. Contradiction and conclusion.}
Recall that, by the initial choice of
$\alpha,\nu$, there exists $m_*>0$
such that
$$
m_\varepsilon\ges
\Xi_\varepsilon(t_*,x_*,x_*)
\ges m_*
$$
for all sufficiently small $\varepsilon$.
{Moreover, since $a_*<t_\varepsilon<T$,}
$$
{e^{-t_\varepsilon}\ges e^{-T}},
\qquad
\frac{\nu}{(t_\varepsilon-a_*)^2}
\ges\frac{\nu}{(T-a_*)^2}.
$$
Hence \eqref{eq:jump-comparison-contradiction}
implies
$$
0\les
-e^{-T}m_*
-\frac{\nu}{(T-a_*)^2}
+r_{\varepsilon,\delta}.
$$
Choosing the split radius and matrix
tolerance so that
$r_{\varepsilon,\delta}\to0$
as $\varepsilon\downarrow0$,
we obtain a contradiction.

Therefore,
$
w_1(t,x)\les w_2(t,x)
$
for all $t>0$ and $x\in\dbR^n$.
By continuity, the inequality also holds
at $t=0$.

\ms

If instead $\widetilde w_1$ is uniformly
Lipschitz continuous in space, with constant
$L_*$, {comparison of the positive
maximum with
$(t_\varepsilon,y_\varepsilon,y_\varepsilon)$
gives
$$
|x_\varepsilon-y_\varepsilon|
\les 2\varepsilon\Big(
L_*+\alpha
\sup_{[0,T]\times\overline B_R}|D\Gamma|
\Big).
$$
The constant $K$ is adjusted accordingly,
and the reference $z$- and $k$-arguments
in Lemma~\ref{lem:jump-finite-scale-envelope}
are taken from $\widetilde w_1$ at
$(t_\varepsilon,x_\varepsilon)$.
Their uniform bounds follow from the
spatial Lipschitz continuity of
$\widetilde w_1$ and {\bf(G6)}.
The remaining argument is analogous,
and the details are omitted.}
\end{proof}

\begin{theorem}
\label{thm:jump-viscosity-uniqueness}\sl
Assume {\bf(G1)--(G7)}.
Then the lower value function $W$ is the unique
continuous viscosity solution of
\eqref{equ:5.1}, in the class of functions
with uniform linear growth.
\end{theorem}

\begin{proof}
By Theorem~\ref{thm:jump-viscosity-existence},
$W$ is a viscosity solution of \eqref{equ:5.1}.
Moreover, $W$ is uniformly Lipschitz continuous
in the spatial variable.
Let $w$ be any other continuous viscosity
solution of \eqref{equ:5.1} with uniform
linear growth.
Applying Theorem~\ref{thm:jump-one-sided-comparison}
to $(w,W)$ and $(W,w)$ yields $w=W$.
\end{proof}

\subsection{Viscosity characterization of the upper value function}

We next consider the upper value function $U$.
The upper Isaacs Hamiltonian is defined by
\begin{equation}\label{eq:Hplus}
H^+[\varphi](t,x)
:=\inf_{v\in V}\sup_{u\in U}
H^{u,v}[\varphi](t,x).
\end{equation}
The corresponding HJBI equation is
\begin{equation}\label{equ:5.2}
\begin{cases}
-\partial_t U(t,x)-H^+[U](t,x)=0,
& (t,x)\in[0,\infty)\times\dbR^n,\\[1ex]
\displaystyle
\lim_{t\to\infty}U(t,x)=\Phi(x),
& \text{locally uniformly in }x\in\dbR^n.
\end{cases}
\end{equation}
Using the $\delta$-split nonlocal operators
in \eqref{eq:Bdelta}, define
\begin{equation}\label{eq:Hdelta-plus}
\begin{aligned}
&\mathcal H_\delta^+(t,x,y,\varphi;w)
:=\inf_{v\in V}\sup_{u\in U}
\Big\{
\mathfrak L^{u,v}\varphi(t,x)
+\mathfrak B_\delta^{u,v}(w,\varphi)(t,x)\\
&\hskip 3.5cm
+f\big(
t,x,y,
\sigma(t,x,u,v)^\top D\varphi(t,x),
\mathfrak C_\delta^{u,v}(w,\varphi)(t,x),
u,v
\big)
\Big\}.
\end{aligned}
\end{equation}

Viscosity sub- and supersolutions of
\eqref{equ:5.2} are defined as in
Definition~\ref{def:viscosity-Isaacs-equation},
with $\mathcal H_\delta^-$ replaced by
$\mathcal H_\delta^+$.

\begin{theorem}
\label{thm:jump-upper-viscosity-uniqueness}\sl
Assume {\bf(G1)--(G7)}.
Then the upper value function $U$ is the unique
continuous viscosity solution of
\eqref{equ:5.2}, in the class of functions
with uniform linear growth.
\end{theorem}

\begin{proof}
The viscosity characterization and comparison
principle follow by the same arguments as for
the lower HJBI equation, with
$\sup_u\inf_v$ replaced by $\inf_v\sup_u$.
Since $U$ is uniformly Lipschitz continuous
in the spatial variable, uniqueness follows.
\end{proof}

\begin{corollary}[A contact-wise Isaacs criterion for a game value]
\label{cor:jump-isaacs-value}\sl
Assume {\bf(G1)--(G7)}. Then $W\les U$.
Moreover, $W=U$ if, for some
$\Psi\in\{W,U\}$,
$$
\liminf_{\delta\downarrow0}
\Big[
\mathcal H_\delta^+
(t,x,\Psi(t,x),\varphi;\Psi)
-\mathcal H_\delta^-
(t,x,\Psi(t,x),\varphi;\Psi)
\Big]=0
$$
at every contact point $(t,x)$ and for
every admissible test function $\varphi$
touching $\Psi$ from below if $\Psi=W$,
or from above if $\Psi=U$.
\end{corollary}

\begin{proof}
{Since $\mathcal H_\delta^-\les\mathcal H_\delta^+$,
$W$ is a viscosity subsolution of the upper HJBI equation.
Theorem~\ref{thm:jump-upper-viscosity-uniqueness} and its
upper-Hamiltonian comparison argument give $W\les U$.}

 Fix an admissible test function $\varphi$
touching $\Psi$ at $(t,x)$   as in the statement.  By the local contact property
and the jump-size bound in {\bf(G6)}, there exists
$\delta_*>0$ such that, for every
$0<\delta_2<\delta_1<\delta_*$ and all controls,
$x+\gamma(t,e,x,u,v)$ remains in the contact
neighborhood for $\lambda$-a.e. $e\in E_{\delta_1}$.

Writing $\gamma=\gamma(t,e,x,u,v)$, we have
$$
\begin{aligned}
&\mathfrak B_{\delta_2}^{u,v}(\Psi,\varphi)(t,x)
-\mathfrak B_{\delta_1}^{u,v}(\Psi,\varphi)(t,x) =\int_{E_{\delta_1}\setminus E_{\delta_2}}
\big[(\Psi-\varphi)(t,x+\gamma)
-(\Psi-\varphi)(t,x)\big]\lambda(\mathrm de),\\
&\mathfrak C_{\delta_2}^{u,v}(\Psi,\varphi)(t,x)
-\mathfrak C_{\delta_1}^{u,v}(\Psi,\varphi)(t,x) =\int_{E_{\delta_1}\setminus E_{\delta_2}}
\rho(t,e)\big[(\Psi-\varphi)(t,x+\gamma)
-(\Psi-\varphi)(t,x)\big]\lambda(\mathrm de).
\end{aligned}
$$
At a contact from below these differences are
nonnegative, whereas at a contact from above
they are nonpositive. Since $\rho\ges0$ and
$f$ is nondecreasing in its jump argument,
both $\mathcal H_\delta^-$ and
$\mathcal H_\delta^+$ are nondecreasing as
$\delta\downarrow0$ for $\Psi=W$, and
nonincreasing for $\Psi=U$.

At this fixed contact, write
$$
h_\delta^\pm
:=\mathcal H_\delta^\pm
(t,x,\Psi(t,x),\varphi;\Psi),
\qquad
g_\delta:=h_\delta^+-h_\delta^-\ges0.
$$
By assumption, there exists
$\delta_n\downarrow0$ such that
$g_{\delta_n}\to0$.
Fix any sufficiently small $\delta_0>0$
and take $n$ large enough that
$\delta_n<\delta_0$.

If $\Psi=W$, the lower-Hamiltonian
supersolution inequality gives
$$
\partial_t\varphi(t,x)+h_{\delta_0}^+
 \les\partial_t\varphi(t,x)+h_{\delta_n}^+ =\partial_t\varphi(t,x)+h_{\delta_n}^-
+g_{\delta_n}
\les g_{\delta_n}.
$$
Letting $n\to\infty$ shows that $W$
is a viscosity supersolution of the
upper HJBI equation.

If $\Psi=U$, the upper-Hamiltonian
subsolution inequality similarly gives
$$
\partial_t\varphi(t,x)+h_{\delta_0}^-
 \ges\partial_t\varphi(t,x)+h_{\delta_n}^- =\partial_t\varphi(t,x)+h_{\delta_n}^+
-g_{\delta_n}
\ges-g_{\delta_n}.
$$
Thus $U$ is a viscosity subsolution
of the lower HJBI equation.

In either case, the corresponding
comparison principle, using the uniform
spatial Lipschitz continuity of $W$,
yields $U\les W$. Hence $W=U$.

\end{proof}

\begin{remark}\sl
The full Hamiltonians $H^-[\Psi]$ and
$H^+[\Psi]$ need not be classically defined
for $\Psi=W$ or $U$. The condition in
Corollary~\ref{cor:jump-isaacs-value}
therefore uses smooth test functions for
the derivatives and small jumps, while
retaining $\Psi$ in the large-jump terms.
It is weaker than requiring the full
split Isaacs identity
$\mathcal H_\delta^-=\mathcal H_\delta^+$
for all admissible continuation functions
and test functions.

Alternatively, if
$H^-[\varphi](t,x)=H^+[\varphi](t,x)$
for every admissible smooth test function
$\varphi$ and every $(t,x)$, then $W=U$
also follows from the equivalence of the
full-test and $\delta$-split formulations in
Remark~\ref{rem:viscosity-equivalence}
and the comparison principle.
\end{remark}

\appendix
\section{Proof of Lemma~\ref{lem:jump-short-freeze}}
\label{app:proof-jump-short-freeze}

Before the proof, we point out {that the} eighth-moment estimate below is local in time and follows from {\bf(G6)}, without requiring a global $p=8$ assumption.

\begin{proof}
Fix $M,R>0$ and set $L:=M+1$. Throughout the proof,
$C$ denotes a deterministic constant, independent of
$t\in[0,M]$, $|x|\les R$, $0<h\les1$, and the controls;
its value may change from line to line. Write
$X:=X^{t,x;u,v}$, and set  $\mathbb P_t(A):=\mathbb E_t[\mathbf1_A]$.

\medskip

\emph{Step 1. Uniform eighth-moment estimate.}
For $r\in[0,L]$, set
$$
a_X(r):={\bar q(r)+\bar q^0(r)},
\qquad
A_L:=\int_0^L a_X(r)\,\mathrm dr.
$$
By {\bf (G1)} and {\bf (G5)}, $A_L<\infty$.

We first derive the growth estimates needed
for the eighth moment.
By the state Lipschitz conditions,
for every $(r,y,u,v)$,
$$
\begin{aligned}
&|b(r,y,u,v)|
\les a_X(r)(1+|y|),\q
|\sigma(r,y,u,v)|^2
 \les 2a_X(r)(1+|y|^2),\\
& \int_E|\gamma(r,e,y,u,v)|^2
\lambda(\mathrm de)
 \les 2a_X(r)(1+|y|^2).
\end{aligned}
$$
Moreover, the pointwise growth condition
in {\bf (G6)} yields
$
|\gamma(r,e,y,u,v)|
\les C_L(1+|y|),
$ $ \lambda\text{-a.e. }e\in E.
$
Consequently,
$$
\begin{aligned}
&\int_E|\gamma(r,e,y,u,v)|^8
\lambda(\mathrm de) \les
C_L^6(1+|y|)^6
\int_E|\gamma(r,e,y,u,v)|^2
\lambda(\mathrm de)\\
& \les
C_La_X(r)(1+|y|)^6(1+|y|^2) \les
C_La_X(r)(1+|y|^8).
\end{aligned}
$$

Fix
$(u,v)\in\mathcal U_{t,t+h}
\times\mathcal V_{t,t+h}$.
For $N>R$, let
$
\tau_N:=\inf\{s\in[t,t+h]:|X_s|\ges N\},
$
with $\inf\varnothing:=\infty$, and set
$
X_s^N:=X_{s\wedge\tau_N},
$ $ s\in[t,t+h].
$
The stopped equation reads
$$
\begin{aligned}
X_s^N=x
 +\int_t^{s\wedge\tau_N}
b(r,X_r,u_r,v_r)\,\mathrm dr +\int_t^{s\wedge\tau_N}
\sigma(r,X_r,u_r,v_r)\,\mathrm dB_r +\int_t^{s\wedge\tau_N}\!\!\int_E
\gamma(r,e,X_{r-},u_r,v_r)
\,\widetilde\mu(\mathrm dr,\mathrm de).
\end{aligned}
$$
  Then,
we have
\begin{equation}\label{eq:freeze-moment-decomposition}
H_N(s):=\mathbb E_t\Big[
1+\sup_{t\les q\les s}|X_q^N|^8
\Big]
\les C(1+|x|^8)+C(I_b+I_\sigma+I_\gamma),
\end{equation}
where
$$
\begin{aligned}
I_b
&:=\mathbb E_t\Big[
\Big(\int_t^{s\wedge\tau_N}
|b(r,X_r,u_r,v_r)|\,\mathrm dr
\Big)^8\Big]\les
\mathbb E_t\Big[
\Big(\int_t^s
a_X(r)(1+|X_{r\wedge\tau_N}|)
\,\mathrm dr\Big)^8\Big] \\
&\les
A_L^7\int_t^s a_X(r)
\mathbb E_t\big[
(1+|X_r^N|)^8
\big]\,\mathrm dr \les
CA_L^7\int_t^s a_X(r)H_N(r)\,\mathrm dr,\\
I_\sigma
&:=\mathbb E_t\Big[
\sup_{t\les q\les s}
\Big|\int_t^{q\wedge\tau_N}
\sigma(r,X_r,u_r,v_r)\,\mathrm dB_r
\Big|^8\Big]\les   C\mathbb E_t\Big[
\Big(\int_t^s a_X(r)
(1+|X_r^N|^2)\,\mathrm dr\Big)^4
\Big]\\
 &\les CA_L^3\int_t^s a_X(r)
\mathbb E_t\big[
(1+|X_r^N|^2)^4
\big]\,\mathrm dr \les CA_L^3
\int_t^s a_X(r)H_N(r)\,\mathrm dr  ,\\
I_\gamma
&:=\mathbb E_t\Big[
\sup_{t\les q\les s}
\Big|\int_t^{q\wedge\tau_N}\int_E
\gamma(r,e,X_{r-},u_r,v_r)
\,\widetilde\mu(\mathrm dr,\mathrm de)
\Big|^8\Big]\\
&\les
C\mathbb E_t\Big[
\Big(\int_t^{s\wedge\tau_N}\int_E
|\gamma(r,e,X_{r-},u_r,v_r)|^2
\lambda(\mathrm de)\,\mathrm dr
\Big)^4+
\int_t^{s\wedge\tau_N}\int_E
|\gamma(r,e,X_{r-},u_r,v_r)|^8
\lambda(\mathrm de)\,\mathrm dr
\Big]\\
&\les C\mathbb E_t\Big[
\Big(\int_t^s a_X(r)
(1+|X_{r\wedge\tau_N}|^2)
\,\mathrm dr\Big)^4\Big] +C_L\mathbb E_t\Big[
\int_t^{s\wedge\tau_N}
a_X(r)(1+|X_{r-}|^8)
\,\mathrm dr\Big]\\
&\les C_L(1+A_L^3)
\int_t^s a_X(r)H_N(r)\,\mathrm dr .
\end{aligned}
$$
Therefore,
we obtain
$$
H_N(s)
\les C_L(1+|x|^8)
+C_L(1+A_L^7+A_L^3)
\int_t^s a_X(r)H_N(r)\,\mathrm dr.
$$
Since $a_X\in L^1(0,L)$,
Gronwall's inequality gives
$$
\begin{aligned}
H_N(s)
&\les C_L(1+|x|^8)
\exp\Big(
C_L(1+A_L^7+A_L^3)
\int_t^s a_X(r)\,\mathrm dr
\Big)\\
&\les C_L(1+|x|^8)
\exp\big(
C_L(1+A_L^7+A_L^3)A_L
\big) \les C_{M,R}.
\end{aligned}
$$
The constant $C_{M,R}$ is independent of
$N$, $t$, $h$, and $(u,v)$.

Note that the original state equation is
nonexplosive (Lemma \ref{lemma:game-state-estimates}), so
$\tau_N\uparrow\infty$ almost surely
as $N\to\infty$.
By the conditional Fatou lemma,
$$ \mathbb E_t\Big[
1+\sup_{t\les r\les t+h}|X_r|^8
\Big]\ \les
\liminf_{N\to\infty}
\mathbb E_t\Big[
1+\sup_{t\les r\les t+h}|X_r^N|^8
\Big]
\les C_{M,R}.
$$
Taking the supremum over all admissible
controls and the prescribed initial data,
we conclude that
%
\begin{equation}\label{eq:freeze-eighth-moment}
\sup_{\substack{
0\les t\les M,\ |x|\les R,\ 0<h\les1\\
(u,v)\in\mathcal U_{t,t+h}\times\mathcal V_{t,t+h}}}
\left\|
\mathbb E_t\Big[
1+\sup_{t\les r\les t+h}
|X_r^{t,x;u,v}|^8
\Big]
\right\|_\infty
\les C_{M,R}.
\end{equation}

\emph{Step 2. A high-order exit-probability estimate.}
By the joint continuity of $b,\sigma$ in
{\bf (G1)} and the $L_\lambda^2$-continuity
of $\gamma$ in {\bf (G6)}, for every $N>0$,
\begin{equation}\label{eq:freeze-local-bounded}
\sup_{\substack{
0\les r\les L,\ |y|\les N\\
(u,v)\in U\times V}}
\Big[
|b(r,y,u,v)|+|\sigma(r,y,u,v)|
+|\gamma(r,\cdot,y,u,v)|_{\lambda,2}
\Big]\les C_{L,N}.
\end{equation}

Let $C_L>0$ be the constant in the
pointwise jump bound of {\bf (G6)} on $[0,L]$.
Choose deterministic radii
$$
N_0:=R+2,\qquad
N_{j+1}:=N_j+C_L(1+N_j)+2,
\quad j=0,1,2,
$$
and define the stopping times
$\ds
\tau_j:=\inf\{r\in[t,t+h]:|X_r|\ges N_j\},
$ $ j=0,1,2,3,
$
with the convention $\inf\varnothing=\infty$.
For $N>0$, let
$
B_N:=\{x\in\dbR^n:|x|<N\}
$.

For $j=0,1,2$, on $\{\tau_j\les t+h\}$,
the definition of $\tau_j$ and the
c\`adl\`ag property of $X$ imply
$|X_{\tau_j-}|\les N_j$.
Since the jumps of $X$ are generated by
the Poisson random measure,   {\bf (G6)} yields
$
|\Delta X_{\tau_j}|
\les C_L(1+|X_{\tau_j-}|),
$ $ \mathbb P\text{-a.s.}
$
Consequently,
\begin{equation}\label{eq:freeze-no-overshoot}
|X_{\tau_j}|
 \les |X_{\tau_j-}|+|\Delta X_{\tau_j}| \les N_j+C_L(1+N_j)
=N_{j+1}-2.
\end{equation}
%

To estimate the exit probabilities, set
$\tau_{-1}:=t$ and $X_{\tau_{-1}}:=x$.
Since $N_0-|x|\ges2$, it follows from
\eqref{eq:freeze-no-overshoot} that,
for $j=-1,0,1,2$,
$$
N_{j+1}-|X_{\tau_j}|\ges2
\quad\text{on }\{\tau_j\les t+h\}.
$$
By \eqref{eq:freeze-local-bounded}
and the conditional BDG inequalities,
on $\{\tau_j\les t+h\}$,
$$
\mathbb E\Big[
\sup_{\tau_j\les r\les(t+h)\wedge\tau_{j+1}}
|X_r-X_{\tau_j}|^2
\,\big|\,\mathcal F_{\tau_j}\Big]
\les Ch.
$$
Since $\{\tau_{j+1}\les t+h\}$ implies
$
|X_{\tau_{j+1}}-X_{\tau_j}|\ges2,
$
the conditional Markov inequality and
the tower property yield
$$
\mathbb P_t(\tau_{j+1}\les t+h)
 =\mathbb E_t\Big[
\mathbf1_{\{\tau_j\les t+h\}}
\mathbb P(\tau_{j+1}\les t+h
\mid\mathcal F_{\tau_j})\Big] \les Ch\,\mathbb P_t(\tau_j\les t+h),
\quad j=-1,0,1,2.
$$
As $\mathbb P_t(\tau_{-1}\les t+h)=1$,
iteration gives
\begin{equation}\label{eq:freeze-rare-exit}
\mathbb P_t(\tau_3\les t+h)\les Ch^4.
\end{equation}

Finally, applying the conditional
$L^2$ maximal inequalities directly
to the state equation stopped at $\tau_3$,
and using \eqref{eq:freeze-local-bounded},
we obtain
\begin{equation}\label{eq:freeze-stopped-increment}
\mathbb E_t\Big[
\sup_{t\les r\les t+h}
|X_{r\wedge\tau_3}-x|^2
\Big]\les Ch.
\end{equation}
All constants depend only on $L$, $R$,
and the local coefficient bounds, and are
independent of $t$, $h$, and the predictable
controls $(u,v)$.

\medskip
\noindent
\emph{Step 3. Freezing the generator at zero
BSDE variables.}
Set
$$
\Delta F_r^0
:= F^\varphi(r,X_r,0,0,0,u_r,v_r) -F^\varphi(r,x,0,0,0,u_r,v_r).
$$
Let $N:=N_3$ and
$A_h:=\{\tau_3\les t+h\}$.
By the continuity of $F^\varphi$ and
compactness, there exist a constant
$C_{L,N,\varphi}>0$ and a modulus of continuity
$\omega(\eta)\to0$ as $\eta\downarrow0$ such that
$$
\begin{aligned}
&|F^\varphi(r,a,0,0,0,u,v)
-F^\varphi(r,x,0,0,0,u,v)|
\les\omega(\eta)
\end{aligned}
$$
whenever $r\in[0,L]$, $|x|\les R$,
$|a|\les N$, $|a-x|\les\eta$, and
$(u,v)\in U\times V$.

On $A_h^c$, the state remains in $B_N$
throughout $[t,t+h]$.
By the conditional Markov inequality
and \eqref{eq:freeze-stopped-increment},
for every $\eta>0$,
\begin{equation}\label{eq:freeze-local-g}
\begin{aligned}
&\mathbb E_t\Big[
\Big(\mathbf1_{A_h^c}
\int_t^{t+h}|\Delta F_r^0|
\,\mathrm dr\Big)^2
\Big] \les h^2\Big[
\omega(\eta)
+C\mathbb P_t\Big(
\sup_{t\les r\les t+h}
|X_{r\wedge\tau_3}-x|>\eta
\Big)^{{\frac{1}{2}}}
\Big]^2\\
& \les Ch^2\(
\omega(\eta)^2+\frac{h}{\eta^2}
\).
\end{aligned}
\end{equation}
On $A_h$, the growth estimate
\eqref{eq:freeze-g-growth}, conditional
H\"older's inequality,
\eqref{eq:freeze-eighth-moment},
and \eqref{eq:freeze-rare-exit} yield
\begin{equation}\label{eq:freeze-tail-g}
\begin{aligned}
&\mathbb E_t\Big[
\Big(\mathbf1_{A_h}
\int_t^{t+h}|\Delta F_r^0|
\,\mathrm dr\Big)^2
\Big]\les
C_R\(\int_t^{t+h}
a_\varphi(r)\,\mathrm dr\)^2
\mathbb E_t\Big[
\Big(1+\sup_{t\les r\les t+h}|X_r|^2
\Big)^4\Big]^{{\frac{1}{2}}}
\mathbb P_t(A_h)^{{\frac{1}{2}}}\\
&\les
Ch^2\(\int_t^{t+h}
a_\varphi(r)\,\mathrm dr\)^2.
\end{aligned}
\end{equation}
Since $a_\varphi\in L^1(0,L)$,
$$
A_\varphi(h):=
\sup_{0\les t\les M}
\int_t^{t+h}a_\varphi(r)\,\mathrm dr
\longrightarrow0
\quad\text{as }h\downarrow0.
$$

Finally, taking $\eta=h^{{\frac{1}{4}}}$ in
\eqref{eq:freeze-local-g} and combining
it with \eqref{eq:freeze-tail-g}, we obtain
\begin{equation}\label{eq:freeze-g-oh}
 \mathbb E_t\Big[
\Big(\int_t^{t+h}|\Delta F_r^0|
\,\mathrm dr\Big)^2
\Big]
\les h^2\varepsilon_0(h)^2,
\end{equation}
where $\varepsilon_0(h)
:=C\big[
\omega(h^{{\frac{1}{4}}})
+h^{{\frac{1}{4}}}+A_\varphi(h)
\big]\longrightarrow0,$
 as  $h\downarrow0$.
The estimate is uniform in
$t\in[0,M]$, $|x|\les R$, and all
admissible controls $(u,v)$.

\medskip
\noindent
\emph{Step 4. Estimates for the nonzero BSDE variables.}
By the joint continuity of
$(r,a,u,v)\mapsto F^\varphi(r,a,0,0,0,u,v)$
and the compactness of
$[0,L]\times\overline B_R\times U\times V$,
we have
$$
\sup_{\substack{
0\les r\les L,\ |a|\les R\\
(u,v)\in U\times V}}
|F^\varphi(r,a,0,0,0,u,v)|
\les C_{L,R,\varphi}.
$$
Together with the time-integrable
Lipschitz conditions in {\bf (G2)},
the conditional $L^2$ estimate
for the frozen BSDE \eqref{eq:jump-Y2},
whose terminal value is zero, gives
\begin{equation}\label{eq:freeze-Y2-small}
\begin{aligned}
&
\mathbb E_t\Big[
 \sup_{t\les r\les t+h}|Y_r^2|^2
+\int_t^{t+h}|Z_r^2|^2\,\mathrm dr +\int_t^{t+h}|K_r^2|_{\lambda,2}^2
\,\mathrm dr
\Big]\\
&\les
C\mathbb E_t\Big[
\Big(\int_t^{t+h}
|F^\varphi(r,x,0,0,0,u_r,v_r)|
\,\mathrm dr\Big)^2
\Big]\les C_{L,R,\varphi}h^2,
\end{aligned}\end{equation}
where the constant in the BSDE estimate
depends only on the time integrals of
the Lipschitz coefficients, while
$C_{L,R,\varphi}$ also depends on the
local uniform bound of $F^\varphi$.

Next, set
$
\Delta F_r:=
F^\varphi(r,X_r,Y_r^2,Z_r^2,K_r^2,u_r,v_r)
 -F^\varphi(r,x,Y_r^2,Z_r^2,K_r^2,u_r,v_r).
$
By \eqref{eq:jump-F-difference}, we have
$$
|\Delta F_r|
\les
|\Delta F_r^0|
+2\ell_y(r)|Y_r^2|
+2\ell_z(r)|Z_r^2| +2\ell_k(r)|\rho(r,\cdot)|_{\lambda,2}
|K_r^2|_{\lambda,2}.
$$

For $0<h\les1$, define
$$
\theta_M(h)
:=\sup_{0\les t\les M}\Big\{ \int_t^{t+h}\ell_y(r)\,\mathrm dr +\Big(\int_t^{t+h}\ell_z(r)^2
\,\mathrm dr\Big)^{{\frac{1}{2}}} +\Big(\int_t^{t+h}
\ell_k(r)^2|\rho(r,\cdot)|_{\lambda,2}^2
\,\mathrm dr\Big)^{{\frac{1}{2}}}
\Big\}.
$$
By {\bf (G2)} and H\"older's inequality,
all the integrands are in $L^1(0,L)$.
Hence, the absolute continuity of the integral
implies
$
\theta_M(h)\to0
$
as $h\downarrow0$.

{Using conditional Cauchy--Schwarz and
\eqref{eq:freeze-Y2-small} for the \emph{squares} of the
integrals, we obtain the following $L^2$ estimates.}
\begingroup
$$\begin{aligned}
&\dbE_t\Big[\Big(\int_t^{t+h}
\ell_y(r)|Y_r^2|\,\mathrm dr\Big)^2\Big]
\les \Big(\int_t^{t+h}\ell_y(r)\,\mathrm dr\Big)^2
\dbE_t\Big[\sup_{t\les r\les t+h}|Y_r^2|^2\Big]
\les Ch^2\Big(\int_t^{t+h}\ell_y(r)\,\mathrm dr\Big)^2,\\
&\dbE_t\Big[\Big(\int_t^{t+h}
\ell_z(r)|Z_r^2|\,\mathrm dr\Big)^2\Big]
\les \Big(\int_t^{t+h}\ell_z(r)^2\,\mathrm dr\Big)
\dbE_t\Big[\int_t^{t+h}|Z_r^2|^2\,\mathrm dr\Big]
\les Ch^2\int_t^{t+h}\ell_z(r)^2\,\mathrm dr,\\
&\dbE_t\Big[\Big(\int_t^{t+h}
\ell_k(r)|\rho(r,\cdot)|_{\lambda,2}
|K_r^2|_{\lambda,2}\,\mathrm dr\Big)^2\Big]
\les Ch^2\int_t^{t+h}
\ell_k(r)^2|\rho(r,\cdot)|_{\lambda,2}^2\,\mathrm dr.
\end{aligned}$$
\endgroup
Combining these estimates with
\eqref{eq:freeze-g-oh}, we conclude that
\begin{equation}\label{eq:freeze-F-oh}
\dbE_t\[\Big |
\int_t^{t+h}|\Delta F_r|\,\mathrm dr\Big|^2
\]
 \les
  {h^2\bigl[\varepsilon_0(h)+C\theta_M(h)\bigr]^2
 =:h^2\varepsilon_F(h)^2},
\end{equation}
where
$
\varepsilon_F(h)\to0
$
as $h\downarrow0$, uniformly in
$t\in[0,M]$, $|x|\les R$, and the
admissible controls $(u,v)$.

Since the two short-time BSDEs have the same
zero terminal value, the conditional $L^2$
stability estimate for BSDEs with jumps,
together with \eqref{eq:freeze-F-oh}, yields
$$
|Y_t^{1;u,v}-Y_t^{2;u,v}|^2
 \les C\dbE_t\[\Big |
\int_t^{t+h}|\Delta F_r|\,\mathrm dr\Big|^2
\]  \les Ch^2\varepsilon_F(h)^2.
$$
Thus \eqref{eq:jump-short-freeze} follows by setting
$
\varepsilon_{M,R,\varphi}(h):=C\varepsilon_F(h),
$
which is deterministic and tends to zero
as $h\downarrow0$, uniformly in $t,x,u,v$.

\end{proof}


\begin{thebibliography}{99}

\bibitem{Applebaum2009}
D.~Applebaum,
\emph{L\'evy Processes and Stochastic Calculus},
2nd ed., Cambridge University Press, Cambridge, 2009.

\bibitem{BarlesBuckdahnPardoux1997}
G.~Barles, R.~Buckdahn and E.~Pardoux,
Backward stochastic differential equations and integral-partial
differential equations,
\emph{Stochastics Stochastics Rep.}
\textbf{60} (1997), 57--83.

\bibitem{BarlesImbert2008}
G.~Barles and C.~Imbert,
Second-order elliptic integro-differential equations:
Viscosity solutions' theory revisited,
\emph{Ann. Inst. H. Poincar\'e C Anal. Non Lin\'eaire}
\textbf{25} (2008), 567--585.

\bibitem{Biswas2012}
I.~H.~Biswas,
On zero-sum stochastic differential games with jump-diffusion
driven state: A viscosity solution framework,
\emph{SIAM J. Control Optim.}
\textbf{50} (2012), 1823--1858.

\bibitem{BriandHu1998}
P.~Briand and Y.~Hu,
Stability of BSDEs with random terminal time
and homogenization of semilinear elliptic PDEs,
\emph{J. Funct. Anal.}
\textbf{155} (1998), 455--494.

\bibitem{BuckdahnHuLi2011}
R.~Buckdahn, Y.~Hu and J.~Li,
Stochastic representation for solutions of Isaacs' type
integral-partial differential equations,
\emph{Stochastic Process. Appl.}
\textbf{121} (2011), 2715--2750.

\bibitem{BL2008}
R.~Buckdahn and J.~Li,
Stochastic differential games and viscosity solutions of
Hamilton--Jacobi--Bellman--Isaacs equations,
\emph{SIAM J. Control Optim.}
\textbf{47} (2008), 444--475.

\bibitem{BuckdahnLiZhao2021}
R.~Buckdahn, J.~Li and N.~Zhao,
Representation of limit values for
nonexpansive stochastic differential games,
\emph{J. Differential Equations}
\textbf{276} (2021), 187--227.

\bibitem{ChenWang2000}
Z.~Chen and B.~Wang,
Infinite time interval BSDEs and the convergence
of $g$-martingales,
\emph{J. Austral. Math. Soc. Ser. A}
\textbf{69} (2000), 187--211.

\bibitem{CIL1992}
M.~G.~Crandall, H.~Ishii and P.-L.~Lions,
User's guide to viscosity solutions of second order partial
differential equations,
\emph{Bull. Amer. Math. Soc. (N.S.)}
\textbf{27} (1992), 1--67.

\bibitem{DarlingPardoux1997}
R.~W.~R.~Darling and E.~Pardoux,
Backwards SDE with random terminal time and applications
to semilinear elliptic PDE,
\emph{Ann. Probab.}
\textbf{25} (1997), 1135--1159.



\bibitem{HuTessitore2007}
Y.~Hu and G.~Tessitore,
BSDE on an infinite horizon and elliptic PDEs
in infinite dimension,
\emph{Nonlinear Differ. Equ. Appl.}
\textbf{14} (2007), 825--846.

\bibitem{HW2026}
S.~Huang and Q.~Wei,
Infinite-horizon non-autonomous zero-sum stochastic
recursive differential games and HJBI equations,
arXiv:2605.12926, 2026.

\bibitem{JacodShiryaev2003}
J.~Jacod and A.~N.~Shiryaev,
\emph{Limit Theorems for Stochastic Processes},
2nd ed., Springer, Berlin, 2003.

\bibitem{LiPeng2009}
J.~Li and S.~Peng,
Stochastic optimization theory of backward stochastic
differential equations with jumps and viscosity solutions
of Hamilton--Jacobi--Bellman equations,
\emph{Nonlinear Anal.}
\textbf{70} (2009), 1776--1796.

\bibitem{LiWei2014}
J.~Li and Q.~Wei,
$L^p$ estimates for fully coupled FBSDEs with jumps,
\emph{Stochastic Process. Appl.}
\textbf{124} (2014), 1582--1611.

\bibitem{LiWei2015}
J.~Li and Q.~Wei,
Stochastic differential games for fully coupled FBSDEs with jumps,
\emph{Appl. Math. Optim.}
\textbf{71} (2015), 411--448.

\bibitem{LiYu2025}
J.~Li and Z.~Yu,
Well-posedness of infinite horizon FBSDEs with non-zero
terminals and LQ problems with random coefficients,
\emph{Appl. Math. Optim.}
\textbf{91} (2025), Paper No.~64, 32 pp.

\bibitem{LiZhao2019}
J.~Li and N.~Zhao,
Representation of asymptotic values for
nonexpansive stochastic control systems,
\emph{Stochastic Process. Appl.}
\textbf{129} (2019), 634--673.

\bibitem{LiuLiFan2020}
Y.~Liu, D.~Li and S.~Fan,
$L^p$ ($p>1$) solutions of BSDEs with
generators satisfying some non-uniform
conditions in $t$ and $\omega$,
\emph{Chinese Ann. Math. Ser. B}
\textbf{41} (2020), 479--494.

\bibitem{LiuWei2026}
L.~Liu and Q.~Wei,
Reflected stochastic recursive control problems with jumps:
Dynamic programming and stochastic verification theorems,
\emph{ESAIM Control Optim. Calc. Var.}
\textbf{32} (2026), Article No.~36.

\bibitem{LuoLiLiWei2025}
S.~Luo, W.~Li, X.~Li and Q.~Wei,
Stochastic representation for solutions of a system of
coupled HJB-Isaacs equations with integral--differential
operators,
\emph{Stochastic Process. Appl.}
\textbf{179} (2025), Article No.~104502, 25 pp.

\bibitem{LLW2025}
S.~Luo, X.~Li and Q.~Wei,
Infinite time horizon stochastic recursive control problems with jumps:
Dynamic programming and stochastic verification theorems,
\emph{SIAM J. Control Optim.}
\textbf{63} (2025), 796--821.

\bibitem{LuoLiWei2026}
S.~Luo, X.~Li and Q.~Wei,
$L^p$-estimate for infinite horizon decoupled
forward--backward stochastic differential equations with jumps,
\emph{Appl. Math. Optim.}
\textbf{94} (2026), Paper No.~19, 23 pp.



\bibitem{Pardoux1998}
E.~Pardoux,
Backward stochastic differential equations and viscosity
solutions of systems of semilinear parabolic and elliptic
PDEs of second order,
in \emph{Stochastic Analysis and Related Topics VI},
Progr. Probab., Vol.~42,
Birkh\"auser, Boston, 1998, pp.~79--127.

\bibitem{PardouxPeng1990}
E.~Pardoux and S.~Peng,
Adapted solution of a backward stochastic differential equation,
\emph{Systems Control Lett.}
\textbf{14} (1990), 55--61.

\bibitem{Peng1991}
S.~Peng,
Probabilistic interpretation for systems of
quasilinear parabolic partial differential equations,
\emph{Stochastics Stochastics Rep.}
\textbf{37} (1991), 61--74.

\bibitem{Peng1997}
S.~Peng,
Backward stochastic differential equations---stochastic optimization
theory and viscosity solutions of HJB equations,
in \emph{Topics on Stochastic Analysis},
J.~Yan, S.~Peng, S.~Fang and L.~Wu, eds.,
Science Press, Beijing, 1997, pp.~85--138 (in Chinese).

\bibitem{PengShi2000}
S.~Peng and Y.~Shi,
Infinite horizon forward--backward stochastic differential equations,
\emph{Stochastic Process. Appl.}
\textbf{85} (2000), 75--92.

\bibitem{Royer2006}
M.~Royer,
Backward stochastic differential equations with jumps and related
non-linear expectations,
\emph{Stochastic Process. Appl.}
\textbf{116} (2006), 1358--1376.


\bibitem{Situ2000}
R.~Situ,
\emph{Backward Stochastic Differential Equations with Jumps
and Applications},
Guangdong Science \& Technology Press, Guangzhou, 2000.

\bibitem{Situ2005}
R.~Situ,
\emph{Theory of Stochastic Differential Equations with Jumps
and Applications},
Springer, 2005.

\bibitem{TangLi1994}
S.~Tang and X.~Li,
Necessary conditions for optimal control of stochastic systems
with random jumps,
\emph{SIAM J. Control Optim.}
\textbf{32} (1994), 1447--1475.

\bibitem{WeiYu2021}
Q.~Wei and Z.~Yu,
Infinite horizon forward-backward SDEs and open-loop
optimal controls for stochastic linear-quadratic
problems with random coefficients,
\emph{SIAM J. Control Optim.}
\textbf{59} (2021), 2594--2623.

\bibitem{Wu2003}
Z.~Wu,
Fully coupled FBSDE with Brownian motion and
Poisson process in stopping time duration,
\emph{J. Aust. Math. Soc.}
\textbf{74} (2003), 249--266.



\bibitem{Yu2017}
Z.~Yu,
Infinite horizon jump-diffusion forward-backward
stochastic differential equations and their
application to backward linear-quadratic problems,
\emph{ESAIM Control Optim. Calc. Var.}
\textbf{23} (2017), 1331--1359.

\end{thebibliography}
\end{document}